\documentclass[a4paper,11pt]{article}
\usepackage{amsmath}
\usepackage{dsfont}
\usepackage[usenames,dvipsnames]{pstricks}
\usepackage{epsfig}
\usepackage{pst-grad} 
\usepackage{pst-plot} 
\usepackage{amssymb}
 \title{A central limit theorem for random walk in a random environment on marked Galton-Watson trees}
 \author{Gabriel Faraud\footnote{Weierstraß-Institut
Mohrenstr. 39,
10117 Berlin, Germany.
Research partially supported by the ANR project MEMEMO. 
{\it faraud@math.univ-paris13.fr}}}

 \newtheorem{theorem}{Theorem}[section]
\newtheorem{proposition}{Proposition}[section]

\newtheorem{corollary}[theorem]{Corollary}
\newtheorem{lemma}[theorem]{Lemma}
\newtheorem{statement}{Statement}[section]
 \begin{document}
 \begin{titlepage}\maketitle
\noindent Abstract : Models of random walks in a random environment were introduced at first by Chernoff in 1967 in order to study biological mechanisms. The original model has been intensively studied since then and is now well understood. In parallel, similar models of random processes in a random environment have been studied. In this article we focus on a model of random walk on random marked trees, following a model introduced by R. Lyons and R. Pemantle (1992).
Our point of view is a bit different yet, as we consider a very general way of constructing random trees with random transition probabilities on them. We prove an analogue of R. Lyons and R. Pemantle's recurrence criterion in this setting, and we study precisely the asymptotic behavior, under restrictive assumptions. Our last result is a generalization of a result of Y. Peres and O. Zeitouni (2006) concerning biased random walks on Galton-Watson trees.\\



\noindent Submitted 09.12.2009, Accepted for publication 22.12.2010.
\end{titlepage}
\section{Introduction and statement of results.}

Models of random walks in a random environment were introduced at first by Chernov in 1967 (\cite{chernov1967rmc}) in order to study biological mechanisms. The original model has been intensively studied since then and is now well understood. On the other hand, more recently, several attempts have been made to study extensions of this original model, for example in higher dimensions, continuous time, or different space. 

It is remarkable that the random walk in $\mathbb{Z}^d$, $d>1$, is still quite
mysterious, in particular no precise criterion for recurrence/transience has ever been found. 

In the case of trees, however, a recurrence criterion exists, and even estimates for the asymptotic behavior have been proven. To present our model and the existing results, we begin with some notations concerning trees. 
Let $T$ be a tree rooted at some vertex $e$. For each vertex $x$ of $T$ we call $N({x})$ the number of his children $\{x_{1},x_{2},...,x_{N(x)}\}$, and $\overleftarrow{x}$ his father. For two vertices $x,y\in T$, we call $d(x,y)$ the distance between $x$ and $y$, that is the number of edges on the shortest path from $x$ to $y$, and $|x|:=d(e,x)$. Let $T_{n}$ be the set of vertices such that $|x|=n$, and $T^*=T\setminus\{e\}$. We also note $x<y$ when $x$ is an ancestor of $y$.

We call a marked tree a couple $(T,A),$ where $A$ is a random application from the vertices of $T$ to $\mathbb{R}_{+}^*$.
Let $\mathbb{T}$ be the set of marked trees.
We introduce the filtration ${\cal G}_{n}$ on $\mathbb{T}$ defined as
$${\cal G}_n=\sigma\{N(x),A(x_{i}),1\leq i \leq n, |x|<n,x\in T \}.$$ 
Following \cite{Neveu:1986}, given a probability measure $q$ on
$\mathbb{N}\otimes \mathbb{R}_{+}^{*\mathbb{N^*}}$, there exists a
probability measure {\tt MT} on $\mathbb{T}$ such that
\begin{itemize}
\item the distribution of the random variable $(N(e),A(e_{1}),A(e_{2}),...)$ is $q$,
\item given ${\cal G}_n$, the random variables $(N(x),A(x_{1}),A(x_{2}),.....),$ for $x\in T_n$, are independent and their conditional distribution is $q$.
\end{itemize}
We will always assume $m:=E[N(e)]>1,$ ensuring that the tree is infinite with a positive probability.

We now introduce the model of random walk in a random environment.
Given a marked tree $T$, we set for $x\in T^*$, $x_i$ a child of $x$,
$$\omega(x,x_i)=\frac{A(x_i)}{1+\sum_{j=1}^{N(x)} A(x_j)}$$
and $$\omega(x\overleftarrow{x})=\frac{1}{1+\sum_{j=1}^{N(x)} A(x_j)}.$$
Morever we set $\omega(x,y)=0$ whenever $d(x,y)\neq1$,

It is easy to check that $(\omega(x,y))_{ x,y\in
T}$ is a
family of non-negative random variables such that, $$\forall x\in T,\, \sum_{y\in
 T}\omega(x,y)=1,$$ and \begin{equation}\forall x\in T^*,\, A(x)=\frac{\omega(
 \overleftarrow{x},x)}{\omega(\overleftarrow{x},\overleftarrow{\overleftarrow{x}})},\end{equation}
where $\omega(e,\overleftarrow{e})$ is artificially defined as 
$$\frac{1}{\omega(e,\overleftarrow{e})}=\sum_{|x|=1} A(x).$$
Further, $\omega(x,y)\neq0$ whenever $x$ and
$y$ are neighbors. 

 $T$ will be called ``the environment'', and we call ``random walk on $T$'' the Markov chain $(X_{n},\mathbb{P}_T)$ defined by $X_0=e$ and
$$\forall x,y\in T ,\; 
 \mathbb{P}_T(X_{n+1}=y|X_{n}=x)=\omega(x,y)
.$$ 
We call ``annealed probability'' the probability $\mathbb{P}_{\tt MT}={\tt MT}\otimes \mathbb{P}_T$ taking into account the total alea.

We set, for $x\in T$, $C_x=\prod_{e< z \leq x} A(z).$
We can associate to the random walk $X_n$ an electrical network with conductance $C_{x}$ along $[\overleftarrow{x},x]$, and a capacited network with capacity $C_{x}$ along $[\overleftarrow{x},x]$.
We recall the definition of an electrical current on an electrical network. Let $G=(V,E)$ be a graph, $C$ be a symmetric function on $E$, and $A,Z$ be two disjoint subsets of $V$. We define the electrical current between $A$ and $Z$ as a function $i$ that is antisymmetric on $E$ such that, for any $x\in V\backslash(A\cup Z),$ the sum on the edges $e$ starting from $x$ of $i(e)$ equals zero (this is call Kirchhoff's node Law), and, moreover, $i$ satisfies the Kirchhoff's cycle Law, that is, for any cycle $x_{1}, x_{2},\dots,x_{n}=x_{1},$
$$\sum_{i=1}^n \frac{i(x_{i},x_{i+1})}{C(x_{i},x_{i+1})}=0.$$ 
A flow on a capacited network is an antisymmetric function $\theta$ that satisfies the Kirchhoff's node Law, and such that, for all edges $e$, $\theta(e)<C(e)$, (for more precisions on this correspondence we refer to the chapters 2 and 3 of \cite{Lyons:2005}).

We shall also frequently use the convex function $\rho$ defined for $\alpha\geq0$ as
\begin{equation*}
\rho(\alpha)= E_{\tt MT}\left[\sum_{0}^{N(e)} A(e_{i})^\alpha\right]=E_{q}\left[\sum_{0}^{N} A(i)^\alpha\right].
\end{equation*}
\textbf{Remark :} This model is in fact inspired by a model introduced in \cite{Lyons:1992}. In this case the tree $T$ and the $A(x)$ were introduced separately, and the $A(x)$ were supposed to be independent. Here we can include models in which the structure of the tree and the transition probabilities are dependent. A simple example that is covered in our model is the following : Let $T$ be a Galton-Watson tree. We chose an i.i.d. family $(B(x))_{x\in T}$ and set, for every $x\in T$, $1\le i\le N(x),$ $A(x_{i})=B(x)$. This way the transition probabilities to the children of any vertex are all equal, but randomly chosen. In R. Lyons and R. Pemantle's article, a recurrence criterion was shown, our first result is a version of this criterion in our setting.

\begin{theorem} \label{rec}
We suppose that there exists $0\leq\alpha\leq1$ such that $\rho$ is finite in a small neighborhood of $\alpha$, 
$\rho(\alpha)=\inf_{0\leq t\leq1} \rho(t):=p$ and $\rho'(\alpha)=E_q\left[\sum_{i=1}^{N(e)}A(e_i)^{\alpha}\log(A(e_i))\right]$ is finite.
We assume that $\sum_{i=1}^{N(e)} A(e_{i})$ is not identically equal to $1$.

 Then,
\begin{enumerate}
\item if $p<1$ then the RWRE is a.s. positive recurrent, the electrical network has zero conductance a.s., and the capacited network admits no flow a.s..
\item if $p\leq 1$ then the RWRE is a.s. recurrent, the electrical network has zero conductance a.s. and the capacited network admits no flow a.s..
\item if $p>1$, then, given non-extinction, the RWRE is a.s. transient, the electrical network has positive conductance a.s. and the capacited network admits flow a.s..
\end{enumerate}
 (By ``almost surely'' we mean ``for {\tt MT} almost every $T$'').
\end{theorem}
{\bf Remark:} In the case where $\sum_{i=1}^{N(e)} A(e_{i})$ is identically equal to $1$, which belongs to the second case, $|X_{n}|$ is a standard unbiased random walk, therefore $X_{n}$ is null recurrent. However, there exists a flow, given by $\theta(\overleftarrow{x},x)=C_{x}.$

The proof of this result is quite similar to the proof of R. Lyons and R. Pemantle, but there are some differences, coming from the fact that in their setting $i.i.d.$ random variables appear along any ray of the tree, whereas it is not the case here. Results on branching processes will help us address this problem.

Theorem \ref{rec} does not give a full answer in the case $p=1$, but this result can be improved, provided some technical assumptions are fulfilled. We introduce the condition 
$$(H1):\:\forall \alpha\in [0,1], \; E_q\left[\left(\sum_{0}^{N(e)}
 A(e_{i})^\alpha\right) \log^+ \left(\sum_{0}^{N(e)}
 A(e_{i})^\alpha\right)\right]<\infty,$$

In the critical case, we have the following
\begin{proposition}\label{reccrit}We suppose $p=1$, $m>1$ and (H1). We also suppose that $\rho'(1)=E_q\left[\sum_{i=1}^{N(e)}A(e_i)\log(A(e_i))\right]$ is defined and that $\rho$ is finite in a small neighborhood of $1$. Then,

\begin{itemize}

\item if $\rho'(1)<0$, then the walk is a.s. null recurrent,
 conditionally on the system's survival,

\item if $\rho'(1)= 0$ and for some
 $\delta>0,$
$$E_{\tt MT}[N(e)^{1+\delta}]<\infty,$$ then the walk is
 a.s. null recurrent, conditionally on the system's survival,
\item if $\rho'(1)>0,$ and if for some $\eta>0,$ $\omega(x,\overleftarrow{x})>\eta$ almost surely, then the walk is almost surely positive recurrent. 

\end{itemize} 
\end{proposition}
{\bf Remark:} The distinction between the case $\rho'(1)= 0$ and $\rho'(1)>0$ is quite unexpected. 
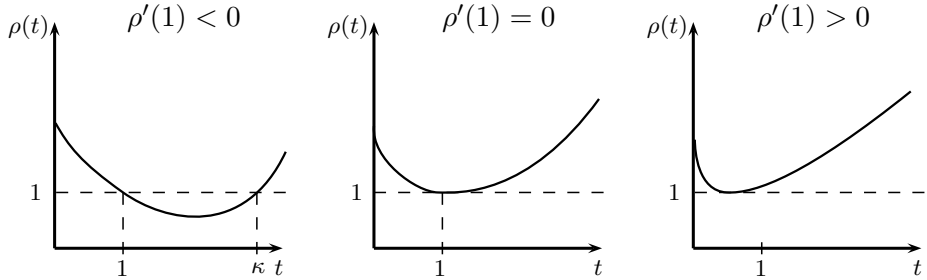
\begin{figure}[h!]

\scalebox{1} 
{
\begin{pspicture}(0,-1.9529166)(12.29,1.9929167)
\rput(5.04,-1.2470833){\psaxes[linewidth=0.04,labels=none,ticks=none,ticksize=0.10583333cm]{->}(0,0)(0,0)(3,3)}
\rput(0.84,-1.2470833){\psaxes[linewidth=0.04,labels=none,ticks=none,ticksize=0.10583333cm]{->}(0,0)(0,0)(3,3)}
\rput(9.24,-1.2470833){\psaxes[linewidth=0.04,labels=none,ticks=none,ticksize=0.10583333cm]{->}(0,0)(0,0)(3,3)}
\psline[linewidth=0.02cm,linestyle=dashed,dash=0.16cm 0.16cm](0.86,-0.50708336)(3.88,-0.50708336)
\psline[linewidth=0.02cm,linestyle=dashed,dash=0.16cm 0.16cm](5.04,-0.50708336)(8.06,-0.50708336)
\psline[linewidth=0.02cm,linestyle=dashed,dash=0.16cm 0.16cm](9.26,-0.50708336)(12.28,-0.50708336)
\rput(0.59,-0.50708336){\footnotesize $1$}
\rput(4.81,-0.50708336){\footnotesize $1$}
\rput(8.99,-0.50708336){\footnotesize $1$}
\rput(0.5,1.6929166){\footnotesize $\rho(t)$}
\rput(4.7,1.6929166){\footnotesize $\rho(t)$}
\rput(8.9,1.6929166){\footnotesize $\rho(t)$}
\rput(3.79,-1.5070833){\footnotesize $t$}
\rput(7.99,-1.5070833){\footnotesize $t$}
\rput(12.19,-1.5070833){\footnotesize $t$}
\rput(2.54,1.7979167){$\rho'(1)<0$}
\rput(6.68,1.7979167){$\rho'(1)=0$}
\rput(10.82,1.7979167){$\rho'(1)>0$}
\psline[linewidth=0.02cm](10.14,-1.1670833)(10.14,-1.3270833)
\rput(1.73,-1.5070833){\footnotesize $1$}
\rput(5.91,-1.5070833){\footnotesize $1$}
\rput(10.13,-1.5070833){\footnotesize $1$}
\psbezier[linewidth=0.03](0.8520847,0.41291666)(1.02,0.11291666)(1.18,-0.107083336)(1.74,-0.50708336)(2.3,-0.90708333)(3.34,-1.1270833)(3.88,0.032916665)
\psline[linewidth=0.02cm,linestyle=dashed,dash=0.16cm 0.16cm](1.74,-1.3270833)(1.74,-0.50708336)
\psline[linewidth=0.02cm,linestyle=dashed,dash=0.16cm 0.16cm](3.5,-1.3270833)(3.5,-0.50708336)

\rput(3.52,-1.5070833){\footnotesize $\kappa$}
\psbezier[linewidth=0.03](5.04,0.3734329)(5.0,-0.06708334)(5.5935097,-0.5005014)(5.88,-0.50708336)(6.1664906,-0.51366526)(7.06,-0.56708336)(8.0,0.73291665)
\psbezier[linewidth=0.03](9.26,0.19291666)(9.26,0.032916665)(9.3,-0.48708335)(9.68,-0.50708336)(10.06,-0.52708334)(10.74,-0.26708335)(12.1,0.8329167)
\psline[linewidth=0.02cm,linestyle=dashed,dash=0.16cm 0.16cm](5.94,-0.50708336)(5.94,-1.3270833)
\end{pspicture} 
}
\caption{Possible shapes for $\rho$ in the critical case}
\end{figure}

The study of the critical case turns out to be quite interesting, indeed several different behaviors appear in this case. The quantity $\kappa=\inf\{t>1,\, \rho(t)>1\},$ associated to $q$ is of particular interest.
When $\rho'(1)\geq 0,$ for regular trees and identically distributed $A(x)$, 
Y. Hu and Z. Shi showed (\cite{hu2007smr}) that there exist constants $0<c_1\leq c_2<\infty$ such that
$$c_1\leq \liminf_{n\rightarrow \infty} \frac{\max_{0<s<n}|X_s|}{(\log n)^3}\leq \limsup_{n\rightarrow \infty} \frac{\max_{0<s<n}|X_s|}{(\log n)^3}\leq c_2, \; \mathbb{P}-a.s..$$
It was recently proven by G. Faraud, Y. Hu and Z. Shi that $\frac{\max_{0<s<n}|X_s|}{(\log n)^3}$ actually converges to an explicit constant (see \cite{gyzalmost}). Interestingly, this constant has a different form when $\rho'(1)= 0$ and when $\rho'(1)>0$.
 
In the case $\rho'(1)<0$, Y. Hu and Z. Shi showed (\cite{hu2007sbr}) that
$$\lim_{n\rightarrow \infty} \frac {\log\max_{0<s<n}|X_s|}{\log n}=1-\frac{1}{\min\{\kappa,2\}}, \; \mathbb{P}-a.s..$$
Results in the case $p<1$ have also been obtained by Y. Hu and Z. Shi (\cite{hu2007sbr}), and the case $p>1$ has been studied by E. Aidekon (\cite{aidekon2008trw}).

Let us go back to the critical case. Our aim is to study what happens when $\kappa$ is large. When $\kappa\geq 2$, the walk behaves asymptotically like $n^{\frac{1}{2}}$. Our aim is to get a more precise estimate in this case. However we are not able to cover the whole regime $\kappa\in [2,\infty]$.

We first introduce the ellipticity assumptions 
\begin{equation} \label{ellipticity} \exists \text{ } 0<\varepsilon_{0}<
 \infty;\, \forall i,\, \varepsilon_{0}\leq A(e_{i})\leq \frac{1}{\varepsilon_0}, \, q-a.s. \end{equation}
and we assume that $(A(e_{i}))_{1\le i \le N(e)}$ is of the form $(A'(i)\mathds{1}_{(i\le N(e)})_{i\ge 1},$ where $(A'(i))_{i\ge 1}$ is a i.i.d. family independent of $N(e)$ and that $E_{q}[N(e)^{\kappa+1}]<\infty.$ (H2)

\noindent\textbf{Remark :} We actually only need this assumption to show Lemma \ref{espwetoile}, we can, for example, alternatively suppose that
 \begin{equation}\label{conditionnement}\exists N_0\;; N(e)\leq N_0,\, q-p.s. \text{ et } P_q[N\geq 2|A(e_1)]\geq\frac{1}{N_0}
\end{equation}

Note furthermore that those conditions imply (H1).

\begin{theorem}\label{clt}
Suppose $N(e)\geq 1$, $q-a.s.$, (\ref{ellipticity}), (\ref{conditionnement}).

If $p=1$, $\rho'(1)<0$ and $\kappa\in(8,\infty]$,
then there is a deterministic constant $\sigma>0$ such that,
for {\tt MT} almost every tree $T,$ 
the process $\{|X_{\lfloor n t\rfloor}|/\sqrt{\sigma^2 n}\}$ converges in law
to the absolute value of a standard Brownian motion, as $n$ goes to infinity. 
\end{theorem}
\textbf{Remark :}
This result is a generalization of a central limit theorem proved by Y. Peres and O. Zeitouni \cite{peres-zeitouni} in the case of a biased standard random walk on a Galton-Watson tree. In this case, $A(x)$ is a constant equal to $\frac{1}{m}$, therefore $\kappa=\infty$. Our
proof follows the same lines as theirs.

In the annealed setting, things happen to be easier, and we can weaken the assumption on $\kappa$.
\begin{theorem}\label{aclt}
Suppose $N(e)\geq 1$, $q-a.s.$, (\ref{ellipticity}), (\ref{conditionnement}).
If $p=1$, $\rho'(1)<0$ and $\kappa\in(5,\infty]$,
then there is a deterministic constant $\sigma>0$ such that,
under $\mathbb{P}_{\tt MT}$, 
the process $\{|X_{\lfloor n t\rfloor}|/\sqrt{\sigma^2 n}\}$ converges in law
to the absolute value of a standard Brownian motion, as $n$ goes to infinity. 
\end{theorem}

\noindent{\bf Remark :} As we will see, the annealed CLT will even be true for $\kappa \in (2,\infty),$ on a different kind of tree, following a distribution that can be described as ``the invariant distribution'' for the Markov chain of the ``environment seen from the particle''.

We thank P. Mathieu for indicating to us the technique of C. Kipnis and S.R.S. Varadhan (\cite{kipnis1986central}), that was quite an inspiration for us.

Our article will be organized as follows
\begin{itemize}
\item In section \ref{sec2} we show Theorem \ref{rec}.
\item In section \ref{sec3} we introduce a new law on trees, with particular properties.
\item In section \ref{sec4} we show a Central Limit Theorem for random walks on trees following the ``new law''.
\item In section \ref{sec5} we expose a coupling between the original law and the new one.
\item In section \ref{sec6} we show some lemmas.
\item In section \ref{sec7} we show Theorem \ref{aclt}
\end{itemize}

\section{Proof of Theorem \ref{rec}.} \label{sec2}

Let us first introduce an associated martingale, which will be
of frequent use in the sequence.

Let $\alpha \in\mathbb{R}^+$ and
\begin{equation*}
 Y_{n}^{(\alpha)} = \sum_{x \in T_{n}} \prod_{e< z \leq x} A(z)^\alpha = \sum_{x \in T_{n}} C_{x}^\alpha.
\end{equation*}
$Y_{n}^{(\alpha)}$ is known as Mandelbrot's Cascade. 

It it is easy to see that if $\rho(\alpha)<\infty$ then
$\frac{Y_{n}^{(\alpha)}}{\rho(\alpha)^n}$
is a non-negative martingale, with a.s. limit $Y^{(\alpha)}$.

We have the following theorem, due to J.D. Biggins (1977) (see \cite{Biggins:1977, biggins-kyprianou}) that allows us to know when $Y^{(\alpha)}$ is non trivial.
\begin{statement}[Biggins] \label{biggins} Let $\alpha\in \mathbb{R}^+$.
Suppose $\rho$ is finite in a small neighborhood of $\alpha$, and $\rho'(\alpha)$ exists and is finite, then the following are equivalent
\begin{itemize}
\item given non-extinction, $Y^{(\alpha)}>0$ a.s.,
\item $P_{\tt MT}[Y^{(\alpha)}=0]<1$,
\item $E_{\tt MT}[Y^{(\alpha)}]=1$,
\item $E_q\left[\left(\sum_{0}^{N(e)}
 A(e_{i})^\alpha\right) \log^+ \left(\sum_{0}^{N(e)}
 A(e_{i})^\alpha\right)\right]<\infty$, and

 (H2):= $\alpha \rho'(\alpha)/\rho(\alpha)< \log \rho(\alpha)$,
\item $\frac{Y^{(\alpha)}}{\rho(\alpha)}$ converges in $L^1$.
\end{itemize}

\end{statement} 
This martingale is related to some branching random walk, and has been
intensively studied
(\cite{Mandelbrot:1974,Biggins:1977,biggins-kyprianou,Liu:2000,Liu:2001,menshikov-petritis}).
We will see that it is closely related to our problem.

Let us now prove Theorem \ref{rec}. 
We shall use the following lemma, whose proof is similar to the proof
presented in page 129 of \cite{Lyons:1992} and omitted.
\begin{lemma}
$$\min_{0\leq t\leq1} E\left[\sum_{x\in T_{1}}A(x)^t\right]=\max_{0<y\leq1}\inf_{t>0} y^{1-t}E\left[\sum_{x\in T_{1}}A(x)^t\right].$$ \label{lem}
\end{lemma}

\textbf{(1)} Let us begin with the subcritical case, We suppose there exists some $0<\alpha <1$ such that $\rho(\alpha)=\inf_{0\leq t<1} \rho(t) <1.$
Then, following \cite{Kemeny:1976} (Prop 9-131), and standard electrical/capacited network theory, if the conductances have
finite sum, then the random walk is positive recurrent, the electrical network has zero conductance a.s., and the capacited network admits no flow a.s.. We have
$$\sum_{x\in T^{*}}C_{x}^\alpha = \sum_{n=0}^\infty \sum_{x\in T_{n}} C_{x}^\alpha=\sum_{n} {\rho(\alpha)^n} Y_{n}^{(\alpha)}.$$ Since $Y_{n}^{(\alpha)}$ is bounded (actually it converges to $0$), we have
$$\sum_{x\in T^{*}}C_{x}^\alpha <\infty, \;{\tt MT}-a.s.. $$
This implies that $a.s.$, for all but finitely many $x$, $C_{x}<1$, and then $C_{x}\leq C_{x}^\alpha$, which gives the result.

\textbf{(2)} As before, we have $\alpha $ such that
$\rho(\alpha)=\inf_{0\leq t\leq1} \rho(t) \leq 1$.
 We have to distinguish two cases. Either $\rho'(1)\geq 0$, therefore it is easy to see that, for $\alpha$, (H2) is not verified, so 
$$\sum_{x\in T_{n}} C_{x}^\alpha =Y_{n}^{(\alpha)} \rightarrow 0,$$ when n
goes to $\infty$ . Then for $n$ large enough, $C_{x}<1$ for every $x\in T_{n}$,
whence $$\sum_{x\in T_{n}} C_{x} \rightarrow 0,$$ then by the
\emph{max-flow min-cut} theorem, the associated capacited network admits no
flow a.s., this implies that no electrical current flows, and that the
random walk is recurrent {\tt MT}-a.s..

We now deal with the case where $\rho'(1)<0$, then $\alpha=1$. The proof is similar to \cite{Lyons:1992}, but, as it is quite short, we give it for the sake of clarity.
We use the fact that, if the capacited network admits no flow from $e$, then the walk is recurrent.

We call $F$ the maximum flows from $e$ in $T$, and for $x\in T, \;|x|=1$, we call $F_x$ the maximum flow in the subtree $T_x=\{y\in T, x\leq y \}$, with capacity $\frac{C_y}{A(x)}$ along the edge $(\overleftarrow{x},x)$. It is easy to see that $F$ and $F_x$ have the same distribution, and that
\begin{equation}\label{recurrence}F=\sum_{|x|=1}A(x)(F_x\wedge 1).\end{equation}
Taking the expectation yields
$$E[F]=E[F_x\wedge 1]=E[F\wedge 1],$$
therefore $ess\sup F\leq 1$.
By independence, we obtain from (\ref{recurrence}) that 
$$ess \sup F =(ess\sup \sum_{|x|=1}A(x) )(ess \sup F ).$$
This implies that $F=0$ almost surely, as $(ess\sup \sum_{|x|=1}A(x) )>1,$ when $\sum_{|x|=1}A(x)$ is not identically equal to $1$ .

\textbf{(3)} We shall use the fact that, if the water flows when $C_{x}$ is reduced exponentially in $|x|$, then the electrical current flows, and the random walk is transient a.s. (see \cite{Lyons:1989}).

We have $$\inf_{\alpha\in [0,1]}E\left[\sum_{0}^{N(e)} A(e_{i})^\alpha\right]=p >1$$
($p$ can be infinite, in which case the proof still applies).

We introduce the measure $\mu_{n}$ defined as
$$\mu_{n} (A)=E[\sharp (A\cap\{\log C_{x}\}_{x\in T_{n}})],$$
where $\sharp$ denotes the cardinality.

One can easily check that
$$\phi_{n}(\lambda):=\int_{-\infty}^{+\infty}{e^{\lambda t} d\mu_{n}(t)}=E\left[\sum_{x\in T_{n}}C_{x}^\lambda\right]=\rho(\lambda)^n.$$
Let $y\in (0,1]$ be such that $p= \inf_{t>0} y^{1-t}E[\sum_{x\in T_{1}}A({x})^t]$.
Then, using Cramer-Chernov theorem (and the fact that the probability
measure $\mu_n/m^n$ has the same Laplace transform as the sum of $n$ independent random variables with law $\mu_{1}/m$), we have $$\frac{1}{n}\log{\mu_{n}([n(-\log y),\infty))}\rightarrow \log (p/y).$$

Now, if we set $1/y<q<p/y$, there exists $k$ such that 
$$E[\sharp\{x\in T_{k}|C_{x}>y^k\}]>q^k.$$
Then the end of the proof is similar to the proof in \cite{Lyons:1992}.
We chose a small $\epsilon>0$ such that,
$$E[\sharp\{x\in T_{k}|C_{x}>y^k,\text{ and }\forall e<z\leq x,A({z})>\epsilon\}]>q^k.$$
Let $T^{k}$ be the tree whose vertices are $\{x\in T_{k n},n\in \mathbb{N}\}$ such that $x=\overleftarrow{y}$ in $T^{k}$ iff $x\leq y$ in $T$ and $|y|=x+k$.
We form a random subgraph $T^k(\omega)$ by deleting the edges $(x,y)$ where
$$\prod_{x<z\leq y}A(z)<q^k\text{ or }\exists x<z\leq y, A(z)<\epsilon.$$
Let $\Gamma_{0}$ be the connected component of the root. The tree $\Gamma_0$ is
a Galton-Watson tree, such that the expected number of children of a
vertex is $q^k>1$, hence with a positive probability $\Gamma_0$ is
infinite and has branching number over $q^k$.

 Using Kolmogoroff's 0-1 Law, conditionally to the survival there is almost surely a infinite connected component, not necessarily containing the root. 
 This connected component has branching number at least $q^k$. Then we can
 construct almost surely a subtree $T'$ of $T$, with branching number
 over $q$, such that $\forall x\in T',\,A(x)>\epsilon$ and if $|x|=n k,
 |y|=(n+1)k$ and $x<y$ then 
$\prod_{x<z\leq y}A(z)>q^k$. This implies the result.
\\

We now turn to the proof of Proposition \ref{reccrit}. 
Let $\pi$ be an invariant measure for the Markov chain $(X_{n},P_{T})$ (that is a measure on $T$ such that, $\forall x\in T,$ $\pi(x)=\sum_{y\in T}\pi(y)\omega(y,x)$), then one can easily check that
$$\pi(x)=\frac{\pi(e)\omega(e,\overleftarrow{e})}{\omega(x,\overleftarrow{x})}\prod_{0<
 z \leq x}A(z),$$
with the convention that a product over an empty set is equal to $1$.

Then almost surely there exists a constant $c>0$ (dependant of the tree) such that
$$\pi(x)>c\, C_x.$$
Thus $$\sum_{x\in T} \pi(x) >c \sum_{n} Y_{n}^{(1)}.$$

-If $\rho'(1)<0$, then (H2) is verified and $Y>0$
a.s. conditionally to the survival of the system, thus the invariant
measure is infinite and the walk is null recurrent.

-If $\rho'(1)=0,$ we use a recent result from Y. Hu and Z. Shi. 
In \cite{yzpolymer} it was shown that, under the assumptions of Theorem \ref{reccrit}, there exists a sequence $\lambda_n$ such that 
$$0<\liminf_{n\rightarrow \infty}\frac{\lambda_n}{n^{1/2}}\leq\limsup_{n\rightarrow \infty}\frac{\lambda_n}{n^{1/2}}<\infty$$ and 
 $\lambda_{n} Y_{n}^{(1)}\rightarrow_{n\rightarrow\infty} Y$, with $Y>0$ conditionally on the system's survival. The result follows easily.

-If $\rho'(1)>0$, there exists $0<\alpha<1$ such that $\rho(\alpha)=1,$ $\rho'(\alpha)=0.$ We set, for every $x\in T$, $\tilde{A}(x):=A(x)^{\alpha}$.
We set accordingly $\tilde{C}(x)=\prod_{0< z \leq x} \tilde{A}(z),$ and 
$$\tilde{\rho}(t):=E_q\left[\sum_{i=1}^{N(e)}\tilde{A}(e_i)^t\right]=\rho(\alpha t).$$
Note that $\tilde{\rho}(1)=1=\inf_{0<t\leq1}\rho(t)$ and $\tilde{\rho}'(1)=0.$
Note that under the ellipticity condition $\omega(x,\overleftarrow{x})>\eta$, for some constant $c>0$
$$\sum_{x\in T} \pi(x)<c\sum_{x\in T} C_x=\sum_{x\in T} \tilde{C}_x^{1/\alpha}.$$
Using Theorem 1.6 of \cite{yzpolymer} with $\beta=1/\alpha$ and $\tilde{C}_x=e^{-V(x)}$, we get that for any $\frac{2}{3}\alpha<r<\alpha$, 
$$E_{\tt MT}\left[\left(\sum_{x\in T_n}C_x\right)^r\right]=n^{-\frac{3r}{2\alpha} + o(1)}.$$
Note that as $r<1$,
$$\left(\sum_n Y_n^{(1)}\right)^r\leq \sum_n \left(Y_n^{(1)}\right)^r,$$
whence, using Fatou's Lemma, 
 $$E_{\tt MT}\left[\left(\sum_{x\in T}C_x\right)^r\right]<\infty.$$
 This finishes the proof.
\section{The {\tt IMT} law.} \label{sec3}

We consider trees with a marked ray, which are composed of a semi infinite ray, called $Ray=\{v_0=e,v_1=\overleftarrow{v_0},v_2=\overleftarrow{v_1}...\}$ such that to each $v_i$ is attached a tree. That way $v_i$ has several children, one of which being $v_{i-1}$.

 As we did for usual trees, we can ``mark'' these trees with
 $\{A(x)\}_{x\in T}$. Let $\tilde{\mathbb{T}}$ be the set of such trees.

Let ${\cal F}_n$ be the sigma algebra $\sigma(N_x, A_{x_i}, v_n\leq x)$ and
${\cal F}_\infty=\sigma({\cal F}_n,n\geq 0)$. While unspecified, ``measurable'' will mean
``${\cal F}_\infty$ - measurable''. 

Let $\hat{q}$ be the law on $\mathbb{N}\times \mathbb{R}_{+}^{*\mathbb{N^*}}$ defined by 
$$\frac{d\hat{q}}{d q}=\sum_1^{N(e)} A(e_i).$$
\textbf{Remark :} For this definition to have any
sense, it is fundamental that $E_q[\sum_1^{N(e)} A_i]=1$, which is provided by the
assumptions $\rho'(1)<0$ and $p=1$. 
\\

Following \cite{peres-zeitouni}, let us introduce some laws on marked trees with a marked ray.
Fix a vertex $v_0$ (the root) and a semi infinite
ray, called $Ray$ emanating from it.
To each vertex $v \in Ray$ we attach independently a set of marked vertices with law
 $\hat{q}$, except to the root $e$ to which we attach a set of children with law $(q+\hat{q})/2$.
We chose one of these vertices, with probability $\frac{A(v_i)}{\sum
 A(v_{i})}$, and identify it with the child of v on $Ray$. Then we attach a tree with law
 \texttt{MT} to the vertices not on $Ray$. 
We call {\tt IMT} the law obtained.
\begin{figure}[h!]
\centering
\scalebox{1} 
{
\begin{pspicture}(0,-4.38)(9.76,4.4)
\definecolor{color125b}{rgb}{0.6,0.6,0.6}
\psline[linewidth=0.04cm](6.92,2.12)(1.3,-3.14)
\psline[linewidth=0.04cm,linestyle=dotted,dotsep=0.16cm](1.26,-3.2)(0.0,-4.36)
\psline[linewidth=0.04cm](6.94,2.1)(6.1,3.36)
\psline[linewidth=0.04cm](6.92,2.12)(6.94,3.34)
\psline[linewidth=0.04cm](6.92,2.1)(7.74,3.32)
\psline[linewidth=0.04cm](5.64,0.92)(4.92,2.14)
\psline[linewidth=0.04cm](5.64,0.94)(5.74,2.14)
\psline[linewidth=0.04cm](4.32,-0.28)(3.32,0.94)
\psline[linewidth=0.04cm](4.34,-0.28)(4.0,0.96)
\psline[linewidth=0.04cm](4.34,-0.28)(4.58,0.94)
\psline[linewidth=0.04cm](3.08,-1.48)(2.52,-0.28)
\psline[linewidth=0.04cm](3.08,-1.48)(3.24,-0.28)
\rput{-180.0}(5.04,0.38){\pstriangle[linewidth=0.0020,dimen=outer,fillstyle=solid,fillcolor=color125b](2.52,-0.32)(0.44,1.02)}

\rput(2.54,0.58){\footnotesize ${\tt MT}$}
\rput{-180.0}(6.48,0.38){\pstriangle[linewidth=0.0020,dimen=outer,fillstyle=solid,fillcolor=color125b](3.24,-0.32)(0.44,1.02)}

\rput(3.26,0.58){\footnotesize ${\tt MT}$}
\rput{-180.0}(6.68,2.82){\pstriangle[linewidth=0.0020,dimen=outer,fillstyle=solid,fillcolor=color125b](3.34,0.9)(0.44,1.02)}

\rput(3.36,1.8){\footnotesize ${\tt MT}$}
\rput{-180.0}(8.04,2.82){\pstriangle[linewidth=0.0020,dimen=outer,fillstyle=solid,fillcolor=color125b](4.02,0.9)(0.44,1.02)}

\rput(4.02,1.8){\footnotesize ${\tt MT}$}
\rput{-180.0}(9.16,2.82){\pstriangle[linewidth=0.0020,dimen=outer,fillstyle=solid,fillcolor=color125b](4.58,0.9)(0.44,1.02)}

\rput(4.6,1.8){\footnotesize ${\tt MT}$}
\rput{-180.0}(9.88,5.22){\pstriangle[linewidth=0.0020,dimen=outer,fillstyle=solid,fillcolor=color125b](4.94,2.1)(0.44,1.02)}

\rput(4.96,3.0){\footnotesize ${\tt MT}$}
\rput{-180.0}(11.48,5.22){\pstriangle[linewidth=0.0020,dimen=outer,fillstyle=solid,fillcolor=color125b](5.74,2.1)(0.44,1.02)}

\rput(5.76,3.0){\footnotesize ${\tt MT}$}
\rput{-180.0}(12.24,7.66){\pstriangle[linewidth=0.0020,dimen=outer,fillstyle=solid,fillcolor=color125b](6.12,3.32)(0.44,1.02)}

\rput(6.14,4.22){\footnotesize ${\tt MT}$}
\rput{-180.0}(13.88,7.62){\pstriangle[linewidth=0.0020,dimen=outer,fillstyle=solid,fillcolor=color125b](6.94,3.3)(0.44,1.02)}

\rput(6.96,4.2){\footnotesize ${\tt MT}$}
\rput{-180.0}(15.48,7.62){\pstriangle[linewidth=0.0020,dimen=outer,fillstyle=solid,fillcolor=color125b](7.74,3.3)(0.44,1.02)}

\rput(7.76,4.2){\footnotesize ${\tt MT}$}

\rput(7.14,1.86){\footnotesize $v_0$}

\rput(5.9,0.66){\footnotesize $v_1$}

\rput(4.62,-0.54){\footnotesize $v_2$}

\rput(3.32,-1.7){\footnotesize $v_3$}
\psline[linewidth=0.02cm,linestyle=dashed,dash=0.16cm 0.16cm](0.94,3.34)(9.38,3.34)
\psline[linewidth=0.02cm,linestyle=dashed,dash=0.16cm 0.16cm](0.92,2.12)(9.36,2.12)
\psline[linewidth=0.02cm,linestyle=dashed,dash=0.16cm 0.16cm](0.94,0.92)(9.38,0.92)
\psline[linewidth=0.02cm,linestyle=dashed,dash=0.16cm 0.16cm](0.94,-0.28)(9.38,-0.28)
\psline[linewidth=0.02cm,linestyle=dashed,dash=0.16cm 0.16cm](0.94,-1.48)(9.38,-1.48)
\psline[linewidth=0.02cm,linestyle=dashed,dash=0.16cm 0.16cm](0.92,-2.68)(9.36,-2.68)

\rput(1.25,3.14){\footnotesize $h=1$}

\rput(1.25,1.94){\footnotesize $h=0$}
\rput(1.3,0.72){\footnotesize $h=-1$}
\rput(1.3,-0.48){\footnotesize $h=-2$}
\rput(1.3,-1.68){\footnotesize $h=-3$}
\rput(0.84,-2.9){\footnotesize $h=-4$}
\psline[linewidth=0.02cm,linestyle=dashed,dash=0.16cm 0.16cm](7.3,2.6)(8.38,2.52)
\psline[linewidth=0.02cm,linestyle=dashed,dash=0.16cm 0.16cm](6.98,2.92)(8.38,2.56)
\psline[linewidth=0.02cm,linestyle=dashed,dash=0.16cm 0.16cm](6.78,2.42)(8.34,2.48)

\rput(9.09,2.58){\footnotesize $(q+\hat{q})/2$}
\psline[linewidth=0.02cm,linestyle=dashed,dash=0.16cm 0.16cm](5.34,1.46)(7.24,1.1)
\psline[linewidth=0.02cm,linestyle=dashed,dash=0.16cm 0.16cm](5.7,1.66)(7.06,1.14)
\psline[linewidth=0.02cm,linestyle=dashed,dash=0.16cm 0.16cm](6.48,1.68)(7.06,1.14)

\rput(7.79,1.12){\footnotesize $\hat{q}$}
\rput(2.72,-2.46){{\it Ray}}
\psline[linewidth=0.02cm,arrowsize=0.05291667cm 2.0,arrowlength=1.4,arrowinset=0.4]{<-}(2.6,-1.94)(2.72,-2.22)
\end{pspicture} 
}
\caption{The {\tt IMT} law.}

\end{figure}
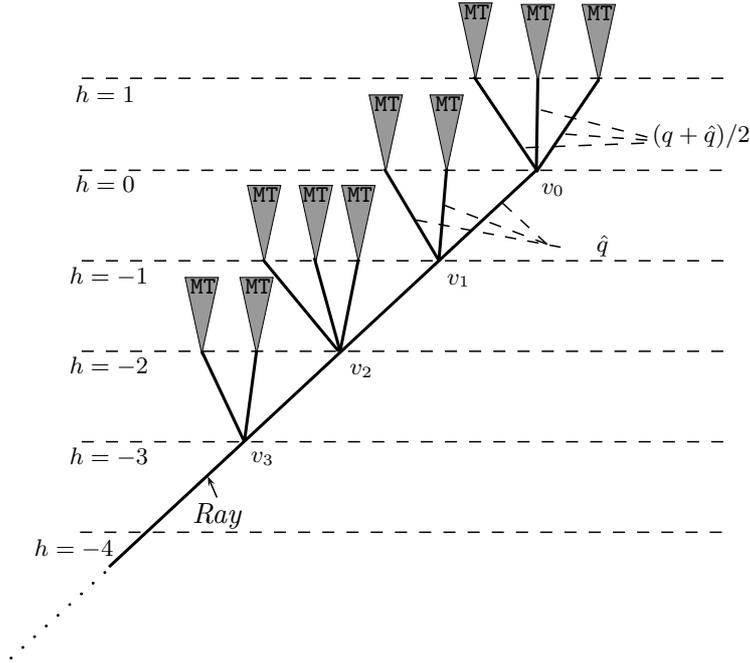

We call $\theta^v T$ be the tree $T$ ``shifted'' to $v$, that is, $\theta^v T$ has the same structure and labels as $T$, but its root is moved to vertex $v$.

Note that as before, given a tree $T$ in $\tilde{\mathbb{T}}$, we can define in a unique way a family $\omega(x,y)$ such that $\omega(x,y)=0$ unless $d(x,y)=1$, $$\forall x\in T,\, \sum_{y\in
 T}\omega(x,y)=1,$$ and \begin{equation}\label{relation}\forall x\in T,\, A(x)=\frac{\omega(
 \overleftarrow{x},x)}{\omega(\overleftarrow{x},\overleftarrow{\overleftarrow{x}})}.\end{equation}
We call random walk on $T$ the Markov chain $(X_{t}, \mathbb{P}_T)$ on $T$, starting from $v_0$ and with transition probabilities $(\omega(x,y))_{x,y\in T}$.
 
Let $T_t=\theta^{X_t}T$ denote the walk seen from the particle. $T_t$ is clearly a Markov chain on $\tilde{\mathbb{T}}$. 
We set, for
any probability measure $\mu$ on $\tilde{\mathbb{T}}$, $\mathbb{P}_\mu=\mu \otimes \mathbb{P}_{T}$ the annealed law of the random walk in a random environment on trees following the law $\mu$.
We have the following
\begin{lemma}\label{stat}
{\tt IMT} is a stationnary and reversible measure for the Markov process $T_t$, in the
sense that, for every $F: \tilde{\mathbb{T}}^2\rightarrow \mathbb{R}$ measurable,
$$\mathbb{E}_{\tt IMT} [F(T_0,T_1)]=\mathbb{E}_{\tt IMT} [F(T_1,T_0)].$$
\end{lemma} 
\textbf{Proof :}
Suppose $G$ is a ${\cal F}_n$-measurable function, that is, $G$ only depends on the (classical) marked tree of the descendants of $v_n$, to which we will refer as $T^{-n}$ and on the position of $v_0$ in the $n-th$ level of $T^{-n}$. We shall write accordingly $G(T)=G(T^{-n},v_0)$

We first show the following
\begin{lemma}\label{changement loi}
If $G$ is ${\cal F}_n$ measurable, then

\begin{equation}E_{\tt IMT}[G(T)]=E_{\tt MT}\left[\sum_{x\in T_{n}} C_x G(T,x)\left(\frac{1+\sum A(x_i)}{2}\right)\right].\label{IMT}
\end{equation}

\end{lemma}
\textbf{Remark :} These formulae seem to create a dependency on n, which
is actually irrelevant, since $E_q[\sum_{i=1}^{N(e)} A(e_i)]=1.$ 

 \noindent \textbf{Proof :}
This can be seen by an induction over $n$, using the fact that 
$$E_{\tt IMT}[G(T^{-n},v_0)]=E_q\left[\sum_{i=1}^{N}A(e_i)E\left[G(T'(i,N,A(e_j)),v_0)|i,N,A(e_j)\right]\right],$$
where $T'(x,N,A(e_i))$ is a tree composed of a vertex $v_n$ with $N$ children marked with the $A(e_i)$, and on each of this children is attached a tree with law {\tt MT}, except on the $i$-th, where we attach a tree whose law is the same as $T^{-(n-1)}$.

Iterating this argument we have

$$E_{\tt IMT}[G(T^{-n},v_0)]=E_{\tt MT}\left[\sum_{x\in T_n}C_x E\left[G(T''(x,T),x)|x,T\right]\right],$$
where the $n$ first levels of $T''(x,T)$ are similar to those of T, to each $y \in T''_n$, $x\neq y$ is attached a tree with law {\tt MT},
and to $x$ is attached a set of children with law $(\hat{q}+q)/2$, upon which
we attach {\tt MT} trees. The result follows. 
\\

Let us go back to the proof of Lemma \ref{stat}.
Using the definition of the random walk, we get

$$\mathbb{E}_{\tt IMT} [F(T_0,T_1)]=E_{\tt IMT}\left[\sum_{x\in T}\omega(v_0,x) F(T,\theta^x T)\right].$$
 \noindent Suppose $F$ is ${\cal F}_{(n-2)}\times {\cal F}_{(n-2)}$ measurable; then $T\rightarrow F(T,\theta^x T)$ is at least $F_{(n-1)}$ measurable. Then we can use (\ref{IMT}) to get
$$\mathbb{E}_{\tt IMT} [F(T_0,T_1)]=E_{\tt MT}\left[\sum_{x\in T_n}C_x\left(\frac{1+\sum A({x_i})}{2}\right)\sum_{y\in T}\omega(x,y)F(T,\theta^y T)\right].$$
It is easily verified that $$\forall x,y \in T, \omega(x,y)\frac{1+\sum A(x_i)}{2}C_x=\omega(y,x)\frac{1+\sum A(y_i)}{2}C_y.$$
 
 \noindent Using this equality, we get
\begin{eqnarray*}
\mathbb{E}_{\tt IMT} [F(T_0,T_1)]&=&E_{\tt
 MT}\left[\sum_{x\in T_n}\sum_{y\in T}\omega(y,x)C_y\left(\frac{1+\sum
 A({y_i})}{2}\right)F((T,x),(T,y))\right]\\
&=&E_{\tt
 MT}\left[\sum_{y\in T_{n+1}}\omega(y,\overleftarrow{y})C_y\left(\frac{1+\sum
 A({y_i})}{2}\right)F((T,\overleftarrow{y}),(T,y))\right]\\
&+&E_{\tt
 MT}\left[\sum_{y\in T_{n-1}}\sum_i\omega(y,y_i)C_y\left(\frac{1+\sum
 A({y_i})}{2}\right)F((T,y_i),(T,y))\right].
\end{eqnarray*}

Using (\ref{IMT}) and the fact that $F$ is ${\cal F}_{(n-2)}\times
{\cal F}_{(n-2)}$-measurable, we get
\begin{eqnarray*}\mathbb{E}_{\tt IMT} [F(T_0,T_1)]&=&E_{\tt IMT}\left[\omega(e,\overleftarrow{e})F(\theta^{\overleftarrow{e}}T,T)\right]+E_{\tt IMT}\left[\sum_i\omega(e,e_i)F(\theta^{e_i}T,T)\right]\\&=&\mathbb{E}_{\tt IMT} \left[F(T_1,T_0)\right].\end{eqnarray*}
This finishes the proof of (\ref{stat}).

\section{The Central Limit Theorem for the RWRE on {\tt IMT} Trees.} \label{sec4}
In this section we introduce and show a central limit theorem for random walk on a tree following the law {\tt IMT}.
For $T\in \tilde{\mathbb{T}}$, let $h$ be the horocycle distance on
$T$ (see Figure 2). $h$ can be defined recursively by
\begin{equation*}
\left\{ \begin{aligned}h(v_0)&=0 \\
h(\overleftarrow{x})&=h(x)-1, \;\forall x\in T
\end{aligned}\right..
\end{equation*} 
 We have the following
\begin{theorem}\label{theo}
Suppose $p=1$, $\rho'(1)<0$ and $\kappa \in [5,\infty]$, as
well as assumptions (\ref{ellipticity}) and (H2) or (\ref{conditionnement}).
There exists a deterministic constant $\sigma$ such that, for ${\tt IMT}-a.e.\text{ } T,$ 
the process $\{h(X_{\lfloor n t\rfloor})/\sqrt{\sigma^2 n}\}$ converges in distribution
to a standard Brownian motion, as $n$ goes to infinity. 
\end{theorem}

The proof of this result consists in the computation of a harmonic function $S_{x}$ on $T$. We will show that the martingale $S_{X_{t}}$ follows an invariance principle, and then that $S_{x}$ stays very close to $h(x)$. 

Let, for $v\in T,$ $$W_v=\lim_n \sum_{x\in T, v<x,d(v,x=n)}\prod_{v< z\leq x}A(z).$$
Statement \ref{biggins} implies that $W_v>0$ $a.s.$ and $E[W_v|\sigma(A(x_{i}),N(x), x<v)]=1.$
Now, let $M_0=0$ and if $X_t=v$, 
$$M_{t+1}-M_t=\left\{\begin{aligned} -W_v \text{ if } X_{t+1}=\overleftarrow{v} \\
W_{v_i},\text{ if } X_{t+1}=v_i\end{aligned}\right..$$
Given $T$, this is clearly a martingale with respect to the filtration
associated to the walk. We introduce the function 
$S_x$ defined as $S_e=0$ and for all $x\in T$,
\begin{equation}\label{defsx}
S_{x_i}=S_x+W_{x_i},
\end{equation} in such a way that $M_t=S_{X_t}$.

Let \begin{equation}\label{defeta}\eta=E_{\tt GW}[W_0^2],\end{equation} which is finite due to Theorem 2.1 of
\cite{Liu:2000} (the assumption needed for this to be true is $\kappa
>2$).
We call $$V_t:=\frac{1}{t}\sum_{i=1}^t \mathbb{E}_T
[(M_{i+1}-M{i})^2|\mathcal{F}_t]$$
the normalized quadratic variation process associated to $M_t$.
We get
$$ \mathbb{E}_T
[(M_{i+1}-M{i})^2|\mathcal{F}_t]=\omega(X_i,\overleftarrow{X_i})W_{X_i}^2+\sum_{j=1}^{N(X_{i})}
\omega(X_i,{X_i}_j)W_{{X_i}_j}^2=G(T_i),$$
where ${X_i}_j$ are the children of $X_{i}$ and $G$ is a $L^1({\tt IMT})$ function on $\tilde{\mathbb{T}}$ (again due to $\kappa>2$). 

Let us define $\sigma$ such that $E_{\tt IMT}[G(T)]:=\sigma^2\eta^2$. We have the following
\begin{proposition}\label{propo}
The process $\{M{\lfloor n t\rfloor}/\sqrt{\sigma^2\eta^2 n}\}$ converges, for ${\tt
 IMT}$ almost every T, 
to a standard Brownian motion, as $n$ goes to infinity. 
\end{proposition}
\textbf{Proof :}
We need the fact that when $t$ goes to infinity,
$$V_t\rightarrow \sigma^2\eta^2.$$ 
This comes from Birkhof's Theorem, using the transformation $\theta$ on $\tilde{\mathbb{T}}$ , which conserves the measure ${\tt IMT}$. 
The only point is to show that this transformation is ergodic, which follows from
the fact that any invariant set must be independent of
${\cal F}_n^p=\sigma(N(x), A({x_i}), v_n\leq x,h(x)<p)$, for all $n,p$, hence
is independent of $F_{\infty}$.

The result follows then from the Central Limit Theorem for martingales.
Our aim is now to show that
$h(X_t)$ and $M_t/\eta$ stay close in some sense, then the central limit theorem for $h(X_t)$ will follow easily. 

Let 
$$\epsilon_{0}<1/100, \delta\in(1/2+1/3+4\epsilon_{0},1-4\epsilon_{0})$$
and for every t, let $\rho_{t}$ be an integer valued random variable uniformly chosen in $[t, t+\lfloor t^\delta \rfloor].$ 

It is important to note that, by choosing $\epsilon_{0}$ small enough, we can get $\delta$ as close to $1$ as we need.

We are going to show the following
\begin{proposition}\label{proposition}
For any $0<\epsilon<\epsilon_{0}$, 
$$\lim_{t\rightarrow \infty} \mathbb{P}_{T}(|M_{\rho_{t}}/\eta - h(X_{\rho_{t}})|\geq \epsilon \sqrt{t})=0, \; {\tt IMT}-a.s.,$$
further,
$$\lim_{t\rightarrow \infty} \mathbb{P}_{T}\left(\sup_{r,s<t,|r-s|<t^\delta}|h(X_r)-h(X_{s})|>t^{1/2-\epsilon}\right)=0, \; {\tt IMT}-a.s..$$
\end{proposition}
Before proving this result, we need some notations. For any vertex $v$ of $T$, let
$$S_v^{\tt Ray} = \sum_{y \text{ on the geodesic connecting $v$ and
 {\tt Ray}}, y\not\in {\tt Ray}} W_y.$$

We need a fundamental result on marked Galton-Watson trees.
For a (classical) tree $T$, and $x$ in $T$, set 
$$S_x=\sum_{e<y\leq x}W_x,$$ with $W_x$ as before, and
$${\bf A}_n^{\epsilon}=\left\{v\in T, d(v,e)=n,\left|\frac{S_v}{n}-\eta\right|>\epsilon\right\}.$$
We have the following
\begin{lemma} \label{cond} Let $2<\lambda<\kappa-1,$ then for some constant $C_1$ depending on $\epsilon,$
\begin{equation}\label{ca} E_{\tt MT}\left[\sum_{x\in {\bf A}_n^{\epsilon}}C_x\right]<C_1 n^{1-\lambda/2}.\end{equation}
\end{lemma}
\textbf{Proof :} 
We consider the set $\mathbb{T}^*$ of trees with a marked path from
the root, that is, an element of $\mathbb{T}^*$ is of the form
$(T,v_0,v_1,...)$, where $T$ is in $\mathbb{T}$,
$v_0=e$ and $v_i=\overleftarrow{v_{i+1}}$.

 We consider the filtration $F_k=\sigma(T,v_1,...v_k).$ 
 Given an integer $n$, we introduce the 
 law $\widehat{\tt MT^*_n}$ on $\mathbb{T}^*$ defined as
follows : we consider a vertex $e$ (the root), to this vertex we attach
a set of marked children with law $\hat{q}$, and we chose one of
those children as $v_1$, with probability $P(x=v_1)=A(x)/\sum
A({e_i})$. To each child of $e$ different from $v_1$ we attach
independently a tree with law ${\tt MT}$, and on $v_1$ we iterate the
 process : we attach a set of children with law $\hat{q}$, we choose
 one of these children to be $v_2$, and so on, until getting to the
 level $n$. Then we attach a tree with law {\tt MT} to $v_n$.
 \begin{figure}[h!]
 \centering
\scalebox{1} 
{
\begin{pspicture}(0,-3.24)(6.02,3.24)
\definecolor{color169b}{rgb}{0.6,0.6,0.6}
\rput(3.63,-2.58){\footnotesize $\hat{q}$}

\psline[linewidth=0.02cm](1.7,-3.0)(0.9,-2.2)
\psline[linewidth=0.04cm](1.7,-3.0)(1.5,-1.42)
\psline[linewidth=0.02cm](1.7,-3.0)(2.1,-2.18)
\psline[linewidth=0.02cm](1.7,-3.0)(2.7,-2.2)
\rput{-180.0}(5.36,-3.46){\pstriangle[linewidth=0.0020,dimen=outer,fillstyle=solid,fillcolor=color169b](2.68,-2.24)(0.44,1.02)}

\rput(2.7,-1.34){\footnotesize ${\tt MT}$}
\rput{-180.0}(4.2,-3.42){\pstriangle[linewidth=0.0020,dimen=outer,fillstyle=solid,fillcolor=color169b](2.1,-2.22)(0.44,1.02)}

\rput(2.12,-1.32){\footnotesize ${\tt MT}$}
\rput{-180.0}(1.8,-3.46){\pstriangle[linewidth=0.0020,dimen=outer,fillstyle=solid,fillcolor=color169b](0.9,-2.24)(0.44,1.02)}

\rput(0.92,-1.34){\footnotesize ${\tt MT}$}
\psline[linewidth=0.02cm](1.5,-1.42)(0.68,-0.6)
\psline[linewidth=0.02cm](1.5,-1.42)(1.3,-0.62)
\psline[linewidth=0.04cm](1.5,-1.42)(2.1,0.18)

\rput(2.0,-3.02){\footnotesize $v_0$}

\rput(1.72,-1.38){\footnotesize $v_1$}
\psline[linewidth=0.02cm,linestyle=dashed,dash=0.16cm 0.16cm](1.92,-2.86)(3.38,-2.6)
\psline[linewidth=0.02cm,linestyle=dashed,dash=0.16cm 0.16cm](1.94,-2.54)(3.38,-2.6)
\psline[linewidth=0.02cm,linestyle=dashed,dash=0.16cm 0.16cm](1.66,-2.3)(3.38,-2.6)
\psline[linewidth=0.02cm,linestyle=dashed,dash=0.16cm 0.16cm](1.46,-2.76)(3.32,-2.6)
\psline[linewidth=0.02cm,linestyle=dashed,dash=0.16cm 0.16cm](1.84,-0.6)(4.5,-0.9)
\psline[linewidth=0.02cm,linestyle=dashed,dash=0.16cm 0.16cm](1.44,-1.08)(4.5,-0.9)
\psline[linewidth=0.02cm,linestyle=dashed,dash=0.16cm 0.16cm](0.94,-0.84)(4.5,-0.9)

\rput(4.79,-0.86){\footnotesize $\hat{q}$}
\rput{-180.0}(1.4,-0.26){\pstriangle[linewidth=0.0020,dimen=outer,fillstyle=solid,fillcolor=color169b](0.7,-0.64)(0.44,1.02)}

\rput(0.72,0.26){\footnotesize ${\tt MT}$}
\rput{-180.0}(2.6,-0.26){\pstriangle[linewidth=0.0020,dimen=outer,fillstyle=solid,fillcolor=color169b](1.3,-0.64)(0.44,1.02)}

\rput(1.32,0.26){\footnotesize ${\tt MT}$}

\rput(2.38,0.16){\footnotesize $v_2$}
\psline[linewidth=0.04cm,linestyle=dashed,dash=0.16cm 0.16cm](2.1,0.18)(2.1,1.4)

\rput(2.65,1.38){\footnotesize $v_{n-1}$}
\psline[linewidth=0.02cm](2.1,1.38)(1.1,2.18)
\psline[linewidth=0.02cm](2.1,1.38)(1.7,2.18)
\psline[linewidth=0.02cm](2.1,1.38)(2.3,2.18)
\psline[linewidth=0.04cm](2.1,1.38)(2.9,2.16)

\rput(3.28,2.16){\footnotesize $v_{n}$}
\rput{-180.0}(2.2,5.34){\pstriangle[linewidth=0.0020,dimen=outer,fillstyle=solid,fillcolor=color169b](1.1,2.16)(0.44,1.02)}

\rput(1.1,3.06){\footnotesize ${\tt MT}$}
\rput{-180.0}(3.4,5.34){\pstriangle[linewidth=0.0020,dimen=outer,fillstyle=solid,fillcolor=color169b](1.7,2.16)(0.44,1.02)}

\rput(1.7,3.06){\footnotesize ${\tt MT}$}
\rput{-180.0}(4.6,5.34){\pstriangle[linewidth=0.0020,dimen=outer,fillstyle=solid,fillcolor=color169b](2.3,2.16)(0.44,1.02)}

\rput(2.32,3.06){\footnotesize ${\tt MT}$}
\rput{-180.0}(5.8,5.34){\pstriangle[linewidth=0.0020,dimen=outer,fillstyle=solid,fillcolor=color169b](2.9,2.16)(0.44,1.02)}

\rput(2.92,3.06){\footnotesize ${\tt MT}$}
\end{pspicture} 
}

 \caption{the law $\widehat{\tt MT^*_n}$.}
 \end{figure}
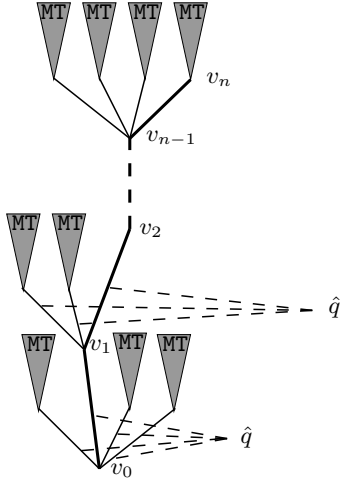
 
The same calculations as in the proof of Lemma \ref{changement loi} allow us to see the
following fact : for $f$ ${F}_n$-measurable,
\begin{equation}\label{chapeau} E_{\widehat{\tt MT^*_n}}[f(T,v_0,...,v_n)]=E_{\tt MT}\left[\sum_{x\in T_n}
 C_x f(T,p(x))\right],\end{equation} where $p(x)$ is the path from $e$ to $x$.
Note that, by construction, under $\widehat{\tt MT_n^*}$ conditionally to
 $\tilde{F}_{n}^*:=(C_{v_i},0\leq i\leq n),$ the trees $T^{(v_i)},0\leq i\leq n$ of the descendants of $v_i$ who are not descendants of $v_{i+1}$ are independent trees, and the law of $T^{(v_i)}$ is the law of a ${\tt MT}$ tree, except for the first level, whose law is $\hat{q}$ conditioned on $v_{i+1},\, A(v_{i+1})$.

For a tree T in $\mathbb{T}^*$ we have
$$W_{v_{k}}=\sum_{v_{k}=\overleftarrow{x}, x\neq v_{k+1}}A(x)W_{x}+A(v_{k+1})W_{v_{k+1}}:=W_{k}^*+A(v_{k+1})W_{v_{k+1}},$$
where $$W_{j}*=\lim_{n\rightarrow \infty} \sum_{x\in T, v_{j}<x, v_{j+1} \not\leq x ,d(v_{j},x)=n}\prod_{v\leq z\leq x}A(z).$$
Iterating this, we obtain
$$W_{v_{k}}=\sum_{j=k}^{n-1}W_{j}^*\prod_{i=k+1}^j A(v_{i})+ W_{v_{n}}\prod_{i=k+1}^n A(v_{i}),$$
with the convention that the product over an empty space is equal to one. We shall use the notation $A_{i}:=A(v_{i})$ for a tree with a marked ray.

 Finally, summing over k, we obtain
\begin{equation}\label{svn}S_{v_{n}}=\sum_{j=0}^{n-1}W_{j}^*\sum_{k=0}^j\prod_{i=k+1}^j A_{i}+ W_{v_{n}}\sum_{k=0}^n\prod_{i=k+1}^n A_{i}.\end{equation}

Let $B_{j}=\sum_{k=0}^j\prod_{i=k+1}^j A_{i}$.
We note for simplicity $W_{v_n}:=W_n^*$.
Note that $$E_{\widehat{\tt MT^*_n}}[W_0]=E_{\tt MT}\left[\left(\sum_{x\in T_n}
 C_x\right)^2\right]:=E_{\tt MT}[M_n^2]$$ converges to $\eta=E_{\tt MT}[W_0^2]$ as $n$ goes to infinity.
Indeed, recalling that $E_{\tt MT}[M_n]=1$, we have
\begin{eqnarray*}
E_{\tt MT}[(M_{n+1}-1)^2]&=&E_q\left[\left(\sum_{i=1}^{N(e)} A(e_i) U_i-1\right)^2\right]\\&=&E_q\left[\left(\sum_{i=1}^{N(e)} A(e_i)
 (U_i-1)+ \sum_{i=1}^{N(e)} A(e_i)-1\right)^2\right],\end{eqnarray*}
where, conditionally to the $A_i$, $U_i$ are i.i.d. random variables,
 with the same law as $M_n$.
We get $$E_{\tt MT}[(M_{n+1}-1)^2]=\rho(2)E_{\tt MT}[(M_{n}-1)^2]+C_2,$$
where $C_2$ is a finite number. It is easy to see then that $E[M_n^2]$
is bounded, and martingale theory implies that $M_n$ converges in $L^2$.
Using the fact that $E_{\widehat{\tt MT^*_n}}[W_{v_{k}}]=E_{\widehat{\tt MT^*_{n-k}}}[W_{0}]$, a ``Cesaro'' argument implies that $E_{\widehat{\tt MT^*_n}}[S_{v_{n}}]/n$
converges to $\eta$ as $n$ goes to infinity.
In view of that and (\ref{chapeau}) it is clear that, for n large enough
 \begin{eqnarray*}E_{\tt MT}\left[\sum_{x\in
	{\bf A}_n^{\epsilon}}C_x\right]&\leq& E_{{\tt MT}}\left[\sum_{x\in
	T_{n}}C_{x}\mathds{1}_{S_x-E_{{\widehat{\tt
	 MT^*_n}}}[S_x]>n\epsilon/2}\right]\\&\leq& P_{\widehat{\tt
	MT^*_n}}\left[\left|S_{v_{n}}-E_{{\widehat{\tt
	 MT^*_{n}}}}[S_{v_{n}}|\tilde{F}_{n}^*]\right|>\frac{n\epsilon}{4}\right]\\&+&P_{\widehat{\tt MT^*_n}}\left[\left|E_{{\widehat{\tt MT^*_n}}}[S_{v_{n}}|\tilde{F}_{n}^*]-E_{{\widehat{\tt MT^*_n}}}[S_{v_{n}}]\right|>\frac{n\epsilon}{4}\right]:=P_{1}+P_{2}.\end{eqnarray*}
Let us first bound $P_{1}$. Let $\tilde{W}_{j}^*:=W_{j}^*-E_{{\widehat{\tt MT^*_n}}}[W_{j}^*|\tilde{F}_{n}^*]$ and $\lambda \in (2,\kappa-1)$. We have
\begin{eqnarray*}E_{n}^{(1)}&:=& E_{\widehat{\tt MT^*_n}}\left[\left|S_{v_{n}}-E_{{\widehat{\tt MT^*_n}}}[S_{v_{n}}|\tilde{F}_{n}^*]\right|^\lambda\right]=E_{{\widehat{\tt MT^*_n}}}\left[\left|\sum_{i=0}^{n}\tilde{W}_{i}^*B_{i}\right|^{\lambda}\right]\\&=&E_{{\widehat{\tt MT^*_n}}}\left[E_{{\widehat{\tt MT^*_n}}}\left[\left|\sum_{i=0}^{n}\tilde{W}_{i}^*B_{i}\right|^{\lambda}|\tilde{F}_{n}^*\right]\right].\end{eqnarray*}
Inequality from page 82 of \cite{Petrov:1975} implies
$$E_{n}^{(1)}\leq C(\lambda)n^{\lambda/2-1} E_{{\widehat{\tt MT^*_n}}}\left[\sum_{i=0}^{n}E_{{\widehat{\tt MT^*_n}}}\left[\left(\tilde{W}_{i}^*B_{i}\right)^{\lambda}|\tilde{F}_{n}^*\right]\right]\leq C_3 n^{\lambda/2-1}
E_{{\widehat{\tt MT^*_n}}}\left[\sum_{i=0}^{n}B_{i}^{\lambda}\right],$$
where we have admitted the following lemma
\begin{lemma}\label{espwetoile}
$\forall \mu<\kappa,$ there exist some constant $C$ such that 
\begin{equation}E_{\tt MT^*_n}[(W_{i}^*)^{\mu}| \tilde{F}_{n}^*]<C.\label{wetoile2}\end{equation}
moreover, there exists some $\varepsilon_{1}>0$ such that
$$E_{\tt MT^*_n}[W_{i}^*| \tilde{F}_{n}^*]>\varepsilon_{1}.$$
\end{lemma}
We postpone the proof of this lemma and finish the proof of Lemma \ref{cond}.
In order to bound $ E_{{\widehat{\tt
 MT^*}}}\left[B_{i}^{\lambda}\right]$ we need to introduce a
result from \cite{biggins-kyprianou} (lemma 4.1).
\begin{statement}[Biggins and Kyprianou]\label{big}
For any $n\geq 1$ and any measurable function $G$,
$$E_{\tt MT}\left[\sum_{x\in T_{n}} C_{x} G(C_{y},e<y\leq x)\right]=E[G(e^{S_{i}};1\leq i \leq n)],$$
where $S_{n}$ is the sum of $n$ i.i.d variables whose common distribution is determined by 
$$E[g(S_1)]=E_{q}\left[\sum_{i=1}^{N(e)} A(e_{i})g(\log A(e_{i}))\right]$$ for any positive measurable function $g$.
\end{statement}
In particular, $E[e^{\lambda S_{1}}]=E_{q}[\sum_{i=1}^{N(e)} A(e_{i})^{\lambda+1}]=\rho(\lambda+1)<1.$ We are now able to compute
$$E_{{\widehat{\tt MT^*_n}}}\left[B_{n}^{\lambda}\right]=E_{\tt MT}\left[\sum_{x\in T_{n}}C_x\left( \sum_{e\leq y\leq x} \prod_{y<z\leq x} A(z)\right)^{\lambda}\right]=E\left[\left(\sum_{k=0}^n e^{S_{n}-S_{k}}\right)^{\lambda}\right].$$
Using Minkowski's Inequality, we get
\begin{equation}\label{bn}E_{{\widehat{\tt MT^*_n}}}\left[B_{n}^{\lambda}\right]\leq \left(\sum_{k=0}^n E\left[e^{\lambda(S_{k}-S_{n})}\right]^{\frac{1}{\lambda}}\right)^{\lambda}\leq \left(\sum_{k=0}^n \rho(\lambda+1)^{\frac{n-k}{\lambda}}\right)^{\lambda}\leq C_4.\end{equation}
We can now conclude, 
$$E_{n}^{(1)}\leq C_5 n^{\lambda/2},$$
and by Markov's Inequality, \begin{equation}\label{p1}P_{1}<C_6/(\epsilon^\lambda n^{{\lambda/2}}).\end{equation}

Now we are going to deal with $$P_{2}=P_{\widehat{\tt MT^*_n}}\left[\left|E_{{\widehat{\tt MT^*_n}}}[S_{v_{n}}|\tilde{F}_{n}^*]-E_{{\widehat{\tt MT^*_n}}}[S_{v_{n}}]\right|>n\epsilon/2\right].$$
Lemma \ref{espwetoile} implies that $E_{{\widehat{\tt MT^*_n}}}[W_{j}^*|\tilde{F}_{n}^*]$ is bounded above and away from zero, and a deterministic function of $A_{j+1}$. We shall note accordingly
\begin{equation}E_{{\widehat{\tt MT^*_n}}}[W_{j}^*|\tilde{F}_{n}^*]:=g(A_{j+1}).\label{wetoile}\end{equation} 
 Recalling (\ref{svn}), we have 
 $$E_{{\widehat{\tt MT^*_n}}}[S_{v_n}|\tilde{F}_{n}^*]=\sum_{j=0}^{n} E_{{\widehat{\tt MT^*_n}}}[W_{j}^*|\tilde{F}_{n}^*]){B}_{j}=\sum_{0\leq j \leq k \leq n}\prod_{i=j}^k A_{i} g(A_{k+1}).$$
with the convention $g(A_{n+1})=1$ and $A_{0}=1$.
We set accordingly 
 $$E_{{\widehat{\tt MT^*_n}}}[S_{v_n}|\tilde{F}_{n}^*]:=F(A_{1},...,A_{n}).$$
Recalling that, due to Statement \ref{big}, under the law
${{\widehat{\tt MT^*_n}}},$ the $A_{i}$ are i.i.d random variables we get $$E_{\widehat{\tt MT^*_n}}[F(A_{1},...,A_{n})]=\sum_{0\leq j \leq k \leq n}\prod_{i=j}^k E_{\widehat{\tt MT^*_n}}[A_{i}]E_{\widehat{\tt MT^*_n}}[g(A_{k+1})].$$
For $m\geq 0$ we call
\begin{multline*}F^m[A_{m+1},...,A_{n}]\\:=\sum_{\substack{0\leq j \leq k \leq n\\k\leq m-1}}\prod_{i=j}^k E_{\widehat{\tt MT^*_n}}[A_{i}]E_{\widehat{\tt MT^*_n}}[g(A_{k+1})]+\sum_{\substack{0\leq j \leq k \leq n\\k\geq m}}\prod_{i=j}^m E_{\widehat{\tt MT^*_n}}[A_i]\prod_{i'=m+1}^{k}A_{i'} g(A_{k+1}).\end{multline*}
Note that $F^0=F$ and $F^n=E_{{\widehat{\tt MT^*_n}}}[S_{v_{n}}]$, thus we can write
\begin{multline*} E_{{\widehat{\tt MT^*_n}}}[S_{v_{n}}|\tilde{F}_{n}^*]-E_{{\widehat{\tt MT^*_n}}}[S_{v_{n}}]=F^0(A_{1},...A_{n})-F^n \\
=F^0(A_{1},...A_{n})-F^1(A_{2},...A_{n})\\
+F^1(A_{2},...A_{n})-F^2(A_{3},...A_{n}) .... \\
+F^{n-1}(A_{n})-F^n. \\
\end{multline*}

We introduce the notations $\rho:=E_{\widehat{\tt MT^*_n}}[A_{1}]=\rho(2)<1,$ and for a random variable $X$, $\tilde{X}:=X-E_{\widehat{\tt MT^*_n}}[X]$. 

The last expression gives us
 \begin{multline*} 
 E_{{\widehat{\tt MT^*_n}}}[S_{v_{n}}|\tilde{F}_{n}^*]-E_{{\widehat{\tt MT^*_n}}}[S_{v_{n}}] \\
 =\tilde{g}(A_{1})+\tilde{A_{1}}(g(A_{2})+A_{2}g(A_{3})+...+\prod_{i=2}^n A_{i}g(A_{n+1})) \\
 +\rho \tilde{g}(A_{2}) + \tilde{A_2}(1+\rho)\left[\sum_{j=3}^n \prod_{i=3}^{j}A_{i}g(A_{j+1})\right]\\
 +\rho^2 \tilde{g}(A_{3}) +\tilde{A_3}(1+\rho+\rho^2)\left[\sum_{j=4}^n \prod_{i=4}^{j}A_{i}g(A_{j+1})\right]
 +... \\ +\rho^{n-1} \tilde{g}(A_{n}) +\tilde{A_n}(1+\rho+\rho^2+...\rho^{n-1}). \\
 \end{multline*}
We deduce easily that
\begin{equation}\label{spine}\left|E_{{\widehat{\tt MT^*_n}}}[S_{v_{n}}|\tilde{F}_{n}^*]-E_{{\widehat{\tt MT^*_n}}}[S_{v_{n}}]\right|<C_7+C_8\left|\sum_{k=1}^n \tilde{A}_k D_k (1+\rho+\rho^2+...\rho^{k-1})\right|,\end{equation}
where $C_{7}, C_8$ are finite constants and 
$$D_k= \sum_{j=k+1}^n \prod_{i=k+1}^{j}A_{i}g(A_{j+1}).$$

To finish the proof of Lemma \ref{cond}, we need to show that for every $\epsilon>0$, $P_{{\widehat{\tt MT^*_n}}}[\sum_{k=1}^n \tilde{A_{k}}D_{k}(1+\rho+\rho^2+...\rho^{k-1})>n\epsilon]<\frac{C(\epsilon)}{n^{\lambda/2-1}}$.

Recalling that $\lambda<\kappa -1$, we can find a small $\nu>0$ such that $\lambda(1+\nu) <\kappa-1.$ 
Then we have, by Minkowski's Inequality 
\begin{multline}\label{minkowski}E_{{\widehat{\tt MT^*_n}}}\left[D_{k}^{\lambda(1+\nu)}\right]\leq \left(\sum_{j=k+1}^{n} \left(E_{{\widehat{\tt MT^*_n}}}\left[C_8\prod_{i=k+1}^{n}A_{i}^{\lambda(1+\mu)} \right] \right)^{1/\lambda(1+\nu)} \right)^{\lambda(1+\nu)}\\
\leq \left(\sum_{j=k}^{n} \left(C_9 \rho(1+{\lambda(1+\mu)})^{n-k+1} \right) ^{1/\lambda(1+\nu)} \right)^{\lambda(1+\nu)} <C_{10}.\\
\end{multline}
Markov's Inequality then implies
\begin{equation}\label{dn}P_{{\widehat{\tt MT^*_n}}}\left[\max_{k\leq n} D_{n}>(\epsilon^2 n)^{\frac{1}{2 (1+\nu)}}\right]\leq C_{11} \frac{n}{n^{\lambda/2}\epsilon^\lambda}.\end{equation}
On the other hand, we call for $0\leq k \leq n$,
$$N_k:=\sum_{j=n-k}^n D_{j}g(A_{n+1})(1+\rho+\rho^2+...\rho^{j-1}).$$
It is easy to check that $N_k$ is a martingale with respect to the filtration $\mathcal{H}_k=\sigma(A_j, n-k\leq j \leq n).$
We can compute the quadratic variation of this martingale
$$\langle N_k \rangle:=\sum_{j=1}^{k}E_{{\widehat{\tt MT^*_n}}}[(N_k-N_{k-1})^2|\mathcal{H}_{k-1}]=\rho(3)\sum_{j=1}^{k} (D_{n-j})^2.$$
On the other hand, the total quadratic variation of $N_k$ is equal to 
$$[ N_k]:=\sum_{j=1}^{k}(N_k-N_{k-1})^2=\sum_{j=1}^{k} (\tilde{A}_{n-j}D_{n-j})^2.$$
It is easy to check that if the event in (\ref{dn}) is fullfilled, then there exists some constant $C_{12}$ such that 
$\langle N_k \rangle < C_{12} n^{1+\frac{1}{2 (1+\nu)}}$ and $[ N_k ] < C_{12} n^{1+\frac{1}{2 (1+\nu)}}$. Therefore, using (\ref{dn}) and Theorem 2.1 of \cite{bercu2008exponential},
 \begin{equation}\label{p2}P_{{\widehat{\tt MT^*_n}}}[|\sum_{k=1}^n \tilde{A_{k}}D_{k}|>n\epsilon]\leq C_{11} \frac{n}{n^{\lambda/2}\epsilon^\lambda}+2\exp{-{\frac{(\epsilon n)^2}{2C_{12} n^{1+\frac{1}{2 (1+\nu)}}}}}.\end{equation}
 
 Putting together (\ref{p1}) and (\ref{p2}), we obtain (\ref{ca}). This finishes the proof of Lemma \ref{cond}.
 In particular, if $\kappa>5$, we can choose $\lambda>4,$ so that $$E_{\tt MT}[\sum_{x\in {\bf A}_n^{\epsilon}}C_x]<n^{-\mu},$$ with $\mu>1$ 
.
The following corollary is a direct consequence of the proof.
\begin{corollary}\label{corollaire} For every $a>0$ and $2<\lambda<\kappa-1$,
$$P_{{\widehat{\tt MT^*_n}}}[|S_{v_k}-k\eta|>a]\leq C_{1} \frac{k^{1-\lambda/2}}{a^\lambda}.$$
\end{corollary}

 We now give the proof of Lemma \ref{espwetoile}.
 As we said in the introduction, for this lemma we need either the assumption (H2) or the assumption (\ref{conditionnement}). We give the proof in both cases.
 Note that, by construction of ${\tt MT}^*,$ as, using Theorem 2.1 of
\cite{Liu:2000}, for every $x$ a child of $v_{i},$ different from $v_{i+1},$ $W(x)$ has finite moments of order $\mu,$
 \begin{eqnarray}E_{\tt MT^*_n}[(W_{i}^*)^{\mu}| \tilde{F}_{n}^*]&=&C_{0}E_{\tt MT^*_n}\left[\left(\sum_{\overleftarrow{x}=v_{i},x\neq v_{i+1}}A(x)\right)^{\mu}\big| \tilde{F}_{n}^*\right]\\&=&C_{0}E_{\hat{q}}\left[\left(\sum_{|x|=1, x\neq v_{1}}A(x)\right)^{\mu}|A(v_{1})\right]\end{eqnarray}
 Note that the upper bound is trivial under assumption (\ref{conditionnement}). We suppose (H2), Let $f$ be a measurable test function, we have by construction
\begin{eqnarray*}&&E_{\hat{q}}\left[\left(\sum_{|x|=1, x\neq v_{1}}A(x)\right)^{\mu}f(A(v_{1}))\right] \\&=&E_{q}\left[\sum_{i=1}^{N(e)}A(e_{i}) \left(\sum_{i\neq j} A(e_{j})\right)^{\mu}f(A(e_{i}))\right]\\
&\leq&\sum_{n=1}^{\infty}P_{q}(N(e)=n)E_{q}\left[\sum_{i=1}^{n}A'(i) \left(\sum_{i\neq j} A'(1)\right)^{\mu}f(A'(1))\right]
 \end{eqnarray*}
 By standard convexity property, we get that the last term is lesser or equal to
\begin{eqnarray*}\sum_{n=1}^{\infty}P_{q}(N(e)=n)E_{q}\left[\sum_{i=1}^{n}A'(i) n^{\mu-1}\sum_{i\neq j} A'(j)^{\mu}f(A'(i))\right]\\\le E_{q}[A'(1)^{\mu}]\sum_{n=1}^{\infty}P_{q}(N(e)=n)n^{\mu+1}E_{q}\left[A'(i) f(A'(i))\right]\\=E_{q}[A'(1)^{\mu}]E_{q}\left[A'(i) f(A'(i))\right]E_{q}\left[N(e)^{\mu+1}\right],\end{eqnarray*}
while, still by construction
\begin{eqnarray*}
E_{\hat{q}}\left[f(A(v_{1}))\right]&=&E_{q}\left[\sum_{i=1}^{N(e)}A(e_{i})f(A(e_{i}))\right]\\&=&\sum_{n=1}^{\infty}P_{q}(N(e)=n)E_{q}[\sum_{i=1}^n A'(i)f(A'(i))]\\&=&E_{q}[N(e)]E_{q}
[A'(1)f(A'(1))].\end{eqnarray*}
Therefore the result is direct. To prove the lower bound we begin with assumption (\ref{conditionnement}). Actually we will only use the second part of this assumption, which is trivially implied by (H2), so the proof will also work for this case.

We have
\begin{eqnarray*}
&\phantom{=}&E_{\hat{q}}\left[\sum_{|x|=1, x\neq v_{1}}A(x) f(A(v_{1}))\right]=E_{q}\left[\sum_{i=1}^{N(e)}A(e_{i})\left(\sum_{i\neq j} A(e_{j})\right)f(A(e_{i}))\right]\\&\ge&\epsilon_{0}\sum_{i=1}^{\infty}E_{q}[A(e_{i})f(A(e_{i}))\mathds{1}_{\{i\le N(e)\}}(N(e)-1)]\\&\ge &\epsilon_{0}\sum_{i=2}^{\infty}E_{q}[A(e_{i})f(A(e_{i}))\mathds{1}_{\{i\le N(e)\}}(N(e)-1)]+\epsilon_{0}E_{q}[A(e_{1})f(A(e_{1}))(N(e)-1)]\\&\ge&\epsilon_{0}\sum_{i=2}^{\infty}E_{q}[A(e_{i})f(A(e_{i}))\mathds{1}_{\{i\le N(e)\}}]+\epsilon_{0}E_{q}[A(e_{1})f(A(e_{1}))P(N(e)>2|A(e_{1}))]\\
&\ge& \frac{\epsilon_{0}}{N_{0}}E_{q}\left[\sum_{i=1}^{N(e)}A(e_{i})f(A(e_{i})\right]=E_{\hat{q}}\left[f(A(v_{1}))\right],
\end{eqnarray*}
indeed for $i\ge 2$, the event $\{i<N(e)\}$ implies $N(e)-1>1$.
This finishes the proof of Lemma \ref{espwetoile}.\\

Let us go back to {\tt IMT} trees.
We consider the following sets
\begin{equation}\label{defb}{\bf B}_n^{\epsilon}=\left\{v\in T, d(v,{\tt Ray})=n,\left|\frac{S_v^{\tt
 Ray}}{n}-\eta\right|>\epsilon\right\}.\end{equation}
We can now prove the following
\begin{lemma}\label{bm}
$$\lim_{t\rightarrow \infty} \mathbb{P}_{T}(X_{\rho_{t}}\in \cup_{n=1}^\infty {\bf B}_n^\epsilon) =0,\, {\tt IMT}-a.s..$$ 
\end{lemma}
\textbf{Proof :}
we recall that a ${\tt IMT}$ tree is composed of a semi-infinite path
from the root : ${\it Ray}=\{ v_{0} = e, v_{1} = \overleftarrow{v_{0}}
... \}$, and that
 $$W_{j}^*=\lim_n \sum_{x\in T, v_{j}<x, v_{j-1} \not\leq x ,d(v_{j},x=n)}\prod_{v\leq z\leq x}A(z).$$ Recalling Lemma \ref{espwetoile}, under {\tt IMT}, conditionally to $\{{\it Ray}, A(v_{i})\}$, $W_{j}^*$ are independent random variables and $E[W_{j}^*]>\varepsilon_0.$ 

Let $1/2<\gamma <\delta$. For a given tree T, we consider the event
$$\Gamma_{t}=\{ \exists u\leq 2t | X_{u}=v_{\lfloor t^\gamma \rfloor}\}. $$
We have 
$$\Gamma_{t}\subset \{\inf_{u\leq2t} M_{u}\leq S_{v_{\lfloor t^\gamma \rfloor}}\},$$
and {\tt IMT} almost surely, for some $\epsilon$, $$S_{v_{\lfloor t^\gamma \rfloor}}\leq -\sum_{0}^{{\lfloor t^\gamma \rfloor}}W_{j}^*<-\epsilon t^{\gamma}, \text{for t large enough}.$$
Since $M_{t}$ is a martingale with bounded normalized quadratic variation $V_{t}$, we get that, for {\tt IMT} almost every tree $T$,
$$\mathbb{P}_{T}(\Gamma_{t})\rightarrow 0.$$
Going back to our initial problem, we have
\begin{eqnarray}\label{gammat}\mathbb{P}_{T}(X_{\rho_{t}}\in \cup_{m=1}^\infty {\bf B}_m^\epsilon) &\leq& \mathbb{P}_{T}(X_{\rho_{t}}\in \cup_{n=1}^\infty {\bf B}_m^\epsilon;\Gamma_{t}^c)+ 
\mathbb{P}_{T}(\Gamma_{t})\\ &\leq&\frac{1}{\lfloor t^\delta \rfloor}\mathbb{E}_{T}\left[\sum_{s=0}^{H_{v_{\lfloor t^\gamma \rfloor}}}\mathds{1}_{X_{s}\in \cup_{m=1}^\infty {\bf B}_m^\epsilon}\right]+ \mathbb{P}_{T}(\Gamma_{t}),\end{eqnarray}
where $H_{v_{\lfloor t^\gamma \rfloor}}$ is the first time the walk hits $v_{\lfloor t^\gamma \rfloor}.$

As before we call $T^{(v_i)}$ the subtree constituted of the vertices $x\in T$ such that $v_i\leq x\not\leq x$.
 The first part of the right hand term of (\ref{gammat}) is equal to
$$\frac{1}{\lfloor t^\delta \rfloor}\mathbb{E}_{T}\left[\sum_{i=0}^{\lfloor t^\gamma \rfloor}\sum_{s=0}^{H_{v_{\lfloor t^\gamma \rfloor}}}\mathds{1}_{X_{s}\in \cup_{m=1}^\infty {\bf B}_m^\epsilon\cap T^{(v_{i})}}\right]\leq\frac{1}{\lfloor t^\delta \rfloor}\sum_{i=0}^{\lfloor t^\gamma \rfloor}\mathbb{E}_{T}\left[\sum_{s=0}^{H_{v_{\lfloor t^\gamma \rfloor}}}\mathds{1}_{X_{s}=v_{i}}\right] N_{i},$$
where $N_{i}$ is the $P_T$-expectation of the number of visits to $\cup_{n=1}^\infty {\bf B}_m^\epsilon\cap T^{(v_{i})}$ during one excursion in $T^{(v_{i})}.$ Lemma \ref{cond} implies that, under {\tt IMT} conditioned on $\{{\it Ray}, A(v_{i})\}$, $N_{i}$ are independent and identically distributed variables, with finite expectation, up to a bounded constant due to the first level of those subtrees. 
We are now going to compute $\mathbb{E}_{T}\left[\sum_{s=0}^{H_{v_{\lfloor t^\gamma \rfloor}}}\mathds{1}_{X_{s}=v_{i}}\right].$
Given $T$, we have $$\sum_{s=0}^{H_{v_{\lfloor t^\gamma \rfloor}}}\mathds{1}_{X_{s}=v_{i}}\leq 1+ M_{i},$$ where $M_{i}$ is the number of times the walk, leaving from $v_{i}$, gets back to $v_{i}$ before hitting $v_{\lfloor t^\gamma \rfloor}.$ $M_{i}$ follows a geometric law, with parameter $p_{i}=\mathbb{P}_{T}^{v_{i}}[H_{v_{\lfloor t^\gamma \rfloor}}<H_{v_{i}}].$

Standard computations for random walks on $\mathbb{Z}$, (see, for example, Theorem 2.1.12 of \cite{Zeitouni:2003}) imply that
$$p_{i}=\frac{\omega(v_{i},v_{i+1})}{1+ \sum_{j=i}^{{\lfloor t^\gamma \rfloor}-1}\prod_{k=j-1}^{{\lfloor t^\gamma \rfloor}}A(v_{k})},$$
and, going back to our initial problem, 
\begin{eqnarray*}\mathbb{P}_{T}(X_{\rho_{t}}\in \cup_{m=1}^\infty {\bf B}_m^\epsilon) &\leq& \mathbb{P}_{T}(\Gamma_{t})+\frac{C_{14}}{\lfloor t^\delta \rfloor}\sum_{i=0}^{\lfloor t^\gamma \rfloor} \left(1+ \sum_{j=i}^{{\lfloor t^\gamma \rfloor}-1}\prod_{k=j-1}^{{\lfloor t^\gamma \rfloor}}A(v_{k})\right)N_{i}\\&\leq& \mathbb{P}_{T}(\Gamma_{t})+V_{t}\frac{C_{14}}{\lfloor t^\delta \rfloor}\sum_{i=0}^{\lfloor t^\gamma \rfloor}N_{i},\end{eqnarray*}
with $V_{t}=1+ \sum_{j=0}^{{\lfloor t^\gamma \rfloor}-1}\prod_{k=j-1}^{{\lfloor t^\gamma \rfloor}}A(v_{k}).$

As in the proof of Lemma \ref{cond}, statement \ref{big} implies that $E_{\tt IMT}[V_{t}^\alpha]<C_{15}$ for some $\alpha>2.$ Now we can choose $\delta$ close to one and $\gamma$ close to $1/2$, and $\mu$ such that $1/\alpha<\mu<\delta-\gamma$ 

Markov's Inequality and the Borel Cantelli Lemma imply that, ${\tt IMT}$-almost surely, there exists $t_{0}$ such that
$\forall t>t_{0}, V_{t}\leq t^{\mu},$ and then,

$$\mathbb{P}_{T}(X_{\rho_{t}}\in \cup_{n=1}^\infty {\bf B}_m^\epsilon) \leq \mathbb{P}_{T}(\Gamma_{t})+\frac{C_{16}}{\lfloor t^{\delta-\mu} \rfloor}\sum_{i=0}^{\lfloor t^\gamma \rfloor}N_{i}.$$
Since ${\delta-\mu}<\gamma$, an application of the law of large numbers finishes the proof of Lemma \ref{bm}.
\\

We are now able to prove the first part of Proposition \ref{proposition}.
Note that under {\tt IMT}, $S_{v_n}$ follows the same law as $S_{v_n}$ in a
$\mathbb{T}^*$ tree 
under $\widehat{{\tt MT}^*_n}$, whence 
$$S_{v_{n}}/n \underset{n\rightarrow\infty}{\rightarrow}- \eta $$ in
probability. Let $Q_{t}$ be the first ancestor of $X_{\rho_{t}}$ on ${\it Ray}$. Statement \ref{big} and standard RWRE theory imply that $Q_t$ is transient, therefore 
$$ S_{Q_{t}}/h(Q_{t}) \underset{t\rightarrow\infty}{\rightarrow} \eta,$$
so that, for any positive $\epsilon_{1},$ for large t, 
\begin{equation}\label{rt} |S_{Q_{t}}/\eta- h(Q_{t})|\leq\epsilon_{1} \sup_{s\leq 2t} |M_{t}|.\end{equation}
 We can now compute
$$|M_{\rho_{t}}/\eta - h(X_{\rho_{t}})|=|S^{\tt Ray}_{X_{\rho_{t}}}/\eta-d(X_{\rho_{t}},{\tt Ray})+S_{Q_{t}}/\eta-h(Q_{t})|.$$
In view of (\ref{rt}) on the event $\{X_{\rho_{t}}\not \in \cup_{n=1}^\infty {\bf B}_m^\epsilon\},$ we have
$$|M_{\rho_{t}}/\eta - h(X_{\rho_{t}})|\leq 2\epsilon_{1} \sup_{s\leq 2t} |M_{s}|.$$
The process $V_{t}$ being bounded ${\tt IMT}\; a. s.$, a standard martingale inequality implies
$$\lim_{\epsilon_{1}\rightarrow 0} \limsup_{t \rightarrow \infty} \mathbb{P}_{\it T}^0(\sup_{s\leq t} |M_{s}|>\epsilon \sqrt{t}/(2\epsilon_{1}))=0.$$
It follows that
$$\lim_{t\rightarrow \infty} \mathbb{P}_{T}(|M_{\rho_{t}}/\eta - h(X_{\rho_{t}})|\geq \epsilon \sqrt{t})=0, {\tt IMT}-a.s.$$
We are now going to prove the second part of Proposition \ref{proposition}.
The course of the proof is similar to \cite{peres-zeitouni}.
We have the following lemma
\begin{lemma} for any $u$, $t\geq1$, \label{lem3.2}
$$\mathbb{P}_{\tt MT}(|X_{i}|\geq u\text{ for some } i\leq t)\leq 2t e^{-u^2/2t}.$$
\end{lemma}
\textbf{Proof :} 
We consider the graph ${\it T}^*$ obtained by truncating the tree ${\it T}$ after the level $u-1$, and adding an extra vertex $e^*,$ connected to all vertices in $T_{u-1}.$
We construct a random walk $X_{s}^*$ on ${\it T}^*$ as following 
$$\mathbb{P}^0_{\it T} (X_{i+1}^*=y|X_{i}^*=x)=\left\{\begin{aligned}& \omega(x,y) \text{ if } |x|<u-1 \text{ or } |x|=u-1,|y|=u-2 \\
&1-\omega(x,\overleftarrow{x}) \text{ if } |x|=u-1,y=e^* 
\\ &\tilde{\omega}(e^*,y) \text{ if } x=e^*,|y|=u-1 
\end{aligned}\right. .$$
We can choose $\tilde{\omega}(e^*,y)$ arbitrarily, provided $\sum_{y\in T_{u-1}} \tilde{\omega}(e^*,y)=1,$ so we will use this choice to ensure the existence of an invariant measure : indeed, if $\pi$ is an invariant measure for the walk, one can easily check that, for any $x$ such that $|x|\leq u-1,$ calling $x^{(1)}$ the first vertex on the path from $e$ to $x$,
$$\pi(x)=\frac{\pi(e)\omega(e,x^{(1)})}{\omega(x,\overleftarrow{x})}\prod_{x^{(1)}< z \leq x}A(z).$$
Further, we need that, for every $x\in T_{u-1},$
$$\pi(x)(1-\omega(x,\overleftarrow{x}))=\pi(e^*)\tilde{\omega}(e^*,x).$$
Summing over x, and using $\sum_{y\in T_{u}} \tilde{\omega}(e^*,y)=1,$ we get
\begin{align*}\pi(e^*)&=\pi(e)\sum_{x\in T_{u-1}} \omega(e,x^{(1)}) \prod_{x^{(1)}< z \leq x}A(z) \frac{\sum \omega(x,x_{i})}{\omega(x,\overleftarrow{x})}\\&\leq \pi(e)\sum_{x\in T_{u}} \prod_{x^{(1)}< z \leq x}A(z)\leq \pi(e) Y_{u}.\end{align*}
Then,
$$
\mathbb{P}_{\tt MT}(\exists i\leq t,\,X_{i}\geq u)\leq
\mathbb{P}_{\tt MT}(\exists i\leq t,\,X_{i}^*= e^*)\leq \sum_{i=1}^{t} \mathbb{P}_{\tt MT}(X_{i}^*= e^*).$$
By the Carne-Varnopoulos Bound (see \cite{Lyons:2005}, Theorem 12.1), 
$$\mathbb{P}_{T}(X_{i}^*= e^*)\leq 2 \sqrt{Y_{u}}e^{-u^2/2i}.$$
Since, by Jensen's Inequality, $E_{\tt MT}(\sqrt{Y_{n}})\leq1$,
$$\mathbb{P}_{\tt MT}(X_{i}\geq u\text{ for some } i\leq t)\leq 2t e^{-u^2/2t}.$$
We have the following corollary, whose proof is omitted
\begin{corollary}
$$\mathbb{P}_{\tt IMT}(|h(X_{i})|\geq u \text{ for some }i\leq t)\leq4t^3 e^{-(u-1)^2/2t}.$$
\end{corollary}
\textbf{Proof :}
see \cite{peres-zeitouni}, Corollary 2.
\\

We can now finish the proof of the second part of Proposition \ref{proposition}. 
Under $\mathbb{P}_{\tt IMT}$, the increments $h(X_{i+1})-h(X_{i})$ are stationnary, therefore, for any $\epsilon$ and $r,s\leq t$ with $|s-r|\leq t^{\delta},$ 
$$\mathbb{P}_{\tt IMT}(|h(X_{r})-h(X_{s})|\geq t^{1/2-\epsilon})\leq\mathbb{P}_{\tt IMT}(|h(X_{r-s})|\geq t^{1/2-\epsilon})\leq4t^3e^{-t^{1-\delta-2\epsilon}}.$$
Whence, by Markov's Inequality, for all t large,
$$P_{\tt IMT}\left(\mathbb{P}^0_{T}\left(|h(X_{r-s})|\geq t^{1/2-\epsilon}\right)\geq e^{-t^{1-\delta-\epsilon}}\right)\leq e^{-t^{1-\delta-\epsilon}}.$$
Consequently,
$$P_{\tt IMT}\left(\mathbb{P}^0_{T}\left(\sup_{r,s\leq t,|r-s|\leq t^{\delta}}|h(X_{r})-h(X_{s})|\geq t^{1/2-\epsilon}\right)\geq e^{-t^{1-\delta-\epsilon}}\right)\leq e^{-t^{1-\delta-\epsilon}}.$$
The Borel-Cantelli Lemma completes the proof.
\\

We are now able to finish the proof of Theorem \ref{theo}.
Due to Proposition \ref{propo}, the process $\{M{\lfloor n t\rfloor}/\sqrt{\sigma^2\eta^2 n}\}$ converges, for ${\tt
 IMT}$ almost every T, 
to a standard Brownian motion, as $n$ goes to infinity. Further, by Theorem 14.4 of \cite{Billingsley:1999}, 
$\{M{\rho_{n t}}/\sqrt{\sigma^2\eta^2 n}\}$ converges, for ${\tt
 IMT}$ almost every T, 
to a standard Brownian motion, as $n$ goes to infinity. Proposition \ref{proposition} implies that the sequence of processes $\{Y_{t}^n\}=\{h(X_{\rho_{n t}})/\sqrt{\sigma^2n}\}$ 
is tight and its finite dimensional distributions converge to those of a standard Brownian motion, therefore it converges in distribution to a standard Brownian motion, and, applying again Theorem 14.4 of \cite{Billingsley:1999},
so does $\{h(X_{\lfloor n t\rfloor}/\sqrt{\sigma^2 n}\}.$

\section{Proof of Theorem \ref{clt}.} \label{sec5}
In this section we finish the proof of Theorem \ref{clt}. 
Our argument relies on a coupling between random walks on {\tt MT} and on {\tt IMT} trees, quite similar to the coupling exposed in \cite{peres-zeitouni}. 
\begin{figure}
\scalebox{1} 
{
\begin{pspicture}(0,-7.699368)(13.749352,7.699368)
\definecolor{color1365}{rgb}{0.8,0.8,0.8}
\definecolor{color319}{rgb}{0.6,0.6,0.6}
\definecolor{color1310}{rgb}{0.4,0.4,0.4}
\psline[linewidth=0.02cm](1.869233,4.3106422)(0.50979084,5.155324)
\psline[linewidth=0.02cm](1.869233,4.3106422)(1.869233,5.155324)
\psline[linewidth=0.02cm](1.869233,4.3106422)(3.2286754,5.155324)
\psline[linewidth=0.02cm,linecolor=color1365](1.869233,5.155324)(1.3594422,6.0000052)
\psline[linewidth=0.02cm,linecolor=color1365](1.869233,5.155324)(2.3790238,6.0000052)
\psline[linewidth=0.02cm,linecolor=color1365](1.869233,5.155324)(1.869233,6.0000052)
\psline[linewidth=0.02cm,linecolor=color1365](0.50979084,5.155324)(0.8496514,6.0000052)
\psline[linewidth=0.02cm,linecolor=color1365](0.50979084,5.155324)(0.16993028,6.0000052)
\psline[linewidth=0.02cm,linecolor=color1365](3.2286754,5.155324)(2.8888147,6.0000052)
\psline[linewidth=0.02cm,linecolor=color1365](3.2286754,5.155324)(3.2286754,6.0000052)
\psline[linewidth=0.02cm,linecolor=color1365](3.2286754,5.155324)(3.5685358,6.0000052)
\psline[linewidth=0.02cm,linecolor=color1365](0.16993028,6.0000052)(0.16993028,6.844687)
\psline[linewidth=0.02cm,linecolor=color1365](0.8496514,6.0000052)(0.6797211,6.844687)
\psline[linewidth=0.02cm,linecolor=color1365](0.8496514,6.0000052)(1.0195817,6.844687)
\psline[linewidth=0.02cm,linecolor=color1365](1.3594422,6.0000052)(1.3594422,6.844687)
\psline[linewidth=0.02cm,linecolor=color1365](1.869233,6.0000052)(1.6993028,6.844687)
\psline[linewidth=0.02cm,linecolor=color1365](1.869233,6.0000052)(2.0391634,6.844687)
\psline[linewidth=0.02cm,linecolor=color1365](2.3790238,6.0000052)(2.3790238,6.844687)
\psline[linewidth=0.02cm,linecolor=color1365](2.3790238,6.0000052)(2.2090936,6.844687)
\psline[linewidth=0.02cm,linecolor=color1365](2.3790238,6.0000052)(2.5489542,6.844687)
\psline[linewidth=0.02cm,linecolor=color1365](2.8888147,6.0000052)(2.8888147,6.844687)
\psline[linewidth=0.02cm,linecolor=color1365](3.2286754,6.0000052)(3.3986056,6.844687)
\psline[linewidth=0.02cm,linecolor=color1365](3.2286754,6.0000052)(3.2286754,6.844687)
\psline[linewidth=0.02cm,linecolor=color1365](3.2286754,6.0000052)(3.058745,6.844687)
\psline[linewidth=0.02cm,linecolor=color1365](3.5685358,6.0000052)(3.5685358,6.844687)
\psline[linewidth=0.02cm,linecolor=color1365](3.5685358,6.0000052)(3.738466,6.844687)
\psline[linewidth=0.02cm,linecolor=color1365](0.16993028,6.844687)(0.0,7.6893682)
\psline[linewidth=0.02cm,linecolor=color1365](0.16993028,6.844687)(0.16993028,7.6893682)
\psline[linewidth=0.02cm,linecolor=color1365](0.16993028,6.844687)(0.33986056,7.6893682)
\psline[linewidth=0.02cm,linecolor=color1365](0.6797211,6.844687)(0.50979084,7.6893682)
\psline[linewidth=0.02cm,linecolor=color1365](0.6797211,6.844687)(0.8496514,7.6893682)
\psline[linewidth=0.02cm,linecolor=color1365](1.0195817,6.844687)(1.0195817,7.6893682)
\psline[linewidth=0.02cm,linecolor=color1365](1.3594422,6.844687)(1.1895119,7.6893682)
\psline[linewidth=0.02cm,linecolor=color1365](1.3594422,6.844687)(1.3594422,7.6893682)
\psline[linewidth=0.02cm,linecolor=color1365](1.3594422,6.844687)(1.5293725,7.6893682)
\psline[linewidth=0.02cm,linecolor=color1365](1.6993028,6.844687)(1.6993028,7.6893682)
\psline[linewidth=0.02cm,linecolor=color1365](2.0391634,6.844687)(1.869233,7.6893682)
\psline[linewidth=0.02cm,linecolor=color1365](2.0391634,6.844687)(2.0391634,7.6893682)
\psline[linewidth=0.02cm,linecolor=color1365](2.2090936,6.844687)(2.2090936,7.6893682)
\psline[linewidth=0.02cm,linecolor=color1365](2.3790238,6.844687)(2.5489542,7.6893682)
\psline[linewidth=0.02cm,linecolor=color1365](2.3790238,6.844687)(2.3790238,7.6893682)
\psline[linewidth=0.02cm,linecolor=color1365](2.5489542,6.844687)(2.7188845,7.6893682)
\psline[linewidth=0.02cm,linecolor=color1365](2.8888147,6.844687)(2.8888147,7.6893682)
\psline[linewidth=0.02cm,linecolor=color1365](3.058745,6.844687)(3.058745,7.6893682)
\psline[linewidth=0.02cm,linecolor=color1365](3.2286754,6.844687)(3.3986056,7.6893682)
\psline[linewidth=0.02cm,linecolor=color1365](3.2286754,6.844687)(3.2286754,7.6893682)
\psline[linewidth=0.02cm,linecolor=color1365](3.3986056,6.844687)(3.5685358,7.6893682)
\psline[linewidth=0.02cm,linecolor=color1365](3.5685358,6.844687)(3.738466,7.6893682)
\psline[linewidth=0.02cm,linecolor=color1365](3.738466,6.844687)(4.0783267,7.6893682)
\psline[linewidth=0.02cm,linecolor=color1365](3.738466,6.844687)(3.9083962,7.6893682)
\psline[linewidth=0.02cm](1.869233,0.31064233)(0.50979084,1.1553239)
\psline[linewidth=0.02cm](1.869233,0.31064233)(1.869233,1.1553239)
\psline[linewidth=0.02cm](1.869233,0.31064233)(3.2286754,1.1553239)
\psline[linewidth=0.02cm,linecolor=color1365](1.869233,1.1553239)(1.3594422,2.0000052)
\psline[linewidth=0.02cm,linecolor=color1365](1.869233,1.1553239)(2.3790238,2.0000052)
\psline[linewidth=0.02cm,linecolor=color1365](1.869233,1.1553239)(1.869233,2.0000052)
\psline[linewidth=0.02cm](0.50979084,1.1553239)(0.8496514,2.0000052)
\psline[linewidth=0.02cm](0.50979084,1.1553239)(0.16993028,2.0000052)
\psline[linewidth=0.02cm,linecolor=color1365](3.2286754,1.1553239)(2.8888147,2.0000052)
\psline[linewidth=0.02cm,linecolor=color1365](3.2286754,1.1553239)(3.2286754,2.0000052)
\psline[linewidth=0.02cm,linecolor=color1365](3.2286754,1.1553239)(3.5685358,2.0000052)
\psline[linewidth=0.02cm](0.16993028,2.0000052)(0.16993028,2.8446867)
\psline[linewidth=0.02cm,linecolor=color1365](0.8496514,2.0000052)(0.6797211,2.8446867)
\psline[linewidth=0.02cm,linecolor=color1365](0.8496514,2.0000052)(1.0195817,2.8446867)
\psline[linewidth=0.02cm,linecolor=color1365](1.3594422,2.0000052)(1.3594422,2.8446867)
\psline[linewidth=0.02cm,linecolor=color1365](1.869233,2.0000052)(1.6993028,2.8446867)
\psline[linewidth=0.02cm,linecolor=color1365](1.869233,2.0000052)(2.0391634,2.8446867)
\psline[linewidth=0.02cm,linecolor=color1365](2.3790238,2.0000052)(2.3790238,2.8446867)
\psline[linewidth=0.02cm,linecolor=color1365](2.3790238,2.0000052)(2.2090936,2.8446867)
\psline[linewidth=0.02cm,linecolor=color1365](2.3790238,2.0000052)(2.5489542,2.8446867)
\psline[linewidth=0.02cm,linecolor=color1365](2.8888147,2.0000052)(2.8888147,2.8446867)
\psline[linewidth=0.02cm,linecolor=color1365](3.2286754,2.0000052)(3.3986056,2.8446867)
\psline[linewidth=0.02cm,linecolor=color1365](3.2286754,2.0000052)(3.2286754,2.8446867)
\psline[linewidth=0.02cm,linecolor=color1365](3.2286754,2.0000052)(3.058745,2.8446867)
\psline[linewidth=0.02cm,linecolor=color1365](3.5685358,2.0000052)(3.5685358,2.8446867)
\psline[linewidth=0.02cm,linecolor=color1365](3.5685358,2.0000052)(3.738466,2.8446867)
\psline[linewidth=0.02cm](0.16993028,2.8446867)(0.0,3.6893682)
\psline[linewidth=0.02cm](0.16993028,2.8446867)(0.16993028,3.6893682)
\psline[linewidth=0.02cm](0.16993028,2.8446867)(0.33986056,3.6893682)
\psline[linewidth=0.02cm,linecolor=color1365](0.6797211,2.8446867)(0.50979084,3.6893682)
\psline[linewidth=0.02cm,linecolor=color1365](0.6797211,2.8446867)(0.8496514,3.6893682)
\psline[linewidth=0.02cm,linecolor=color1365](1.0195817,2.8446867)(1.0195817,3.6893682)
\psline[linewidth=0.02cm,linecolor=color1365](1.3594422,2.8446867)(1.1895119,3.6893682)
\psline[linewidth=0.02cm,linecolor=color1365](1.3594422,2.8446867)(1.3594422,3.6893682)
\psline[linewidth=0.02cm,linecolor=color1365](1.3594422,2.8446867)(1.5293725,3.6893682)
\psline[linewidth=0.02cm,linecolor=color1365](1.6993028,2.8446867)(1.6993028,3.6893682)
\psline[linewidth=0.02cm,linecolor=color1365](2.0391634,2.8446867)(1.869233,3.6893682)
\psline[linewidth=0.02cm,linecolor=color1365](2.0391634,2.8446867)(2.0391634,3.6893682)
\psline[linewidth=0.02cm,linecolor=color1365](2.2090936,2.8446867)(2.2090936,3.6893682)
\psline[linewidth=0.02cm,linecolor=color1365](2.3790238,2.8446867)(2.5489542,3.6893682)
\psline[linewidth=0.02cm,linecolor=color1365](2.3790238,2.8446867)(2.3790238,3.6893682)
\psline[linewidth=0.02cm,linecolor=color1365](2.5489542,2.8446867)(2.7188845,3.6893682)
\psline[linewidth=0.02cm,linecolor=color1365](2.8888147,2.8446867)(2.8888147,3.6893682)
\psline[linewidth=0.02cm,linecolor=color1365](3.058745,2.8446867)(3.058745,3.6893682)
\psline[linewidth=0.02cm,linecolor=color1365](3.2286754,2.8446867)(3.3986056,3.6893682)
\psline[linewidth=0.02cm,linecolor=color1365](3.2286754,2.8446867)(3.2286754,3.6893682)
\psline[linewidth=0.02cm,linecolor=color1365](3.3986056,2.8446867)(3.5685358,3.6893682)
\psline[linewidth=0.02cm,linecolor=color1365](3.5685358,2.8446867)(3.738466,3.6893682)
\psline[linewidth=0.02cm,linecolor=color1365](3.738466,2.8446867)(4.0783267,3.6893682)
\psline[linewidth=0.02cm,linecolor=color1365](3.738466,2.8446867)(3.9083962,3.6893682)
\psline[linewidth=0.02cm,linecolor=color1310](1.869233,-3.6893578)(0.50979084,-2.8446763)
\psline[linewidth=0.02cm,linecolor=color1310](1.869233,-3.6893578)(1.869233,-2.8446763)
\psline[linewidth=0.02cm,linecolor=color1310](1.869233,-3.6893578)(3.2286754,-2.8446763)
\psline[linewidth=0.02cm,linecolor=color1365](1.869233,-2.8446763)(1.3594422,-1.9999946)
\psline[linewidth=0.02cm,linecolor=color1365](1.869233,-2.8446763)(2.3790238,-1.9999946)
\psline[linewidth=0.02cm,linecolor=color1365](1.869233,-2.8446763)(1.869233,-1.9999946)
\psline[linewidth=0.02cm,linecolor=color319](0.50979084,-2.8446763)(0.8496514,-1.9999946)
\psline[linewidth=0.02cm,linecolor=color1310](0.50979084,-2.8446763)(0.16993028,-1.9999946)
\psline[linewidth=0.02cm](3.2286754,-2.8446763)(2.8888147,-1.9999946)
\psline[linewidth=0.02cm](3.2286754,-2.8446763)(3.2286754,-1.9999946)
\psline[linewidth=0.02cm](3.2286754,-2.8446763)(3.5685358,-1.9999946)
\psline[linewidth=0.02cm,linecolor=color1310](0.16993028,-1.9999946)(0.16993028,-1.1553131)
\psline[linewidth=0.02cm,linecolor=color1365](0.8496514,-1.9999946)(0.6797211,-1.1553131)
\psline[linewidth=0.02cm,linecolor=color1365](0.8496514,-1.9999946)(1.0195817,-1.1553131)
\psline[linewidth=0.02cm,linecolor=color1365](1.3594422,-1.9999946)(1.3594422,-1.1553131)
\psline[linewidth=0.02cm,linecolor=color1365](1.869233,-1.9999946)(1.6993028,-1.1553131)
\psline[linewidth=0.02cm,linecolor=color1365](1.869233,-1.9999946)(2.0391634,-1.1553131)
\psline[linewidth=0.02cm,linecolor=color1365](2.3790238,-1.9999946)(2.3790238,-1.1553131)
\psline[linewidth=0.02cm,linecolor=color1365](2.3790238,-1.9999946)(2.2090936,-1.1553131)
\psline[linewidth=0.02cm,linecolor=color1365](2.3790238,-1.9999946)(2.5489542,-1.1553131)
\psline[linewidth=0.02cm](2.8888147,-1.9999946)(2.8888147,-1.1553131)
\psline[linewidth=0.02cm,linecolor=color1365](3.2286754,-1.9999946)(3.3986056,-1.1553131)
\psline[linewidth=0.02cm,linecolor=color1365](3.2286754,-1.9999946)(3.2286754,-1.1553131)
\psline[linewidth=0.02cm,linecolor=color1365](3.2286754,-1.9999946)(3.058745,-1.1553131)
\psline[linewidth=0.02cm](3.5685358,-1.9999946)(3.5685358,-1.1553131)
\psline[linewidth=0.02cm](3.5685358,-1.9999946)(3.738466,-1.1553131)
\psline[linewidth=0.02cm,linecolor=color1310](0.16993028,-1.1553131)(0.0,-0.3106317)
\psline[linewidth=0.02cm,linecolor=color1310](0.16993028,-1.1553131)(0.16993028,-0.3106317)
\psline[linewidth=0.02cm,linecolor=color1310](0.16993028,-1.1553131)(0.33986056,-0.3106317)
\psline[linewidth=0.02cm,linecolor=color1365](0.6797211,-1.1553131)(0.50979084,-0.3106317)
\psline[linewidth=0.02cm,linecolor=color1365](0.6797211,-1.1553131)(0.8496514,-0.3106317)
\psline[linewidth=0.02cm,linecolor=color1365](1.0195817,-1.1553131)(1.0195817,-0.3106317)
\psline[linewidth=0.02cm,linecolor=color1365](1.3594422,-1.1553131)(1.1895119,-0.3106317)
\psline[linewidth=0.02cm,linecolor=color1365](1.3594422,-1.1553131)(1.3594422,-0.3106317)
\psline[linewidth=0.02cm,linecolor=color1365](1.3594422,-1.1553131)(1.5293725,-0.3106317)
\psline[linewidth=0.02cm,linecolor=color1365](1.6993028,-1.1553131)(1.6993028,-0.3106317)
\psline[linewidth=0.02cm,linecolor=color1365](2.0391634,-1.1553131)(1.869233,-0.3106317)
\psline[linewidth=0.02cm,linecolor=color1365](2.0391634,-1.1553131)(2.0391634,-0.3106317)
\psline[linewidth=0.02cm,linecolor=color1365](2.2090936,-1.1553131)(2.2090936,-0.3106317)
\psline[linewidth=0.02cm,linecolor=color1365](2.3790238,-1.1553131)(2.5489542,-0.3106317)
\psline[linewidth=0.02cm,linecolor=color1365](2.3790238,-1.1553131)(2.3790238,-0.3106317)
\psline[linewidth=0.02cm,linecolor=color1365](2.5489542,-1.1553131)(2.7188845,-0.3106317)
\psline[linewidth=0.02cm,linecolor=color1365](2.8888147,-1.1553131)(2.8888147,-0.3106317)
\psline[linewidth=0.02cm,linecolor=color1365](3.058745,-1.1553131)(3.058745,-0.3106317)
\psline[linewidth=0.02cm,linecolor=color1365](3.2286754,-1.1553131)(3.3986056,-0.3106317)
\psline[linewidth=0.02cm,linecolor=color1365](3.2286754,-1.1553131)(3.2286754,-0.3106317)
\psline[linewidth=0.02cm,linecolor=color1365](3.3986056,-1.1553131)(3.5685358,-0.3106317)
\psline[linewidth=0.02cm,linecolor=color1365](3.5685358,-1.1553131)(3.738466,-0.3106317)
\psline[linewidth=0.02cm](3.738466,-1.1553131)(4.0783267,-0.3106317)
\psline[linewidth=0.02cm](3.738466,-1.1553131)(3.9083962,-0.3106317)
\psline[linewidth=0.02cm,fillcolor=color319](10.289352,4.3106318)(5.8893523,4.3106318)
\psline[linewidth=0.02cm,fillcolor=color319](10.289352,4.3106318)(10.689352,4.900632)
\psline[linewidth=0.02cm,fillcolor=color319](10.289352,4.3106318)(10.289352,4.900632)
\psline[linewidth=0.02cm,fillcolor=color319](9.489352,4.3106318)(9.889353,4.900632)
\psline[linewidth=0.02cm,fillcolor=color319](9.489352,4.3106318)(9.489352,4.900632)
\psline[linewidth=0.02cm,fillcolor=color319](9.489352,4.3106318)(9.089353,4.900632)
\psline[linewidth=0.02cm,fillcolor=color319](8.689352,4.3106318)(8.889353,4.900632)
\psline[linewidth=0.02cm,fillcolor=color319](8.689352,4.3106318)(8.489352,4.900632)
\psline[linewidth=0.02cm,fillcolor=color319](7.8893523,4.3106318)(7.8893523,4.900632)
\psline[linewidth=0.02cm,fillcolor=color319](7.0893526,4.3106318)(6.6893525,4.900632)
\psline[linewidth=0.02cm,fillcolor=color319](7.0893526,4.3106318)(7.0893526,4.900632)
\psline[linewidth=0.02cm,fillcolor=color319](7.0893526,4.3106318)(7.4893527,4.900632)
\psline[linewidth=0.02cm,fillcolor=color319](6.0893526,4.3106318)(6.0893526,4.900632)
\psline[linewidth=0.02cm,fillcolor=color319,linestyle=dotted,dotsep=0.16cm](5.8893523,4.3106318)(5.2893524,4.3106318)
\usefont{T1}{ptm}{m}{n}
\rput(12.289352,6.205632){$\tilde{U}_{0}$}
\usefont{T1}{ptm}{m}{n}
\rput(4.7593527,6.205632){$U_0$}
\usefont{T1}{ptm}{m}{n}
\rput(4.7593527,2.205632){$U_1$}
\psdots[dotsize=0.12](9.489352,4.900632)
\psline[linewidth=0.02cm,fillcolor=color319,arrowsize=0.05291667cm 2.0,arrowlength=1.4,arrowinset=0.4]{->}(10.089353,5.500632)(9.489352,4.900632)
\usefont{T1}{ptm}{m}{n}
\rput(10.599353,5.600632){\footnotesize ${\tilde{X}}_{\tau_1}$}
\psline[linewidth=0.02cm,fillcolor=color319](10.489352,0.3106317)(6.0893526,0.3106317)
\psline[linewidth=0.02cm,fillcolor=color319](10.489352,0.3106317)(10.889353,0.900632)
\psline[linewidth=0.02cm,fillcolor=color319](10.489352,0.3106317)(10.489352,0.900632)
\psline[linewidth=0.02cm,fillcolor=color319](9.689352,0.3106317)(10.089353,0.900632)
\psline[linewidth=0.02cm,fillcolor=color319](9.689352,0.3106317)(9.689352,0.900632)
\psline[linewidth=0.02cm,fillcolor=color319](9.689352,0.3106317)(9.289352,0.900632)
\psline[linewidth=0.02cm,fillcolor=color319](8.889353,0.3106317)(9.089353,0.900632)
\psline[linewidth=0.02cm,fillcolor=color319](8.889353,0.3106317)(8.689352,0.900632)
\psline[linewidth=0.02cm,fillcolor=color319](8.089353,0.3106317)(8.089353,0.900632)
\psline[linewidth=0.02cm,fillcolor=color319](7.2893524,0.3106317)(6.8893523,0.900632)
\psline[linewidth=0.02cm,fillcolor=color319](7.2893524,0.3106317)(7.2893524,0.900632)
\psline[linewidth=0.02cm,fillcolor=color319](7.2893524,0.3106317)(7.6893525,0.900632)
\psline[linewidth=0.02cm,fillcolor=color319](6.2893524,0.3106317)(6.2893524,0.900632)
\psline[linewidth=0.02cm,fillcolor=color319,linestyle=dotted,dotsep=0.16cm](6.0893526,0.3106317)(5.4893527,0.3106317)
\psline[linewidth=0.02cm,fillcolor=color319](9.689352,0.900632)(10.089353,1.500632)
\psline[linewidth=0.02cm,fillcolor=color319](9.689352,0.900632)(9.689352,1.500632)
\psline[linewidth=0.02cm,fillcolor=color319](9.689352,0.900632)(9.289352,1.500632)
\psline[linewidth=0.02cm,fillcolor=color319](9.289352,1.500632)(9.089353,2.100632)
\psline[linewidth=0.02cm,fillcolor=color319](9.289352,1.500632)(9.489352,2.100632)
\psline[linewidth=0.02cm,fillcolor=color319](9.089353,2.100632)(9.089353,2.700632)
\psline[linewidth=0.02cm,fillcolor=color319](9.089353,2.700632)(8.889353,3.300632)
\psline[linewidth=0.02cm,fillcolor=color319](9.089353,2.700632)(9.089353,3.300632)
\psline[linewidth=0.02cm,fillcolor=color319](9.089353,2.700632)(9.289352,3.300632)
\usefont{T1}{ptm}{m}{n}
\rput(12.489352,2.205632){$\tilde{U}_{1}$}
\psline[linewidth=0.02cm,fillcolor=color319](10.489352,-3.6893682)(6.0893526,-3.6893682)
\psline[linewidth=0.02cm,fillcolor=color319](10.489352,-3.6893682)(10.889353,-3.0993679)
\psline[linewidth=0.02cm,fillcolor=color319](10.489352,-3.6893682)(10.489352,-3.0993679)
\psline[linewidth=0.02cm,fillcolor=color319](9.689352,-3.6893682)(10.089353,-3.0993679)
\psline[linewidth=0.02cm,fillcolor=color319](9.689352,-3.6893682)(9.689352,-3.0993679)
\psline[linewidth=0.02cm,fillcolor=color319](9.689352,-3.6893682)(9.289352,-3.0993679)
\psline[linewidth=0.02cm,fillcolor=color319](8.889353,-3.6893682)(9.089353,-3.0993679)
\psline[linewidth=0.02cm,fillcolor=color319](8.889353,-3.6893682)(8.689352,-3.0993679)
\psline[linewidth=0.02cm,fillcolor=color319](8.089353,-3.6893682)(8.089353,-3.0993679)
\psline[linewidth=0.02cm,fillcolor=color319](7.2893524,-3.6893682)(6.8893523,-3.0993679)
\psline[linewidth=0.02cm,fillcolor=color319](7.2893524,-3.6893682)(7.2893524,-3.0993679)
\psline[linewidth=0.02cm,fillcolor=color319](7.2893524,-3.6893682)(7.6893525,-3.0993679)
\psline[linewidth=0.02cm,fillcolor=color319](6.2893524,-3.6893682)(6.2893524,-3.0993679)
\psline[linewidth=0.02cm,fillcolor=color319,linestyle=dotted,dotsep=0.16cm](6.0893526,-3.6893682)(5.4893527,-3.6893682)
\psline[linewidth=0.02cm,fillcolor=color319](9.689352,-3.0993679)(10.089353,-2.499368)
\psline[linewidth=0.02cm,fillcolor=color319](9.689352,-3.0993679)(9.689352,-2.499368)
\psline[linewidth=0.02cm,fillcolor=color319](9.689352,-3.0993679)(9.289352,-2.499368)
\psline[linewidth=0.02cm,fillcolor=color319](9.289352,-2.499368)(9.089353,-1.8993679)
\psline[linewidth=0.02cm,fillcolor=color319](9.289352,-2.499368)(9.489352,-1.8993679)
\psline[linewidth=0.02cm,fillcolor=color319](9.089353,-1.8993679)(9.089353,-1.299368)
\psline[linewidth=0.02cm,fillcolor=color319](9.089353,-1.299368)(8.889353,-0.699368)
\psline[linewidth=0.02cm,fillcolor=color319](9.089353,-1.299368)(9.089353,-0.699368)
\psline[linewidth=0.02cm,fillcolor=color319](9.089353,-1.299368)(9.289352,-0.699368)
\usefont{T1}{ptm}{m}{n}
\rput(12.489352,-1.794368){$\tilde{U}_{2}$}
\psline[linewidth=0.02cm,fillcolor=color319,arrowsize=0.05291667cm 2.0,arrowlength=1.4,arrowinset=0.4]{->}(10.689352,2.100632)(10.089353,1.500632)
\usefont{T1}{ptm}{m}{n}
\rput(11.199352,2.200632){\footnotesize $\tilde{X}_{\tau_2}$}
\psdots[dotsize=0.12](10.089353,1.500632)
\psline[linewidth=0.02cm,fillcolor=black](10.089353,-2.499368)(9.889353,-1.8993679)
\psline[linewidth=0.02cm,fillcolor=black](10.089353,-2.499368)(10.089353,-1.8993679)
\psline[linewidth=0.02cm,fillcolor=black](10.089353,-2.499368)(10.289352,-1.8993679)
\psline[linewidth=0.02cm,fillcolor=black](9.889353,-1.8993679)(9.889353,-1.299368)
\psline[linewidth=0.02cm,fillcolor=black](10.289352,-1.8993679)(10.289352,-1.299368)
\psline[linewidth=0.02cm,fillcolor=black](10.289352,-1.8993679)(10.489352,-1.299368)
\psline[linewidth=0.02cm,fillcolor=black](10.489352,-1.299368)(10.489352,-0.699368)
\psline[linewidth=0.02cm,fillcolor=black](10.489352,-1.299368)(10.689352,-0.699368)
\psline[linewidth=0.02cm,linecolor=color1310](1.869233,-7.6893578)(0.50979084,-6.844676)
\psline[linewidth=0.02cm,linecolor=color1310](1.869233,-7.6893578)(1.869233,-6.844676)
\psline[linewidth=0.02cm,linecolor=color1310](1.869233,-7.6893578)(3.2286754,-6.844676)
\psline[linewidth=0.02cm](1.869233,-6.844676)(1.3594422,-5.9999948)
\psline[linewidth=0.02cm,linecolor=color1365](1.869233,-6.844676)(2.3790238,-5.9999948)
\psline[linewidth=0.02cm](1.869233,-6.844676)(1.869233,-5.9999948)
\psline[linewidth=0.02cm,linecolor=color1310](0.50979084,-6.844676)(0.8496514,-5.9999948)
\psline[linewidth=0.02cm,linecolor=color1310](0.50979084,-6.844676)(0.16993028,-5.9999948)
\psline[linewidth=0.02cm,linecolor=color1310](3.2286754,-6.844676)(2.8888147,-5.9999948)
\psline[linewidth=0.02cm,linecolor=color1310](3.2286754,-6.844676)(3.2286754,-5.9999948)
\psline[linewidth=0.02cm,linecolor=color1310](3.2286754,-6.844676)(3.5685358,-5.9999948)
\psline[linewidth=0.02cm,linecolor=color1310](0.16993028,-5.9999948)(0.16993028,-5.155313)
\psline[linewidth=0.02cm,linecolor=color1365](0.8496514,-5.9999948)(0.6797211,-5.155313)
\psline[linewidth=0.02cm,linecolor=color1365](0.8496514,-5.9999948)(1.0195817,-5.155313)
\psline[linewidth=0.02cm](1.3594422,-5.9999948)(1.3594422,-5.155313)
\psline[linewidth=0.02cm,linecolor=color1365](1.869233,-5.9999948)(1.6993028,-5.155313)
\psline[linewidth=0.02cm,linecolor=color1365](1.869233,-5.9999948)(2.0391634,-5.155313)
\psline[linewidth=0.02cm,linecolor=color1365](2.3790238,-5.9999948)(2.3790238,-5.155313)
\psline[linewidth=0.02cm,linecolor=color1365](2.3790238,-5.9999948)(2.2090936,-5.155313)
\psline[linewidth=0.02cm,linecolor=color1365](2.3790238,-5.9999948)(2.5489542,-5.155313)
\psline[linewidth=0.02cm,linecolor=color1310](2.8888147,-5.9999948)(2.8888147,-5.155313)
\psline[linewidth=0.02cm,linecolor=color1365](3.2286754,-5.9999948)(3.3986056,-5.155313)
\psline[linewidth=0.02cm,linecolor=color1365](3.2286754,-5.9999948)(3.2286754,-5.155313)
\psline[linewidth=0.02cm,linecolor=color1365](3.2286754,-5.9999948)(3.058745,-5.155313)
\psline[linewidth=0.02cm,linecolor=color1310](3.5685358,-5.9999948)(3.5685358,-5.155313)
\psline[linewidth=0.02cm,linecolor=color1310](3.5685358,-5.9999948)(3.738466,-5.155313)
\psline[linewidth=0.02cm,linecolor=color1310](0.16993028,-5.155313)(0.0,-4.3106318)
\psline[linewidth=0.02cm,linecolor=color1310](0.16993028,-5.155313)(0.16993028,-4.3106318)
\psline[linewidth=0.02cm,linecolor=color1310](0.16993028,-5.155313)(0.33986056,-4.3106318)
\psline[linewidth=0.02cm,linecolor=color1365](0.6797211,-5.155313)(0.50979084,-4.3106318)
\psline[linewidth=0.02cm,linecolor=color1365](0.6797211,-5.155313)(0.8496514,-4.3106318)
\psline[linewidth=0.02cm,linecolor=color1365](1.0195817,-5.155313)(1.0195817,-4.3106318)
\psline[linewidth=0.02cm](1.3594422,-5.155313)(1.1895119,-4.3106318)
\psline[linewidth=0.02cm](1.3594422,-5.155313)(1.3594422,-4.3106318)
\psline[linewidth=0.02cm](1.3594422,-5.155313)(1.5293725,-4.3106318)
\psline[linewidth=0.02cm,linecolor=color1365](1.6993028,-5.155313)(1.6993028,-4.3106318)
\psline[linewidth=0.02cm,linecolor=color1365](2.0391634,-5.155313)(1.869233,-4.3106318)
\psline[linewidth=0.02cm,linecolor=color1365](2.0391634,-5.155313)(2.0391634,-4.3106318)
\psline[linewidth=0.02cm,linecolor=color1365](2.2090936,-5.155313)(2.2090936,-4.3106318)
\psline[linewidth=0.02cm,linecolor=color1365](2.3790238,-5.155313)(2.5489542,-4.3106318)
\psline[linewidth=0.02cm,linecolor=color1365](2.3790238,-5.155313)(2.3790238,-4.3106318)
\psline[linewidth=0.02cm,linecolor=color1365](2.5489542,-5.155313)(2.7188845,-4.3106318)
\psline[linewidth=0.02cm,linecolor=color1365](2.8888147,-5.155313)(2.8888147,-4.3106318)
\psline[linewidth=0.02cm,linecolor=color1365](3.058745,-5.155313)(3.058745,-4.3106318)
\psline[linewidth=0.02cm,linecolor=color1365](3.2286754,-5.155313)(3.3986056,-4.3106318)
\psline[linewidth=0.02cm,linecolor=color1365](3.2286754,-5.155313)(3.2286754,-4.3106318)
\psline[linewidth=0.02cm,linecolor=color1365](3.3986056,-5.155313)(3.5685358,-4.3106318)
\psline[linewidth=0.02cm,linecolor=color1365](3.5685358,-5.155313)(3.738466,-4.3106318)
\psline[linewidth=0.02cm,linecolor=color1310](3.738466,-5.155313)(4.0783267,-4.3106318)
\psline[linewidth=0.02cm,linecolor=color1310](3.738466,-5.155313)(3.9083962,-4.3106318)
\psline[linewidth=0.02cm,fillcolor=color319,arrowsize=0.05291667cm 2.0,arrowlength=1.4,arrowinset=0.4]{->}(8.689352,-2.499368)(8.089353,-3.0993679)
\usefont{T1}{ptm}{m}{n}
\rput(9.009352,-2.399368){\footnotesize $\tilde{X}_{\tau_3}$}
\psdots[dotsize=0.12](8.089353,-3.0993679)
\usefont{T1}{ptm}{m}{n}
\rput(4.7593527,-1.794368){$U_2$}
\psline[linewidth=0.02cm,fillcolor=color319](10.489352,-7.6893682)(6.0893526,-7.6893682)
\psline[linewidth=0.02cm,fillcolor=color319](10.489352,-7.6893682)(10.889353,-7.099368)
\psline[linewidth=0.02cm,fillcolor=color319](10.489352,-7.6893682)(10.489352,-7.099368)
\psline[linewidth=0.02cm,fillcolor=color319](9.689352,-7.6893682)(10.089353,-7.099368)
\psline[linewidth=0.02cm,fillcolor=color319](9.689352,-7.6893682)(9.689352,-7.099368)
\psline[linewidth=0.02cm,fillcolor=color319](9.689352,-7.6893682)(9.289352,-7.099368)
\psline[linewidth=0.02cm,fillcolor=color319](8.889353,-7.6893682)(9.089353,-7.099368)
\psline[linewidth=0.02cm,fillcolor=color319](8.889353,-7.6893682)(8.689352,-7.099368)
\psline[linewidth=0.02cm,fillcolor=color319](8.089353,-7.6893682)(8.089353,-7.099368)
\psline[linewidth=0.02cm,fillcolor=color319](7.2893524,-7.6893682)(6.8893523,-7.099368)
\psline[linewidth=0.02cm,fillcolor=color319](7.2893524,-7.6893682)(7.2893524,-7.099368)
\psline[linewidth=0.02cm,fillcolor=color319](7.2893524,-7.6893682)(7.6893525,-7.099368)
\psline[linewidth=0.02cm,fillcolor=color319](6.2893524,-7.6893682)(6.2893524,-7.099368)
\psline[linewidth=0.02cm,fillcolor=color319,linestyle=dotted,dotsep=0.16cm](6.0893526,-7.6893682)(5.4893527,-7.6893682)
\psline[linewidth=0.02cm,fillcolor=color319](9.689352,-7.099368)(10.089353,-6.499368)
\psline[linewidth=0.02cm,fillcolor=color319](9.689352,-7.099368)(9.689352,-6.499368)
\psline[linewidth=0.02cm,fillcolor=color319](9.689352,-7.099368)(9.289352,-6.499368)
\psline[linewidth=0.02cm,fillcolor=color319](9.289352,-6.499368)(9.089353,-5.899368)
\psline[linewidth=0.02cm,fillcolor=color319](9.289352,-6.499368)(9.489352,-5.899368)
\psline[linewidth=0.02cm,fillcolor=color319](9.089353,-5.899368)(9.089353,-5.299368)
\psline[linewidth=0.02cm,fillcolor=color319](9.089353,-5.299368)(8.889353,-4.699368)
\psline[linewidth=0.02cm,fillcolor=color319](9.089353,-5.299368)(9.089353,-4.699368)
\psline[linewidth=0.02cm,fillcolor=color319](9.089353,-5.299368)(9.289352,-4.699368)
\usefont{T1}{ptm}{m}{n}
\rput(12.489352,-5.794368){$\tilde{U}_{3}$}
\psline[linewidth=0.02cm,fillcolor=black](10.089353,-6.499368)(9.889353,-5.899368)
\psline[linewidth=0.02cm,fillcolor=black](10.089353,-6.499368)(10.089353,-5.899368)
\psline[linewidth=0.02cm,fillcolor=black](10.089353,-6.499368)(10.289352,-5.899368)
\psline[linewidth=0.02cm,fillcolor=black](9.889353,-5.899368)(9.889353,-5.299368)
\psline[linewidth=0.02cm,fillcolor=black](10.289352,-5.899368)(10.289352,-5.299368)
\psline[linewidth=0.02cm,fillcolor=black](10.289352,-5.899368)(10.489352,-5.299368)
\psline[linewidth=0.02cm,fillcolor=black](10.489352,-5.299368)(10.489352,-4.699368)
\psline[linewidth=0.02cm,fillcolor=black](10.489352,-5.299368)(10.689352,-4.699368)
\psline[linewidth=0.02cm,fillcolor=black](8.089353,-7.099368)(8.089353,-6.499368)
\psline[linewidth=0.02cm,fillcolor=black](8.089353,-7.099368)(7.8893523,-6.499368)
\psline[linewidth=0.02cm,fillcolor=black](7.8893523,-6.499368)(7.8893523,-5.899368)
\psline[linewidth=0.02cm,fillcolor=black](7.8893523,-5.899368)(7.6893525,-5.299368)
\psline[linewidth=0.02cm,fillcolor=black](7.8893523,-5.899368)(7.8893523,-5.299368)
\psline[linewidth=0.02cm,fillcolor=black](7.8893523,-5.899368)(8.089353,-5.299368)
\usefont{T1}{ptm}{m}{n}
\rput(4.7793527,-5.794368){$U_3$}
\psline[linewidth=0.01cm,fillcolor=black,arrowsize=0.05291667cm 2.0,arrowlength=1.4,arrowinset=0.4]{->}(1.6893525,4.300632)(0.5293525,5.020632)
\psline[linewidth=0.01cm,fillcolor=black,arrowsize=0.05291667cm 2.0,arrowlength=1.4,arrowinset=0.4]{->}(0.38935247,5.180632)(0.08935247,5.900632)
\psline[linewidth=0.01cm,fillcolor=black,arrowsize=0.05291667cm 2.0,arrowlength=1.4,arrowinset=0.4]{->}(0.06935247,6.060632)(0.08935247,6.860632)
\psline[linewidth=0.01cm,fillcolor=black,arrowsize=0.05291667cm 2.0,arrowlength=1.4,arrowinset=0.4]{->}(0.26935247,6.820632)(0.26935247,6.020632)
\psline[linewidth=0.01cm,fillcolor=black,arrowsize=0.05291667cm 2.0,arrowlength=1.4,arrowinset=0.4]{->}(0.28935248,5.980632)(0.5293525,5.360632)
\psline[linewidth=0.01cm,fillcolor=black,arrowsize=0.05291667cm 2.0,arrowlength=1.4,arrowinset=0.4]{->}(0.62935245,5.240632)(1.8293525,4.460632)
\usefont{T1}{ptm}{m}{n}
\rput(0.49935248,4.525632){$X_n$}
\psline[linewidth=0.01cm,fillcolor=black,arrowsize=0.05291667cm 2.0,arrowlength=1.4,arrowinset=0.4]{->}(1.9893525,0.26063204)(3.2693524,1.060632)
\psline[linewidth=0.01cm,fillcolor=black,arrowsize=0.05291667cm 2.0,arrowlength=1.4,arrowinset=0.4]{->}(3.3493524,1.280632)(3.6493526,2.000632)
\psline[linewidth=0.01cm,fillcolor=black,arrowsize=0.05291667cm 2.0,arrowlength=1.4,arrowinset=0.4]{->}(3.6493526,2.060632)(3.8293524,2.860632)
\psline[linewidth=0.01cm,fillcolor=black,arrowsize=0.05291667cm 2.0,arrowlength=1.4,arrowinset=0.4]{->}(3.6893525,2.880632)(3.6093526,2.4206321)
\psline[linewidth=0.01cm,fillcolor=black,arrowsize=0.05291667cm 2.0,arrowlength=1.4,arrowinset=0.4]{->}(3.5093524,2.000632)(3.2693524,1.380632)
\psline[linewidth=0.01cm,fillcolor=black,arrowsize=0.05291667cm 2.0,arrowlength=1.4,arrowinset=0.4]{->}(3.2093525,1.3606321)(2.9693525,1.940632)
\psline[linewidth=0.01cm,fillcolor=black,arrowsize=0.05291667cm 2.0,arrowlength=1.4,arrowinset=0.4]{->}(2.8493524,1.900632)(3.1493526,1.180632)
\psline[linewidth=0.01cm,fillcolor=black,arrowsize=0.05291667cm 2.0,arrowlength=1.4,arrowinset=0.4]{->}(3.0893524,1.180632)(1.9093524,0.42063203)
\psline[linewidth=0.01cm,fillcolor=black,arrowsize=0.05291667cm 2.0,arrowlength=1.4,arrowinset=0.4]{->}(1.8493525,-2.659368)(1.4293525,-1.959368)
\psline[linewidth=0.01cm,fillcolor=black,arrowsize=0.05291667cm 2.0,arrowlength=1.4,arrowinset=0.4]{->}(1.9293525,-3.5993679)(1.9293525,-2.819368)
\psline[linewidth=0.01cm,fillcolor=black,arrowsize=0.05291667cm 2.0,arrowlength=1.4,arrowinset=0.4]{->}(1.4293525,-1.879368)(1.4293525,-1.1393679)
\psline[linewidth=0.01cm,fillcolor=black,arrowsize=0.05291667cm 2.0,arrowlength=1.4,arrowinset=0.4]{->}(1.2893524,-1.1593679)(1.2893524,-1.919368)
\psline[linewidth=0.01cm,fillcolor=black,arrowsize=0.05291667cm 2.0,arrowlength=1.4,arrowinset=0.4]{->}(1.2893524,-2.039368)(1.8093525,-2.879368)
\psline[linewidth=0.01cm,fillcolor=black,arrowsize=0.05291667cm 2.0,arrowlength=1.4,arrowinset=0.4]{->}(1.8093525,-2.899368)(1.8093525,-3.5993679)
\end{pspicture} 
}
\caption{the coupling}
\end{figure}
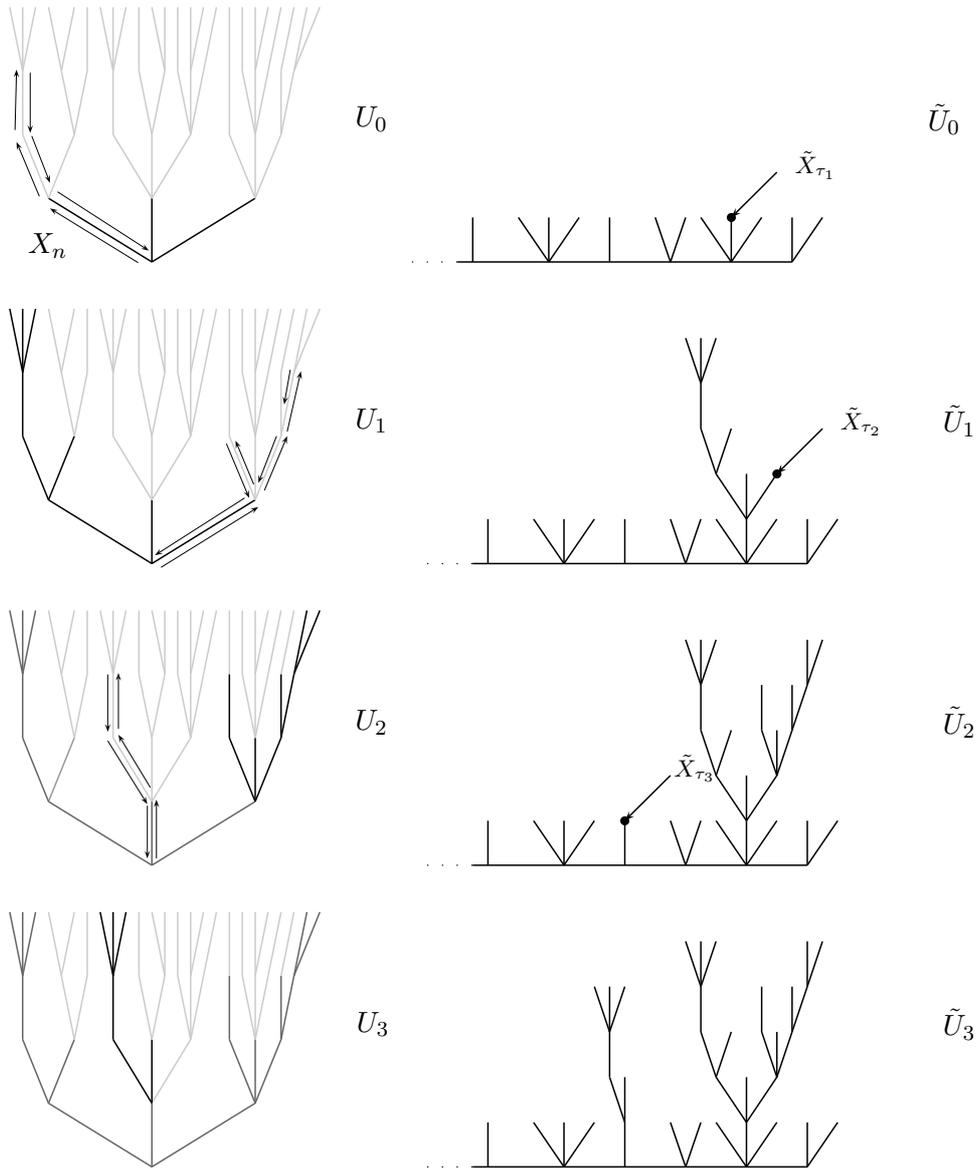
Let us introduce some notations : for $T,S$ two trees, finite or infinite, we set $LT$ the leaves of $T$, that is the vertices of $T$ that have no offspring, $T^o=T/LT$ and for $v\in T$ we denote by $T\circ^{v}S$ the tree obtained by gluing the root of $S$ to the vertex $v$ of $T$, with vertices marked as in their original tree (the vertex coming from both $v$ and the root of $S$ is marked as $v$).
Given a tree $T\in \mathbb{T}$ and a path $\{X_{t}\}$ on $T$ we construct a family of finite trees $T_{i},U_{i}$ as follows : 
let $\tau_{0}=\eta_{0}=0$, and $U_{0}$ the finite tree consisting of the root $e$ of $T$ and its offspring, marked as in $T$.
For $i\geq1,$ let
\begin{align}\tau_{i} & = \min\{t\geq\eta_{i-1}: X_{t}\in LU_{i-1}\} \\ 
\eta_{i}&=\min\{t>\tau_{i};X_{t}\in U_{i-1}^o\}. 
\end{align}
Let $T_{i}$ be the tree ``explored'' by the walk during the excursion
$[\tau_{i},\eta_{i})$, that is to say $T_{i}$ is composed of the
 vertices of $T$ visited by $\{X_{t} ,t\in[\tau_{i},\eta_{i})\}$,
 together with their offspring, marked as in $T$, and the root of $T_i$ is $X_{\tau_{i}}$. Let $U_{i}=U_{i-1}\circ^{X_{\tau_{i}}}T_{i}$ be the tree explored by the walk from the beginning. We call $\{u_{t}^{i}\}_{t=0}^{\eta_{i}-\tau_{i}-1}$ the path in $T_{i}$ defined by $u_{t}^{i}=X_{\tau_{i}+t}$.
If $T$ is distributed according to {\tt MT}, and ${X_{t}}$ is the path of the random walk on $T$, then, the walk being recurrent, $\mathbb{P}_{\tt MT}-$almost surely $T=\lim U_{i}$.
\\

We are now going to construct $\tilde{T}\in\tilde{\mathbb{T}}$, a tree
with a semi-infinite ray emanating from the root, coupled with $T$,
and a path $\{\tilde{X}_{t}\}$ on T, in such a way that, if $T$ is distributed
according to {\tt MT}, and ${X_{t}}$ is the path of the random walk on
$T$, then $\tilde{T}$ will be distributed according to {\tt IMT} and
$\{\tilde{X}_{t}\}$ will follow the law of a random walk on $\tilde{T}$.

Let $\tilde{U}_{o}$ be the tree defined as follows : we choose a vertex denoted by $e$, as the root of $\tilde{U}_{o}$, and a semi-infinite ray $\{e=v_{0},v_{1},...\}.$ To each vertex $v_{i} \in Ray$ different from $e$ we attach independently a set of marked vertices with law
 $\hat{q}$. To $e$ we attach a set of children with distribution $(q+\hat{q})/2$
If $i\geq1$ we chose one of those vertices, with probability $\frac{A(x)}{\sum_{y}
 A(y)}$, and identify it with $v_{i-1}$. We obtain a tree with a semi-infinite ray and a set of children for each vertex $v_{i}$ on $Ray$, one of them being $v_{i-1}$.
 
 We set $\tilde{\tau_{0}}=\tilde{\eta_{0}}=0$. Recalling the relation (\ref{relation}) between the $A_{x}$ and the $\omega(x,y)$, one can easily check that for any vertex $x$, knowing the $\{w(x,y)\}_{{y\in T}}$ is equivalent to knowing $\{A(x_{i})\}_{x_{i}\text{ children of }x}$. Thus, knowing $\tilde{U}_{0}$ one can compute the $\{\omega(x,y)\}_{x\in Ray,y\in \tilde{U}_{0}}$ and define a random walk $\tilde{X}_{t}$ on $\tilde{U}_{0},$ stopped when it gets off ${\it Ray}$. We set accordingly $\tilde{\tau}_{1}=\min\{t>0:\tilde{X}_{t}\in L\tilde{U}_{0}\}.$ 

We are now going to ``glue'' the first excursion of $\{X_{t}\}$.
Let
\begin{align*}&\tilde{U}_{1}=\tilde{U}_{0}\circ^{\tilde{X}_{\tilde{\tau}_{1}}}T_{1},\\
&\tilde{\eta}_{1}=\tilde{\tau}_{1}+\eta_{1}-\tau_{1},\\
& \{\tilde{X}_{t}\}_{t=\tilde{\tau}_{1}}^{\tilde{\eta}_{1}-1}=u_{t-\tilde{\tau}_{1}}^1,\\
&\tilde{X}^{\tilde{\eta_{1}}}=\overleftarrow{\tilde{X}^{\tilde{\eta}_{1}-1}}.\end{align*} One can easily check that $\{\tilde{X}_{t}\}_{t\leq\tilde{\eta}_{1}}$ follows the law of a random walk on $\tilde{U}_{1}$.

We iterate the process, in the following way : for $i>1$, start a random walk $\{\tilde{X}_{t}\}_{t\geq \tilde{\eta}_{i-1}}$ on $\tilde{U}_{i-1}$, and define\begin{align*}
&\tilde{\tau}_{i}=\min\{t>0:\tilde{X}_{t}\in L\tilde{U}_{i-1}\}, \\
&\tilde{U}_{i}=\tilde{U}_{i-1}\circ^{\tilde{X}_{\tilde{\tau}_{i}}}T_{i},\\
&\tilde{\eta}_{i}=\tilde{\tau}_{i}+\eta_{i}-\tau_{i},\\
&\{\tilde{X}_{t}\}_{t=\tilde{\tau}_{i}}^{\tilde{\eta}_{i}-1}=u_{t-\tilde{\tau}_{i}}^i, \\
&\tilde{X}^{\tilde{\eta}_{i}}=\overleftarrow{\tilde{X}^{\tilde{\eta}_{i}-1}}.
\end{align*}

Finally, set $\tilde{U}=\bigcup_0^{\infty} \tilde{U_i}$ and $\tilde{T}$
the tree obtained by attaching independents {\tt MT} trees to each
leaves of $\tilde{U}$. It is a direct consequence of the construction
that
\begin{proposition}
If $T$ is distributed according to {\tt MT} and $X_t$ follows
$\mathbb{P}_T$, then $\tilde{T}$ is distributed according to {\tt IMT},
and $\tilde{X}_t$ follows $\mathbb{P}_{\tilde{T}}$.
\end{proposition}

As a consequence, under proper assumptions on $q$, application of
Proposition \ref{theo} implies that for {\tt MT} almost every $T$ the process $\{h(\tilde{X}_{\lfloor n t\rfloor})/\sqrt{\sigma^2 n}\}$ converges
to a standard Brownian motion, as $n$ goes to infinity.

We introduce $R_t=h(\tilde{X}_t)-\min_{i=1}^t
h(\tilde{X}_i)$. We get immediately that 

\noindent$\{R_{\lfloor
 n t\rfloor}/\sqrt{\sigma^2 n}\}$ converges to a Brownian motion
reflected to its minimum, which has the same law as the absolute value
of a Brownian motion.

In order to prove Theorem \ref{clt}, we need to control the distance
between $R_t$ and $|X_t|$. 
 
Let $I_t=\max\{i:\tau_i\leq t\}$ and
$\tilde{I}_t=\max\{i:\tilde{\tau}_i\leq t\}$ the number of excursions
started by $\{X_t\}$ and $\{\tilde{X}_t\}$ before time $t$. Let
$\Delta_t=\sum_{i=1}^{I_t}(\tau_i-\eta_{i-1})$ and
$\tilde{\Delta}_t=\sum_{i=1}^{\tilde{I}_t}(\tilde{\tau}_i-\tilde{\eta}_{i-1})$,
which measure the time spent by $\{X_t\}$ and $\{\tilde{X}_t\}$ outside the
coupled excursions before time $t$. 
By construction, the distance between $R_t$ and $|X_t|$ comes only
from the parts of the walks outside those excursion. In order to
control these parts, we set for $0\leq\alpha<1/2$ 
$$\Delta_t^{\alpha}=\sum_{i=1}^{I_t}\sum_{s=\eta_{i-1}}^{\tau_i-1}\mathds{1}_{|X_s|\leq
 t^{\alpha}};$$
similarly, 
$$\tilde{\Delta}_t^{\alpha}=\sum_{i=1}^{\tilde{I}_t}\sum_{s=\tilde{\eta}_{i-1}}^{\tilde{\tau}_i-1}\mathds{1}_{d(\tilde{X}_s,Ray)\leq
 t^{\alpha}}.$$
Finally, let $${\bf B}_t=\max_{0\leq r<s\leq t,\tilde{X}_r\in Ray,\tilde{X}_s\in
 Ray}(h(\tilde{X}_s)-h(\tilde{X}_r)),$$ be the maximum amount the walk $\{\tilde{X}_t\}$ moves against the
 drift on $Ray$. 
We have the following 
\begin{proposition}\label{lt}
Under the assumptions of Theorem \ref{clt}, for some $\alpha<1/2$
\begin{equation}\label{lt1}
\lim_{t\rightarrow \infty} \mathbb{P}_T(\Delta_t\neq
\Delta_t^{\alpha})=0, \: {\tt MT}-a.s.,
\end{equation}
and 
\begin{equation}\label{lt2}
\lim_{t\rightarrow \infty} \mathbb{P}_T(\tilde{\Delta}_t\neq
\tilde{\Delta}_t^{\alpha})=0, \: {\tt IMT}-a.s..
\end{equation}
Further, 
\begin{equation}\label{lt3}
\limsup \frac{\Delta_t}{t}=0, \: {\tt MT}-a.s.,
\end{equation}
and
\begin{equation}\label{lt4}
\limsup \frac{\tilde{\Delta}_t}{t}=0, \: {\tt IMT}-a.s..
\end{equation}
Finally,
\begin{equation}\label{lt5}
\limsup \frac{{\bf B}_t}{\sqrt{t}}=0,\: {\tt IMT}-a.s..
\end{equation}
\end{proposition}

Before proving this proposition, note that on the event $\{ \Delta_t=
\Delta_t^{\alpha} \} \cap \{\tilde{\Delta}_t=
\tilde{\Delta}_t^{\alpha}\},$ we have 
$$\min_{s:|s-t|\leq \Delta_t+ \tilde{\Delta}_t}||X_{t}|-R_{s}|\leq 2 t^{\alpha}+{\bf B}_{t}.$$ Therefore we obtain that almost surely, there exists a time change $\theta_t$ such that, for $t$ large enough,
$$\frac{|X_t-R_{\theta_t}|}{\sqrt{t}}\rightarrow_{t\rightarrow \infty} 0$$ and 
$$\frac{|\theta_t-t|}{t}\rightarrow_{t\rightarrow \infty} 0.$$
As we said earlier, Proposition \ref{theo} implies that
$\{R_{\lfloor
 n t\rfloor}/\sqrt{\sigma^2 n}\}$ converges, as $n$ goes to infinity, to the law of the absolute value
of a Brownian motion. Therefore so does $\{R_{\lfloor
 n \theta_t\rfloor}/\sqrt{\sigma^2 n}\}.$ We deduce the result for $|X_t|$.

We now turn to the proof of Lemma \ref{lt}. 
We introduce some notations: for $k\geq 1$, let $a_{k}=\sum_{j=1}^k \tau_{j},$ $b_{k}=\sum_{j=0}^{k-1}\eta_j$ and $J_{k}=[a_{k}-b_{k}+k,a_{k+1}-b_{k+1}+k]$. Note that $\{J_{k}\}_{k\geq1}$ is a partition of 
$\mathbb{N},$ such that the length of $J_{k}$ is equal to the time
spent by the walk between the $k-th$ and the $k+1-th$ excursion. For
$s\in J_{k}$, let ${\bf t}(s)=\eta_{k}+s-(a_{k}-b_{k}+k)$ and
$Y_{0}=0,Y_{1}=X_{\tau_{1}},$ and
$Y_{s}=X_{{\bf t}(s)}$. $\{Y_{s}\}_{s\geq 0}$ is the walk $X_n$
restricted off the excursions, it is clearly not Markovian, nevertheless, it is adapted to the filtration $G_{s}=\sigma(X_{k},k\leq t(s))$. 
For a fixed $t$, we set the sequence $\Theta_{i}$ of stopping times with respect to $G_{s}$ defined by $\Theta_{0}=0$ and 
$$\Theta_{i}=\min\{s>\Theta_{i-1}: \left | |Y_{s}|-|Y_{\Theta_{i-1}}|\right|=\lfloor(\log t)^{3/2}\rfloor\}.$$ 
Similarly, we set, 
for $k\geq1$, $\tilde{a}_{k}=\sum_{j=1}^k \tilde{\tau}_{j},$ $\tilde{b}_{k}=\sum_{j=0}^{k-1}\tilde{\eta}_j$ and $\tilde{J}_{k}=[\tilde{a}_{k}-\tilde{b}_{k}+k,\tilde{a}_{k+1}-\tilde{b}_{k+1}+k]$, and for 
$s\in \tilde{J}_{k}$, we call $\tilde{t}(s)=\tilde{\eta}_{k}+s-(\tilde{a}_{k}-\tilde{b}_{k}+k)$ and
$\tilde{Y}_{0}=0,\tilde{Y}_{1}=\tilde{X}_{\tilde{\tau}_{1}},$ and
$\tilde{Y}_{s}=\tilde{X}_{\tilde{t}(s)}$ the walk $\tilde{X}_n$
restricted off the excursions. We set $\tilde{G}_{s}=\sigma(\tilde{X}_{k},k\leq \tilde{t}(s))$.
For a fixed $t$, we set the sequence of stopping times $\tilde{\Theta}_{i}$ with respect to $\tilde{G}_{s}$ defined by $\tilde{\Theta}_{0}=0$ and 
$$\tilde{\Theta}_{i}=\min\left\{s>\tilde{\Theta}_{i-1}: \left | d(\tilde{Y}_{s},{Ray})-d(\tilde{Y}_{\Theta_{i-1}},Ray)|\right|=\lfloor(\log t)^{3/2}\rfloor\right\}.$$ 
We need the following lemma, whose demonstration will be postponed.
\begin{lemma}\label{lemme}
For all $\epsilon>0$ 
\begin{align}&\label{l1}\lim_{t\rightarrow\infty} \mathbb{P}_{T}\left(\sum_{i=1}^{t^{1/2 +\epsilon}}(\eta_{i}-\tau_{i})<t\right)=0, \: {\tt MT}-a.s.,\\
&\label{l1'}\lim_{t\rightarrow\infty} \mathbb{P}_{T}\left(\sum_{i=1}^{t^{1/2 +\epsilon}}(\tilde{\eta}_{i}-\tilde{\tau}_{i})<t\right)=0, \: {\tt IMT}-a.s.,\\
\label{l2}&\exists
\epsilon'>0:\lim_{t\rightarrow\infty}\mathbb{P}_T\left(\exists s\leq t,
W_{X_s}>t^{1/4-\epsilon'}\right)=0,\: {\tt
 MT}-a.s.,\\
\label{l6}&\text{ and }\lim_{t\rightarrow\infty}\mathbb{P}_T\left(\exists s\leq t,
W_{X_{s}}>t^{1/4-\epsilon'}\right)=0,\: {\tt IMT}-a.s.,\\
\label{l3}&\lim_{t\rightarrow\infty}\mathbb{P}_T(\exists k\leq I_t,
\Theta_{i-1},\Theta_{i}\in J_k,|Y_{\Theta_i}|>|Y_{\Theta_{i-1}}|)=0, \: {\tt MT}-a.s.,\\
\label{l3'}&\lim_{t\rightarrow\infty}\mathbb{P}_T(\exists k\leq I_t,
\tilde{\Theta}_{i-1},\tilde{\Theta_{i}}\in \tilde{J}_k,d(\tilde{Y}_{\tilde{\Theta_i}},Ray)>d(\tilde{Y}_{\tilde{\Theta}_{i-1}},Ray)=0, \: {\tt IMT}-a.s.,\\
&\label{l4}
\lim_{t\rightarrow\infty} \mathbb{P}_{T}(X_s\in \cup_{k=t^{\alpha}-(\log t)^2}^{t^{\alpha}}{\bf A}_{k}^{\epsilon} \text{ for
 some } s\leq t)=0,\, {\tt MT}-a.s.,\\
&\label{l5}
\lim_{t\rightarrow\infty} \mathbb{P}_{T}(X_s\in \cup_{k=t^{\alpha}-(\log t)^2}^{t^{\alpha}}{\bf B}_{k}^{\epsilon} \text{ for
 some } s\leq t)=0, \, {\tt IMT}- a.s..
\end{align}
\end{lemma}
Using this lemma, we can finish the proof of Proposition \ref{lt}.
We shall prove the following statement, which implies (\ref{lt1}) : for some $\alpha\leq1/2$,
\begin{equation}\label{st}
\lim_{t\rightarrow \infty} \mathbb{P}_{T}(\max_{s\in\cup_{k=1}^{I_{t}}J_{k}}|Y_{s}|\geq t^\alpha)=0, \;{\tt MT}-a.s..
\end{equation}
It is a direct consequence of (\ref{l1}) and (\ref{l3}) that, {\tt MT} almost surely, with $\mathbb{P}_T$
probability approaching $1$ as $t$ goes to infinity, $$t(\Theta_{2
 t^{1/2+\epsilon}})>t,$$ whence, using lemma \ref{lemme}, 
\begin{multline*}
\lim_{t\rightarrow \infty}
\mathbb{P}_{T}\big(\max_{s\in\cup_{k=1}^{I_{t}}J_{k}}|Y_{s}|\geq
t^\alpha)\\
\leq \limsup_{t\rightarrow\infty}\sum_{i=0}^{2
 t^{1/2+\epsilon}}\mathbb{P}_{T}\left(\exists j>i:|Y_{\Theta_j}|\geq
t^{\alpha}-(\log t)^2,\,Y_{\Theta_i}=e,\right.\\\left. S_{Y_{\Theta_j}}\geq (\eta-\epsilon_1)
t^{\alpha}/2,\,
|Y_{\Theta_k}|>0, \forall i<k\leq j;\right.\\|S_{X_s}-|X_s||\leq
\epsilon|X_s|, \forall s\leq t \big):=\limsup_{t\rightarrow\infty}\sum_{i=1}^{2
 t^{1/2+\epsilon}}P_{i,t}\,;\\
\end{multline*}
where $\epsilon',\epsilon_1$ are positive numbers that can be
chosen arbitrarily small.

For a fixed $i$ and a fixed $t$, we set $\tilde{M}_s=S_{X_{\Theta_{i+s}}}$, and 
$$
K_t=\min\left\{s>1: X_r=0 \text{ for some }r\in\left[s(\theta_{i+1}),s(\theta_{i+t})\right]\right\}.
$$
The process $\{N_{s}\}=\{ \tilde{M}_{s \wedge K_r}-\tilde{M}_1 \}$ is
a supermartingale with respect to the filtration
$\tilde{G}'_s=\tilde{G}_{\theta_{i+s}}$; indeed as long as the walk
does not come back to the root, the conditional expectation of
$S_{Y_{s+1}}-S_{Y_{s}}$ is lesser or equal to $0$, and by
construction the walk can only return to the root at a $\Theta_i$.

Note that $M_s$ and $N_s$ depend on $t$, whereas this is omitted in
the notation. Let $A_{s}$ be the predictable process such that $N_{s}+A_{s}$ is a martingale. 

Note that, on the event $\{W_{X_{s}}\leq t^{1/4-\epsilon'}, \forall s\leq t \}$ the increments of
$N_s$ are bounded by $t^{1/4-\epsilon'}(\log t)^{3/2}$. One can easily see that the increments of $A_{s}$ are also bounded by $t^{1/4-\epsilon'}(\log t)^{3/2}.$ Therefore
Azuma's Inequality implies
$$P_{i,t}\leq \exp{\left(-t^{2\alpha}/t{1/2+2\epsilon+2(1/4-\epsilon')}\right)}.$$
Recalling that we can choose $\epsilon$ arbitrarily small and
$\alpha$ arbitrarily close to $1/2$, we get the result.\\

The proof of (\ref{lt2}) is quite similar and omitted.
\\

To prove (\ref{lt3}) we introduce \begin{equation}\label{deftepsilon}T^{\epsilon}(t)=\min\{s:|X_s|\geq t^{1/2+\epsilon}\}.\end{equation} By Lemma \ref{lem3.2}, we have
$$\mathbb{P}_{\tt MT}(T_\epsilon(t)<t)\leq t e^{-t^{2\epsilon}}.$$
Using the Borel-Cantelli Lemma, we get that, {\tt MT} almost surely
\begin{equation}\label{ptepsilon}\mathbb{P}_{T}(T_\epsilon(t)<t)\leq e^{-t^{\epsilon}} \text{ for } t>t_0(T).\end{equation}
Let $C_{0,l}$ be the conductance between the root and the level $l$ of the tree. Recalling that for $w$ an offspring of $v$, the conductance associated to the edge $[v,w]$ is $C_w$, Thomson's principle implies that
$$C_{0,l}^{-1}=\inf_{f\text{ unit flow }} \sum_{i=0}^{l} \sum_{v\in T_i} \sum_{w\text{ offspring of v}}\frac{f^2_{v,w}}{C_w}.$$ 
As one can easily check, $f_{v,w}=\frac{C_w W_w}{W_e}$ is a unit flow from the root to $T_l$, so we get
$$C_{0,l}^{-1}\leq\frac{1}{W_e} \sum_{i=1}^l \sum_{v\in T_i}C_w W_w^2.$$
As, conditionally to ${\cal G}_{i}$, $W_w^2$ are independent and identically distributed variables, with finite moment of order two (the assumption needed for that is $\kappa>4$), we have
$$E_{\tt MT}\left[\left( \sum_{v\in T_i}C_w W_w^2- \sum_{v\in T_i}C_w E_{\tt MT}[W_w^2]\right)^2\right]\leq C_{17}\rho(2)^{i},$$ for some constant $C_{17}$, then, using Markov's Inequality, for every $\nu>0$ there exists a constant $C_{18}$ such that
$$P_{\tt MT}\left[ \sum_{v\in T_i}C_w |W_w^2-E[W_{w}^2]|>\nu\right]\leq C_{18}\rho(2)^{i}.$$
This is summable, so by the Borel-Cantelli Lemma, for some constant $C(T)$ dependant only on $T$, we get
$$ \sum_{v\in T_i}C_w W_w^2\leq C(T) \sum_{v\in T_i}C_w.$$
The last part being convergent, thus bounded, we get
\begin{equation}\label{czerol}C_{0,l}^{-1}\leq C(T)l.\end{equation}
If $L_{0}(t)$ denotes the number of visits to the root before time $t$, we get
$$\mathbb{E}_{T}[L_{0}(T_{\epsilon}(t))]=1+C_{0,t^{1/2+\epsilon}}^{-1},$$
indeed $L_{0}(T_{\epsilon}(t))-1$ follows a geometric law with parameter $1- C_{0,t^{1/2+\epsilon}}^{-1}.$

Let $N_{t}(\alpha)=\sum_{k=0}^t \mathds{1}_{|X_{k}|\leq t^\alpha}$ 
On the event that $T_{\epsilon}(t)> t$, we have, using Markov's property,
$$\mathbb{E}_{T}[N_{t}(\alpha);T_{\epsilon}(t)> t]\leq \mathbb{E}_{T}[L_{0}(T_{\epsilon}(t))]\pi\left(\bigcup_{0}^{t^{\alpha}}T_{t}\right)\leq C_{19}(T)t^{1/2+\epsilon+\alpha}.$$ 
Thus as $\mathbb{P}_{T}(T_{\epsilon}(t)\leq t)\leq C_{19}(T) e^{-n^\epsilon},$ using the monotonicity of $N_n(\alpha)$, we obtain $N_{t}(\alpha)/t \rightarrow 0$, from which the result follows, as $\Delta_{t}^\alpha \leq N_{t}^{\alpha}$ and $\mathbb{P}_{T}(\Delta_{t}\neq \Delta_{t}^{\alpha})\rightarrow 0$.
\\

Now we turn to the proof of (\ref{lt4}). 
By the same calculations as in the proof of Lemma \ref{bm}, for $\kappa>5$, we get that 
$\mathbb{E}_{\tt IMT}[\sum_{s\leq t} \mathds{1}_{d(X_{s},{\it Ray})<t^{\alpha }}] \leq t^{1/2+\alpha + \epsilon}$ for any $\epsilon>0$, from which the result follows by an application of Markov's Inequality and the Borel-Cantelli Lemma, using also the fact that the quantity in the expectation is non-decreasing in $n$.
\\

The conductance from $v_k$ to
$v_{k-u}$ is at most $C_{v_{k-u}}$, thus we have the bound
$$\mathbb{P}_{T} ({\bf B}_{t}>u)\leq t \sum_{k=u}^t \Pi_{i=k}^{k-u} A(v_{i}).$$
By Theorem \ref{biggins} and Lemma \ref{changement loi}, the {\tt IMT}-expectation of the right hand side is of order at most $t^2
\rho(2)^{u}$, therefore (\ref{lt5}) follows by standard arguments.
\section{Proof of Lemma \ref{lemme}.}\label{sec6}

It is clear that (\ref{l1}) and (\ref{l1'}) are equivalent. We postpone the proof of these parts to the end of the section.

\noindent \textbf{Proof of (\ref{l2}) :}
 following \cite{peres-zeitouni}, we call ``fresh time'' a time where the walk explore a new vertex, we have
 \begin{eqnarray*}\mathbb{P}_{\tt MT}\left(\exists s\leq t,W_{X_s}>t^{1/4-\epsilon'}\right)&\leq& \sum_{0}^t \mathbb{P}_{\tt MT}[W_{X_{s}}>t^{1/4-\epsilon'};s\text{ is a fresh time}]
 \\&=&\mathbb{P}_{\tt MT}[W_{0}>t^{1/4-\epsilon'}]<C_{20} t/t^{\mu(1/4+\epsilon')},\end{eqnarray*}
for $\mu<\kappa$. If $\kappa>8$, for $\epsilon$ small enough, we can chose $\mu$ such that this is summable. Then the Borel-Cantelli Lemma implies the result. 
\\

{\bf Proof of (\ref{l6})} We are going to use the same arguments, excepted that we have to treat separately the vertices on ${\it Ray}$. More precisely
\begin{multline*}\mathbb{P}_{\tt MT}\left(\exists s\leq t,W_{X_s}>t^{1/4-\epsilon'}\right)\\ \leq \sum_{0}^t \mathbb{P}_{\tt IMT}[W_{X_{s}}>t^{1/4-\epsilon'};s\text{ is a fresh time and }X_s\not\in Ray]
+\\ \mathbb{P}_{\tt MT}\left(\exists s\leq t,W_{v_s}>t^{1/4-\epsilon'}\right).
\end{multline*}
The second term is easily bounded, and the first one is similar to the previous case.
\\

\noindent\textbf{Proof of (\ref{l3}) :}
the event in the probability
in (\ref{l3}) implies that, before time $t$ the walk $X_s$ gets to some
vertex $u$, situated at least at a distance $\lfloor (\log t)^{3/2}\rfloor$,
then back to the ancestor $a(u)$ of $u$ situated at distance $\lfloor
(\log t)^{3/2}\rfloor$ from u, then back again. Decomposing on the
hittings of the root, we can majorate this
probability by
$$\sum_{s\leq t}\mathbb{P}_T(X_t=e)\sum_{k=\lfloor (\log t)^{3/2}\rfloor}^t
\sum_{u\in
 T_k}\mathbb{P}_T(H_u<H_e)\mathbb{P}_T^{a(u)}(H_u<t),$$
where $H_u$ stands for the hitting time of $u$. Using the fact that
the conductance from $0$ to $u$ is bounded by $C_u$, the
probability we are considering is at most
$$n\sum_{k=\lfloor (\log t)^{3/2}\rfloor}^t
\sum_{u\in
 T_k}C_u\mathbb{P}_T^{a(u)}(H_u<t).$$
Denoting by $C(v\rightarrow u)$ the conductance between $v$ and $u$, we
have easily
$$\mathbb{P}_T^{v}(H_u<t)<t\frac{C(v\rightarrow
 u)}{\pi(v)}<c_1 t\frac{C_u}{C_v}.$$
As a direct consequence of Theorem \ref{big}, we have
\begin{multline*}
E_{\tt MT}\left[\sum_{u\in
 T_k}C_u\mathbb{P}_T^{a(u)}(H_u<t)\right] \leq c_1 t^2 E_{\tt MT}\left[\sum_{u\in
 T_k}C_u \frac{C_u}{C_{a(u)}} \right] \\ \leq c_1 t^2\left(E_q\left[\sum
 A_i\exp(\log(A_i))\right]\right)^{\lfloor (\log
 t)^{3/2}\rfloor} \leq c_1 t^2 \rho(2)^{\lfloor (\log
 t)^{3/2}\rfloor}. \end{multline*} 
The result follows by an application of the Borel-Cantelli Lemma.
\\

\noindent\textbf{Proof of (\ref{l3'}) :}
The proof is quite similar to the precedent argument, summing over the different $T^{(v_i)}$.

\noindent\textbf{Proof of (\ref{l4}) :}
using $\kappa > 5,$ by Lemma \ref{cond} we can find an $\varepsilon>0$ such that {\tt IMT}- almost surely the sequence $n^{3/2+\varepsilon}\pi({\bf A}_{n}^\epsilon)$ is summable, thus bounded, so there exists a constant $C'(T)$ such that for each $n$, $C_{e\rightarrow {\bf A}_{n^\epsilon}}\leq C'(T)/n^{3/2+\varepsilon}$.
Recalling from the proof of (\ref{lt3}) the definition of $L_{0}(t),$ and $T_{\epsilon}(t)$ we have
\begin{equation}\mathbb{P}_{T}(X_{t}\in {\bf A}_{t^{\alpha}}^{\epsilon};t\leq T_{\epsilon}(t))\leq \mathbb{E}_{T}[L_{0}(T_{\epsilon}(t))] C'(T)/t^{\alpha(3/2+\varepsilon)} \leq t^{1/2+\epsilon'-\alpha(1+\varepsilon)},\end{equation}
where $\epsilon'$ can be chosen arbitrarily close to $0$. By choosing $\alpha$ close enough to $1/2$, the result follows easily, using (\ref{ptepsilon}).
\\

\noindent\textbf{Proof of (\ref{l5}) :}
we recall from (\ref{defb}) the definition of the sets ${\bf B}_{n}^{\epsilon}.$
By the same argument as in the proof of Lemma \ref{bm}, we get
$$
\lim_{t\rightarrow\infty} \mathbb{P}_{T}(X_s\in {\bf B}_{t^\alpha}^{\epsilon} \text{ for
 some } s\leq t\} 
 \leq H_{\lfloor t^{1/2+\epsilon}\rfloor}\sum_{i=0}^{\lfloor t^{1/2+\epsilon} \rfloor}U_{i}^{t^\alpha},\\
$$
with $H_{t}=1+ \sum_{j=0}^{t-1}\prod_{k=j-1}^{t}A(v_{k}),$ and $U_{i}^{t^\alpha}$ is the probability to get to $B_{t^\alpha}^{\epsilon}$ during one excursion in $T^{v_{i}}$.
By the same argument as in the proof of Lemma \ref{cond}, we get that, almost surely, there exists a constant $C''(T)$ such that
 $$H_{t}\leq C''(T)t^{1/7},$$ whence
$$
\lim_{t\rightarrow\infty} \mathbb{P}_{T}(X_s\in {\bf B}_{t^\alpha}^{\epsilon} \text{ for
 some } s\leq t\} 
 \leq C''(T)t^{1/7}\sum_{i=0}^{\lfloor t^{1/2+\epsilon} \rfloor}U_{i}^{t^\alpha}.$$
Then, denoting $\sum_{t=0}^\infty U_{i}^t t^{1+\epsilon'}:=E_{i},$ the
$E_{i}$ are i.i.d. variables (under {\tt IMT}) with finite expectation
for $\epsilon'$ small enough and
$U_{i}^{t}<\frac{1}{t^{3/2}}E_{i}.$ Then the result follows, using the
law of large numbers.
\\
\noindent{\bf Proof of (\ref{l1}) and (\ref{l1'}):}
Note that, under {\tt MT}, the random variables ${\eta_i-\tau_i}$ are i.i.d.. On the other hand, as a consequence of (\ref{ellipticity}), for some constant $\nu_0>0,$
$$\mathbb{P}_{\tt MT}[{\eta_i-\tau_i}>x]\geq \nu_0 \mathbb{P}_{\tt MT}[T_0>x],$$
Where $T_0$ is the first return to the root.
We recall from (\ref{deftepsilon}) that
\begin{equation*}T^{\epsilon}(t)=\min\{t:|X_s|\geq t^{1/2+\epsilon}\}\end{equation*}
Then, following the proof of Lemma 10 of \cite{peres-zeitouni}, we have,
\begin{equation}\mathbb{P}_{T}[T_0>t]\geq\mathbb{P}_{T}[T_0>T^{\epsilon/2}(t)]\mathbb{P}_{T}[T^{\epsilon/2}(t)\geq t|T_0>T^{\epsilon}(t)].\end{equation}
As a consequence of (\ref{czerol}), for come constant depending on the tree $C_3(T),$
$$\mathbb{P}_{T}[T_0>T^{\epsilon/2}(t)]>C_3(T)t^{-1/2-\epsilon/2}.$$
On the other hand, 
$$\mathbb{P}_{T}[T^{\epsilon/2}(t)< t|T_0>T^{\epsilon/2}(t)]\leq\frac{\mathbb{P}_{T}[T^{\epsilon/2}(t)< t]}{\mathbb{P}_{T}[T_0>T^{\epsilon/2}(t)]}\leq C_4(T)t^{1/2+\epsilon}e^{-t^{\epsilon/2}},$$
${\tt MT}-$almost surely, using (\ref{ptepsilon}) and the previous estimate.
We get then that almost surely, for $t$ large enough (the ``enough'' depending on $T$), 
$$\mathbb{P}_{T}[T^{\epsilon/2}(t)> n|T_0>T^{\epsilon/2}(t)]>1/2.$$
Therefore for some positive constant $C_5(T),$
$$\mathbb{P}_{T}[T_0>t]\geq C_5(T)t^{-1/2-\epsilon/2}.$$ We deduce by taking the expectation that 
$$\mathbb{P}_{\tt MT}[T_0>t]\geq C_{22} t^{-1/2-\epsilon/2},$$
for some positive and deterministic constant $C_{22}$. 
Now $$\mathbb{P}_{\tt MT}\left[\sum_{i=1}^{t^{1/2+\epsilon}}{\eta_i-\tau_i}<t\right]
\leq\left(1-\nu_0 C_{22} t^{-1/2-\epsilon/2}\right)^{t^{1/2+\epsilon}}\leq e^{-C_{23} t^{\epsilon/2}}.$$
An application of the Borel-Cantelli Lemma finishes the proof of (\ref{l1}) and (\ref{l1'}).
This finishes the proof of Lemma \ref {lemme}.

We now turn to our last part, namely the annealed central limit theorem. The proof has many parts in common with the proof in the quenched case, so we feel free to refer to the previous part.
\section{Proof of Theorem \ref{aclt}.}\label{sec7}
We recall from section \ref{sec3} the definition of the ``environment seen from the particle $T_t=\theta^v(T).$
As for the quenched case, we will first show a central limit theorem on ${\tt IMT}$ trees, then in a second part we will use the coupling to deduce the result for ${\tt MT}$ trees
\subsection{The annealed CLT on ${\tt IMT}$ trees}
We will first show the following proposition :
\begin{proposition}\label{annealedclt1}
 Suppose $N(e)\geq 1$, $q-a.s.$, (\ref{ellipticity}).
If $p=1$, $\rho'(1)<0$ and $\kappa\in(2,\infty]$,
then there is a deterministic constant $\sigma>0$ such that,
under $\mathbb{P}_{\tt IMT}$, 
the process $\{h(X_{\lfloor n t\rfloor})/\sqrt{\sigma^2 n}\}$ converges in law
to a standard Brownian motion, as $n$ goes to infinity.
\end{proposition}
\textbf{Remark :} This result is of great theoretical interest, as it is the only context where we are able to cover the whole case $\kappa>2$, we could actually make the proof of Theorem \ref{aclt} without this proposition, but as it has an interest in itself, we give the proof in the general case.
 
\noindent{\bf Proof :}
 Let, as in the quenched setting, $0<\delta<1$ and $\rho_t$ be a random variable, independent of the walk, uniformly chosen in $[t,t+t^{\delta}]$.
We recall from (\ref{defsx}) the definition of $S_x,x\in T$ and from (\ref{defeta}) the definition of $\eta$. We are going to show the following 
 \begin{lemma}\label{alemma}
 Under the assumptions of Theorem \ref{aclt}, 
\begin{equation} \label{7.1}\frac{S_{X_{\rho_t}}}{h(X_{\rho_t})} \rightarrow \eta,\end{equation} in probability.
 \end{lemma}
 We admit for the moment this lemma and finish the proof of Proposition \ref{annealedclt1}.
 We have
 $$h(X_{\rho_t})=\frac{h(X_{\rho_t})}{S_{X_{\rho_t}}} S_{X_{\rho_t}}=\eta S_{X_{\rho_t}} +\left[\frac{h(X_{\rho_t})}{S_{X_{\rho_t}}}-\eta\right]S_{X_{\rho_t}}.$$
 The first term converges to a Brownian motion with variance $\sigma$, by the same arguments as in the quenched setting, while the second one is a $o(S_{X_{\rho_t}})$. The result then follows easily, using the same arguments as in the proof of Theorem \ref{clt}.
\\
 
 We now turn to the proof of Lemma \ref{alemma}.
 The proof is quite similar to the proof of Proposition \ref{proposition}: we take some small $\epsilon>0$, then we estimate the number of visits to the points in $B_n^{\epsilon}$ during one excursion in $T^{v_i}$, and estimate the number of such excursion before time $n$. 
 We rely on the following lemma, similar to Lemma \ref{cond}
 \begin{lemma}\label{conda} Suppose that the assumptions of Theorem \ref{aclt} are true. Then for $1<\lambda<\kappa-1\wedge 2,$ and $n>0$, there exists some constant $C'_{1}$ such that
$$E_{\tt MT}\left[\sum_{x\in A_n^{\epsilon}} C_x\right]<C'_{1}n^{-(\lambda-1)}.$$
 \end{lemma}

\noindent\textbf{Proof :} the proof relies on the same ideas as the proof of Proposition \ref{cond}.
First recall that, for n large enough, 
\begin{multline*}E_{\tt MT}\left[\sum_{x\in
	{\bf A}_n^{\epsilon}}C_x\right]\leq P_{\widehat{\tt
	MT^*_n}}\left[\left|S_{v_{n}}-E_{{\widehat{\tt
	 MT^*_{n}}}}[S_{v_{n}}|\tilde{F}_{n}^*]\right|>\frac{n\epsilon}{4}\right]\\+P_{\widehat{\tt MT^*_n}}\left[\left|E_{{\widehat{\tt MT^*_n}}}[S_{v_{n}}|\tilde{F}_{n}^*]-E_{{\widehat{\tt MT^*_n}}}[S_{v_{n}}]\right|>\frac{n\epsilon}{4}\right]:=P_{1}+P_{2}.\end{multline*}
To bound $P_1$, we recall that, under the law $\widehat{\tt
	MT^*_n}$, 
$$S_{v_{n}}-E_{{\widehat{\tt MT^*_{n}}}}[S_{v_{n}}|\tilde{F}_{n}^*]=\sum_{i=0}^{n}\tilde{W}_{i}^*B_{i},$$ where $W_i$ are centered and independent random variables with bounded moments of order $\lambda+1$ and $$B_{j}=\sum_{k=0}^j\prod_{i=k+1}^j A_{i}.$$ Using Inequality 2.6.20 from page 82 of \cite{petrov1995limit}, we obtain that, for some constant $C_2$ 
$$E\left[\left|S_{v_{n}}-E_{{\widehat{\tt MT^*_{n}}}}[S_{v_{n}}|\tilde{F}_{n}^*]\right|^{\lambda}\right]<C_2\sum_{k=0}^n E[B_k^\lambda].$$
Then, using the same arguments as in the proof of Proposition \ref{cond}, we get that $E[B_k^\lambda]$ is bounded independently of $n$ and $k$, whence 
$$E\left[\left|S_{v_{n}}-E_{{\widehat{\tt MT^*_{n}}}}[S_{v_{n}}|\tilde{F}_{n}^*]\right|^{\lambda}\right]<C_3 n.$$ Using Markov's Inequality, 
there exists $C_4$ such that 
\begin{equation}\label{p1a}P_1 < \frac{C_4}{2} n^{-(\lambda-1)}.\end{equation}
On the other hand, recalling (\ref{spine}),
\begin{equation*}\left|E_{{\widehat{\tt MT^*_n}}}[S_{v_{n}}|\tilde{F}_{n}^*]-E_{{\widehat{\tt MT^*_n}}}[S_{v_{n}}]\right|<C_5+\left|\sum_{k=1}^n \tilde{A}_k D_kg(A_{n+1})(1+\rho+\rho^2+...\rho^{k-1})\right|,\end{equation*}
where $C_{5}$ is a finite constant and 
$$D_k= \sum_{j=k+1}^n \prod_{i=k+1}^{j}A_{i}g(A_{j+1}),$$ where $g$ is a bounded function. We recall that 
$$N_k:=\sum_{j=n-k}^n \tilde{A_{j}}D_{j}(1+\rho+\rho^2+...\rho^{j-1})$$
is a martingale with respect to the filtration $\mathcal{H}_k=\sigma(A_j, n-k\leq j \leq n),$ whence, using Burkholder's Inequality, 
$$E_{{\widehat{\tt MT^*_n}}}[(N_n)^{\lambda}]\leq C_6 E_{{\widehat{\tt MT^*_n}}}\left[\left(\sum_{i=0}^n (D_i)^2\right)^{\lambda/2}\right].$$
We recall that $1<\lambda<(\kappa-1)\wedge 2$, whence, by concavity, the last expression is bounded above by 
$$C_6 E_{{\widehat{\tt MT^*_n}}}\left[\sum_{i=0}^n (D_i)^{\lambda}\right]< C_7 n.$$ 
Therefore, using Markov's Inequality, we get that 
$$P_2<n^{1-\lambda}.$$
This, together with (\ref{p1a}), finishes the proof of Lemma \ref{conda}.\\

We now finish the proof of Lemma \ref{alemma}. Let us go back to {\tt IMT} trees.
We recall the definition of the sets ${\bf B}_n^{\epsilon}$:
\begin{equation}{\bf B}_n^{\epsilon}=\left\{v\in T, d(v,{\tt Ray})=n,\left|\frac{S_v^{\tt
 Ray}}{n}-\eta\right|>\epsilon\right\}.\end{equation}
We are going to prove that
$$\lim_{t\rightarrow \infty} \mathbb{P}_{T}(X_{\rho_{t}}\in \cup_{n=1}^\infty {\bf B}_m^\epsilon) =0,\, {\tt IMT}-a.s..$$ 
We introduce $\gamma>1/2$, and recall the definition of the event
$$\Gamma_{t}=\{ \exists u\leq 2t | X_{u}=v_{\lfloor t^\gamma \rfloor}\}. $$
It is easy to see, using the same arguments as in the proof of Lemma \ref{bm}, that
$$P_{\tt IMT}[\Gamma_t]\underset{t\rightarrow \infty}{\rightarrow}0.$$
Furthermore, we introduce the event
$$\Gamma'_t =\{\exists \,0\leq u \leq t, d(X_u,{\it Ray})>n^{\gamma}\};$$
then it is a direct consequence of Lemma \ref{lem3.2} that 
$$P_{\tt IMT}[\Gamma'_t]\underset{t\rightarrow \infty}{\rightarrow}0.$$
As for the quenched case, we have
\begin{eqnarray}\nonumber\mathbb{P}_{\tt IMT}(X_{\rho_{t}} \in \cup_{m=1}^\infty {\bf B}_m^\epsilon) &\leq& \mathbb{P}_{\tt IMT}(X_{\rho_{t}}\in \cup_{m=1}^{n^{\gamma}} {\bf B}_m^\epsilon;\Gamma_{t}^{c} \cap {\Gamma'}_{t}^{c})+ 
\mathbb{P}_{\tt IMT}(\Gamma_{t})+ \mathbb{P}_{\tt IMT}({\Gamma'}_{t})\\ &\leq&\frac{1}{\lfloor t^\delta \rfloor}E_{\tt IMT}\left[\mathbb{E}_{T}\left[\sum_{s=0}^{H_{v_{\lfloor t^\gamma \rfloor}}}\mathds{1}_{X_{s}\in \cup_{m=1}^{t^{\gamma}} {{{\bf B}^\epsilon}_m}}\right]\right]+ o(1)\label{gammat'},\end{eqnarray}
where $H_{v_{\lfloor t^\gamma \rfloor}}$ is the first time the walk hits $v_{\lfloor t^\gamma \rfloor}.$

We recall that $T^{(v_i)}$ the subtree constituted of the vertices $x\in T$ such that $v_i\leq x, v_{i-1}\not\leq x$, the same computations as in the proof of Lemma \ref{bm} imply 

\begin{equation}\label{pmt}\mathbb{P}_{\tt IMT}(X_{\rho_{t}}\in \cup_{m=1}^\infty {\bf B}_m^\epsilon)\leq\frac{1}{\lfloor t^\delta \rfloor}E_{\tt IMT}\left[\sum_{i=0}^{\lfloor t^\gamma \rfloor}\mathbb{E}_{T}\left[\sum_{s=0}^{H_{v_{\lfloor t^\gamma \rfloor}}}\mathds{1}_{X_{s}=v_{i}}\right] \tilde{N}_{i}\right],\end{equation}
Where $\tilde{N}_{i}$ is the $\mathbb{P}_{T}-$expectation of the number of visits to $\cup_{m=1}^{n^{\delta}} {\bf B}_m^\epsilon\cap T^{(v_{i})}$ during one excursion in $T^{(v_{i})}.$ Lemma \ref{conda}, and the method of \ref{espwetoile} imply that, under {\tt IMT} conditioned on $\{{\it Ray}, A(v_{i})\}$, $\tilde{N}_{i}$ are independent and identically distributed variables, with expectation at most equal to 
$C'_1\sum_{i=0}^{n^{\gamma}} i^{1-\lambda}$ for some $\lambda>1$. By choosing $\gamma$ close enough to $0$, we get
$E_{\tt IMT}[\tilde{N}_i|\{{\it Ray}, A(v_{i})\}]\leq C'_1 n^{1/2-\varepsilon}$ for some $\varepsilon>0$.
We recall that 
$$\mathbb{E}_{T}\left[\sum_{s=0}^{H_{v_{\lfloor t^\gamma \rfloor}}}\mathds{1}_{X_{s}=v_{i}}\right]\leq C'' \left(1+ \sum_{j=0}^{{\lfloor t^\gamma \rfloor}-1}\prod_{k=j-1}^{{\lfloor t^\gamma \rfloor}}A(v_{k})\right).$$
The latter expression has bounded expectation under {\tt IMT}, as an easy consequence of Statement \ref{big} and Lemma \ref{changement loi}.

We deduce that
$$ \mathbb{P}_{\tt IMT}(X_{\rho_{t}}\in \cup_{m=1}^\infty {\bf B}_m^\epsilon)\leq C_5 n^{\frac{1}{2}-\varepsilon +\gamma-\delta}.$$
Since $\gamma$ can be chosen as close to $1/2$ as needed, the exponent can be taken lower than $0$. The end of the proof is then completely similar to the quenched case.

 \subsection{The annealed CLT on ${\tt MT}$ trees.}
We now turn to the proof of Theorem \ref{aclt}. We use the coupling and the notations presented in section \ref{sec5}.
Our main proposition in this part will be the following:
\begin{proposition}
 Under the assumptions of Theorem \ref{aclt}, for some $\alpha<1/2$
\begin{equation}\label{alt1}
\lim_{t\rightarrow \infty} \mathbb{P}_{\tt MT}(\Delta_t\neq
\Delta_t^{\alpha})=0,
\end{equation}
and 
\begin{equation}\label{alt2}
\lim_{t\rightarrow \infty} \mathbb{P}_{\tt IMT}(\tilde{\Delta}_t\neq
\tilde{\Delta}_t^{\alpha})=0.
\end{equation}
Further, under ${\tt MT},$
\begin{equation}\label{alt3}
\limsup \frac{\Delta_t}{t}=0,
\end{equation}
and under ${\tt IMT},$
\begin{equation}\label{alt4}
\limsup \frac{\tilde{\Delta}_t}{t}=0.
\end{equation}
Finally, under {\tt IMT},
\begin{equation}\label{alt5}
\limsup \frac{{\bf B}_t}{\sqrt{t}}=0.\end{equation}
(Here $\limsup$ denotes the limit in law.)
\end{proposition}
Before proving the latter proposition, we introduce some technical estimates, whose proof will be postponed.
\begin{lemma}\label{lemmtech}
For all $\epsilon>0$ 
\begin{align}&\label{la1}\lim_{t\rightarrow\infty} \mathbb{P}_{\tt MT}\left(\sum_{i=1}^{t^{1/2 +\epsilon}}(\eta_{i}-\tau_{i})<t\right)=0,\\
&\label{la1'}\lim_{t\rightarrow\infty} \mathbb{P}_{\tt IMT}\left(\sum_{i=1}^{t^{1/2 +\epsilon}}(\tilde{\eta}_{i}-\tilde{\tau}_{i})<t\right)=0,\\
\label{la3}&\lim_{t\rightarrow\infty}\mathbb{P}_{\tt MT}(\exists k\leq I_t,
\Theta_{i-1},\Theta_{i}\in J_k,|Y_{\Theta_i}|>|Y_{\Theta_{i-1}}|)=0,\\
\label{la3'}&\lim_{t\rightarrow\infty}\mathbb{P}_{\tt IMT}(\exists k\leq I_t,
\tilde{\Theta}_{i-1},\tilde{\Theta}_{i}\in \tilde{J}_k,d(\tilde{Y}_{\tilde{\Theta}_i},Ray)>d(\tilde{Y}_{\tilde{\Theta}_{i-1}},Ray)=0,\\
&\label{la4'}
\lim_{t\rightarrow\infty} \mathbb{P}_{\tt MT}(X_s\in \cup_{k=t^{\alpha}-(\log t)^2}^{t^{\alpha}}{\bf A}_{t}^{\epsilon} \text{ for
 some } s\leq t\}=0,\\
&\label{la5'}
\lim_{t\rightarrow\infty} \mathbb{P}_{\tt IMT}(X_s\in \cup_{k=t^{\alpha}-(\log t)^2}^{t^{\alpha}}{\bf B}_{k}^{\epsilon} \text{ for
 some } s\leq t)=0.\\
 &\label{la6'} \lim_{t\rightarrow\infty} \mathbb{P}_{\tt MT}\left(W_{X_s}> t^{1/4-\varepsilon} \text{ for some } 0\leq s \leq t\right)=0\\
 &\label{la7'} \lim_{t\rightarrow\infty} \mathbb{P}_{\tt IMT}\left(W_{X_s}> t^{1/4-\varepsilon} \text{ for some } 0\leq s \leq t\right)=0
\end{align}
\end{lemma}

We now turn to the proof of (\ref{alt1}).
As a consequence of (\ref{la1}) and (\ref{la3}) that, with $\mathbb{P}_{\tt MT}$
probability approaching $1$ as $n$ goes to infinity, $${\bf t}(\Theta_{2
 t^{1/2+\epsilon}})>t,$$ whence, using Lemmas \ref{lemmtech} and \ref{lem3.2}, 
\begin{multline*}
\lim_{t\rightarrow \infty}
\mathbb{P}_{\tt MT}\left(\max_{s\in\cup_{k=1}^{I_{t}}J_{k}}|Y_{s}|\geq
t^\alpha\right)\\
\leq \limsup_{t\rightarrow\infty}\sum_{i=0}^{2
 t^{1/2+\epsilon}}\mathbb{P}_{\tt MT}\big(\exists j>i:|Y_{\Theta_j}|\geq
t^{\alpha}-(\log t)^2,\,Y_{\Theta_i}=e,\\ S_{Y_{\Theta_j}}\geq (\eta-\epsilon_1)
t^{\alpha}/2,\,
|Y_{\Theta_k}|>0, \forall i<k\leq j;W_{X_s}\leq t^{1/4-\varepsilon}\forall 0\leq s \leq t\\|S_{X_s}-|X_s||\leq
\epsilon t^{1/4-\epsilon'}|X_s|, \forall s\leq t\big):=\limsup_{t\rightarrow\infty}\sum_{i=1}^{2
 t^{1/2+\epsilon}}\mathbb{P}_{i,t}\,;\\
\end{multline*}
where $\epsilon,\epsilon_1$ are positive numbers that can be
chosen arbitrarily small.
We recall that the process $\{N_s\}=\{S_{X_{\theta_{i+s\wedge K_t}}}\}$ is a supermartingale.
and that there exists a previsible and non-decreasing process $A_s$ such that 
$N_s+A_s$ is a martingale.
Furthermore, on the event $\{W_{X_s}\leq t^{1/4-\varepsilon}\forall 0\leq s \leq t\},$ the increments of this martingale are bounded above by $t^{1/2-\varepsilon}.$ Azuma's Inequality implies the result, as in the quenched case.\\

The proof of (\ref{alt2}) is similar and omitted.\\

We recall that in the proofs of (\ref{lt3}),(\ref{lt4}) and (\ref{lt5}) we only used the assumption $\kappa>5$, therefore the proof of (\ref{alt3}),(\ref{alt4}) and (\ref{alt5}) are direct consequence, by dominated convergence.\\

We now turn to the proof of Lemma \ref{lemmtech}. The proofs of (\ref{la1}), (\ref{la1'}), (\ref{la3}), (\ref{la3'}) and (\ref{la4'}) follow directly from equations (\ref{l1}), (\ref{l1'}), (\ref{l3}), (\ref{l3'}) and (\ref{l4}), whose proofs did not use any assumption other than $\kappa>5$, by dominated convergence.\\

To prove (\ref{la5'}), note that, similarly to the proof of \ref{alemma},
\begin{multline*}\mathbb{P}_{\tt IMT}(X_s\in \cup_{k=t^{\alpha}-(\log t)^2}^{t^{\alpha}}{\bf B}_{k}^{\epsilon} \text{ for
 some } s\leq t)\\=E_{\tt IMT}\left[\sum_{i=0}^{\lfloor t^\gamma \rfloor}\mathbb{E}_{T}\left[\sum_{s=0}^{H_{v_{\lfloor t^\gamma \rfloor}}}\mathds{1}_{X_{s}=v_{i}}\right] N'_{i}\right],\end{multline*}
where $N'_{i}$ is the $\mathbb{P}_{T}-$expectation of the number of visits to $\cup_{k=t^{\alpha}-(\log t)^2}^{t^{\alpha}}{\bf B}_{k}^{\epsilon}\cap T^{(v_{i})}$ during one excursion in $T^{(v_{i})}.$ Lemma \ref{cond} and the method of Lemma \ref{espwetoile} imply that, under {\tt IMT} conditioned on $\{{\it Ray}, A(v_{i})\}$, $N'_{i}$ are independent and identically distributed variables, up to a bounded constant, with expection at most equal to 
$C'(\log t)^2 t^{-\alpha(\lambda-1)}$ for some $\lambda>2$.
We also recall that
$$\mathbb{E}_{T}\left[\sum_{s=0}^{H_{v_{\lfloor t^\gamma \rfloor}}}\mathds{1}_{X_{s}=v_{i}}\right]\leq C'' \left(1+ \sum_{j=0}^{{\lfloor t^\gamma \rfloor}-1}\prod_{k=j-1}^{{\lfloor t^\gamma \rfloor}}A(v_{k})\right).$$
has bounded expectation under {\tt IMT}, as an easy consequence of Statement \ref{big} and Lemma \ref{changement loi}. By choosing $\gamma$ close enough to $0$ and $\alpha$ close to $1$, we get the result.

The proofs of (\ref{la6'}) and (\ref{la7'}) are easily deduced from the proofs of (\ref{l2}) and (\ref{l6}), the only difference being that we do not need to apply the Borel-Cantelli Lemma. 

{\bf Acknowledgement:} We would like to warmly thank Y. Hu for his help and encouragement.

\bibliographystyle{plain}
\bibliography{biblio}

\begin{thebibliography}{10}

\bibitem{aidekon2008trw}
E.~Aidekon.
\newblock {Transient random walks in random environment on a Galton--Watson
  tree}.
\newblock {\em Probability Theory and Related Fields}, 142(3):525--559, 2008.

\bibitem{bercu2008exponential}
B.~Bercu and A.~Touati.
\newblock {Exponential inequalities for self-normalized martingales with
  applications}.
\newblock {\em Annals of Applied Probability}, 18(5):1848--1869, 2008.

\bibitem{Biggins:1977}
J.D. Biggins.
\newblock {Martingale convergence in the branching random walk}.
\newblock {\em Journal of Applied Probability}, pages 25--37, 1977.

\bibitem{biggins-kyprianou}
J.D. Biggins and A.E. Kyprianou.
\newblock {Seneta-Heyde norming in the branching random walk}.
\newblock {\em Annals of Probability}, 25(1):337--360, 1997.

\bibitem{Billingsley:1999}
P.~Billingsley.
\newblock {\em Convergence of Probability measures}.
\newblock Wiley, 1999.

\bibitem{chernov1967rmc}
A.A. Chernov.
\newblock {Replication of a multicomponent chain, by the lightning mechanism}.
\newblock {\em Biophysics}, 12(2):336--341, 1967.

\bibitem{gyzalmost}
G.~Faraud, Y.~Hu, and Z.~Shi.
\newblock {An almost sure convergence for stochastichally biased random walk on
  a Galton-Watson tree }.
\newblock {\em http://arxiv.org/abs/1003.5505}, 2010.

\bibitem{hu2007sbr}
Y.~Hu and Z.~Shi.
\newblock {A subdiffusive behaviour of recurrent random walk in random
  environment on a regular tree}.
\newblock {\em Probability theory and related fields}, 138(3):521--549, 2007.

\bibitem{hu2007smr}
Y.~Hu and Z.~Shi.
\newblock {Slow movement of random walk in random environment on a regular
  tree}.
\newblock {\em Annals of Probability}, 35(5):1978--1997, 2007.

\bibitem{yzpolymer}
Y.~Hu and Z.~Shi.
\newblock {Minimal position and critical martingale convergence in branching
  random walks, and directed polymers on disordered trees}.
\newblock {\em Annals of Probability}, 37(2):742--789, 2009.

\bibitem{Kemeny:1976}
J.G Kemeny, A.W. Knapp, and J.L. Snell.
\newblock {\em Denumerable Markov Chains}.
\newblock 2nd ed. Springer, 1976.

\bibitem{kipnis1986central}
C.~Kipnis and S.R.S. Varadhan.
\newblock {Central limit theorem for additive functionals of reversible Markov
  processes and applications to simple exclusions}.
\newblock {\em Communications in Mathematical Physics}, 104(1):1--19, 1986.

\bibitem{Liu:2000}
Q.~Liu.
\newblock {On generalized multiplicative cascades}.
\newblock {\em Stochastic processes and their applications}, 86(2):263--286,
  2000.

\bibitem{Liu:2001}
Q.~Liu.
\newblock {Asymptotic properties and absolute continuity of laws stable by
  random weighted mean}.
\newblock {\em Stochastic processes and their applications}, 95(1):83--107,
  2001.

\bibitem{Lyons:1989}
R.~Lyons.
\newblock The \text{I}sing model and percolation on trees and tree-like graphs.
\newblock {\em Communications in Mathematical Physics}, 125:337--353, 1989.

\bibitem{Lyons:1992}
R.~Lyons and R.~Pemantle.
\newblock Random walks in a random environment and first-passage percolation on
  trees.
\newblock {\em Annals of Probability}, 20:125--136, 1992.

\bibitem{Lyons:2005}
R.~Lyons and Y.~Peres.
\newblock {\em Probability on trees and networks}.
\newblock 2005.

\bibitem{Mandelbrot:1974}
B.~Mandelbrot.
\newblock {Multiplications aleatoires iterees et distributions invariantes par
  moyenne ponderee aleatoire: quelques extensions}.
\newblock {\em Comptes rendus de l'acad\'emie des Sciences}, 278:355--358,
  1974.

\bibitem{menshikov-petritis}
M.V. Menshikov and D.~Petritis.
\newblock {On random walks in random environment on trees and their
  relationship with multiplicative chaos}.
\newblock {\em Mathematics and computer science II (Versailles, 2002)}, pages
  415--422, 2002.

\bibitem{Neveu:1986}
J.~Neveu.
\newblock {Arbres et processus de Galton-Watson}.
\newblock {\em Annales de l'Institut H. Poincare}, 22(2):199--207, 1986.

\bibitem{peres-zeitouni}
Y.~Peres and O.~Zeitouni.
\newblock {A central limit theorem for biased random walks on Galton--Watson
  trees}.
\newblock {\em Probability Theory and Related Fields}, 140(3):595--629, 2008.

\bibitem{Petrov:1975}
V.V. Petrov.
\newblock {\em Sums of Independent Random Variables}.
\newblock Springer-Verlag, 1975.

\bibitem{petrov1995limit}
V.V. Petrov.
\newblock {\em {Limit theorems of probability theory: sequences of independent
  random variables}}.
\newblock Oxford University Press, USA, 1995.

\bibitem{Zeitouni:2003}
O.~Zeitouni.
\newblock {Lecture notes on random walks in random environment. Ecole d'\'et\'e
  de probabilit\'es de Saint-Flour 2001}.
\newblock {\em Lecture Notes in Mathematics}, 1837:189--312, 2003.

\end{thebibliography}
\end{document}